\documentclass[red,11pt,a4paper]{article}
\usepackage{cite}

\usepackage{amsmath}
\usepackage{amscd}
\usepackage{amssymb}
\usepackage{graphicx}
\usepackage{verbatim}

\usepackage{latexsym,empheq}
\usepackage{color}
\usepackage{floatrow}
\newfloatcommand{capbtabbox}{table}[][\FBwidth]

\usepackage{blindtext}

\newcommand{\Grad}{\nabla}

\newcommand{\lr}[1]{\left( #1 \right)}

\newcommand{\eq}[1]{\begin{equation}
\begin{split}
#1
\end{split}
\end{equation}}
\newcommand{\eqh}[1]{\begin{equation*}
\begin{split}
#1
\end{split}
\end{equation*}}
\newcommand{\ep}{\varepsilon}

\newtheorem{thm}{Theorem}
\newtheorem{lemma}[thm]{Lemma}
\newtheorem{hypot}{Hypothesis}
\newtheorem{prop}[thm]{Proposition}

\title{Modelling pattern formation through differential repulsion}
\author{J. Barr\'e$^{1,2}$, P. Degond$^3$, D. Peurichard$^{4}$, E. Zatorska$^5$}
\date{}

\begin{document}
%% create title page and toc
\maketitle
%    \begin{center}
%    \small
%1-Faculty of Mathematics, University of Vienna, Oskar-Morgenstern Platz 1, 1090 Vienna, Austria\\
%e-mail: Diane.Peurichard@univie.ac.at\\
%\end{center}

\centerline{1. Institut Denis Poisson, Universit\'e d'Orl\'eans, CNRS, Universit\'e de Tours.}
\centerline{B.P. 6759, 45067 Orl\'eans cedex 2, France.}
\bigskip
\centerline{2. Institut Universitaire de France, Paris, France.}
 \bigskip
\centerline{3. Department of Mathematics, Imperial College London, }
\centerline{London SW7 2AZ, United Kingdom.}
\bigskip
% \centerline{4. Neural Development Laboratory, and Computational Cell }
% \centerline{and Molecular Biology Laboratory, The Francis Crick Institute, 1 Midland Road, London NW1 1AT, UK}
% \bigskip
\centerline{4. INRIA Team Mamba, INRIA Paris, 2 rue Simone Iff, CS 42112, 75589 Paris, France}
\centerline{Universit\'e Pierre et Marie Curie-Paris 6, UMR 7598 LJLL, BC187}
\centerline{ 4, Place de Jussieu, F-75252 Paris Cedex 5, France.}
\bigskip
\centerline{5. Department of Mathematics, University College London}
\centerline{Gower Street London WC1E 6BT, UK, United Kingdom.}
\bigskip

\begin{abstract}
\noindent Motivated by experiments on cell segregation, we present a two-species model of interacting particles, aiming at a quantitative description of this phenomenon. Under precise scaling hypothesis, we derive from the microscopic model a macroscopic one and we analyze it. In particular, we determine the range of parameters for which segregation is expected. We compare our analytical results and numerical simulations of the macroscopic model to direct simulations of the particles, and comment on possible links with experiments. 
\end{abstract}

\section{Introduction}

The organisation of biological tissues during development is accompanied by the formation of sharp borders between distinct cell populations. During the morphogenesis of numerous tissues/organs, cells of the same type regroup into regions, creating niches with specific identities that drive the differentiation of particular cell types. This spatial organization is ensured via cell-cell signalling leading specific cells/tissues to form at the appropriate location. The maintenance of this cell segregation is key in adult tissue homeostatis, and its disruption can lead tumor cells to spread and form metastasis. This segregation is challenged during tissue growth and morphogenesis due to the high mobility of many cells that can lead to intermingling. Therefore, understanding the mechanisms involved in the generation and maintenance of cell segregation is of tremendous importance in tissue morphogenesis, homeostasis, and in the development of various invasive diseases such as tumors.

Numerous experiments have been conducted to identify the mechanisms of cell segregation. Experiments show that mixing cells from different tissues in vitro leads to their segregation, with initially fuzzy borders that sharpen in time \cite{Taylor2017}. This has been observed for many systems for instance in the development of the wing imaginal disc in Drosophilia \cite{Hsia2017}, in the developing nervous system \cite{Douarin86} etc. So far, three types of mechanisms have been identified to have a role in segregation and border formation, namely (i) differential adhesion, (ii) contact cell repulsion and (iii) cortical tension. Indeed, the segregation of cells derived from different tissues in vitro was initially suggested to occur through a combination of directed cell migration and the selective adhesion of cells of the same type \cite{Steinberg2007}. However, other studies have shown that the contact inhibition of cell migration induced by contact repulsion of cells by Eph receptor and ephrin signaling, and finally the induction of cortical tension by actomyosin contraction are important mechanisms that can restrict intermingling between cell populations. However, it remains unclear whether individually each of these mechanisms account for cell segregation or border sharpening, and to what extent the interplay of these different mechanisms is required to achieve cell segregation. 

Several modelling efforts have been done to identify the main mechanisms involved in cell segregation. The mathematical models for cell-cell interactions are usually agent-based models, where each cell undergoes a random walk exclusion process and interacts with its neighbours. For instance the cellular Potts model for segregation between two cell populations \cite{Glazier1993} predicts cell rearrangements in epithelia based on the minimization of a free energy and has been widely used to explore the rate of cell sorting due to differential adhesion, but does not include cell migration as a mechanism. The varying adhesive and repulsive forces between different cell populations, which can result from Eph/ephrin interactions, have also been modelled by representing cells as spheres which can attract or repel each other, giving rise to empirically observed cell sorting pattern \cite{Taylor2017}. In \cite{Aharon2014}, the authors  develop a mathematical model for Eph/ephrin regulated cell-cell segregation and tissue boundary formation, which features independent random cell motion, Eph/ephrin-dependent attraction/repulsion interactions between neighbouring cells, and cell division. The authors in \cite{Aharon2014} show that the dynamics of Eph/ephrin-mediated cell cluster formation and cell segregation can be captured with these mechanisms. However, this model do not use parameters from measurements of cell behaviour or examine whether cell repulsion is sufficient for border sharpening. Another approach, which simulates cell adhesion, de-adhesion and migration in greater detail \cite{Taylor2011,Taylor2012}, was used to model the time course of segregation of cells differing in cadherin expression. The results from this model were acurate for describing cell segregation mediated by differential expression of cadherins but less accurate when simulating the significantly faster rate of Eph-ephrin mediated cell segregation. However when modified to account for repulsive behaviors \cite{Taylor2017}, the model correctly reproduced the experiments and showed that heterotypic repulsion can account for cell segregation and border sharpening, and is more efficient than decreased heterotypic adhesion.

All these results suggest that cell segregation and border sharpening is the result of a complex interplay between homotypic/heterotypic cell adhesion, de-adhesion and repulsion. But how the balance of these phenomena is precisely linked to the existence/size of the segregated zones remains unclear to this day. In this paper, we aim to provide a mathematical framework which enables to quantitatively link the segregation and border sharpening ability of the tissue to these cell-cell interaction phenomena of interest. As agent-based models do not enable precise mathematical analysis of their solutions due to the lack of theoretical results, we turn towards a continuous -macroscopic- model for which the theoretical study gives precise criteria for phase transitions as functions of key model parameters. As a drawback contrary to microscopic models, macroscopic models lose the information at the individual level.  In order to overcome this weakness, we aim to derive, as rigorously as possible,  the macroscopic model from an agent-based formulation to ensure the good correspondence between the two formulations as it was done in \cite{Barre2016,Degond_Motsch_M3AS08,Degond_Peurichard_2016}. 

The starting point is an individual-based model inspired from \cite{Barre2016,Barre2017} and which bears similarities with the approach \cite{Taylor2017}. We consider two families of cells, each cell being modelled as a point particle which interacts with its close neighbors via local cross-links. The links are modeled by springs that are randomly created and destructed. This enables us to model cell-cell attraction and repulsion, with different spring strengths according to the type of link (intraspecies or interspecies). We let the particles move randomly in space to model random motion of cells. In the mean field limit, assuming large numbers of particles and links as well as propagation of chaos, the corresponding kinetic system consists of two equations for the individual particle distribution functions and two equations for the link densities. In the large-scale limit and in the regime where the link creation/destruction frequency is very large, it was shown in \cite{Degond_Peurichard_2016,Barre2017,Barre2016} that the link density distributions become local functions of the particle distributions. The latter evolve through aggregation diffusion equations. Similar macroscopic model, but with nonlinlinear porous medium type of diffusion has been recently considered analytically and numerically in \cite{CaHuSch}. 
Although the derivation of the macroscopic model in the hydrodynamic limit follows closely the steps of \cite{Barre2017}, the originality of this work lies in the presence of two coupled families of cells which introduce a new level of complexity and make the stability analysis more involved than in previous works. Inspired from the results of \cite{Taylor2017}, we mainly consider repulsive springs and aim to quantify the influence of heterotypic/homotypic repulsion on cell segregation and border sharpening. By addressing the stability of a homogeneous distribution of particles for Hookean repulsive potentials, we obtain a precise condition for the phase transition, which links the system segregation ability to the model parameters and give further insight into the cell segregation processes.

\begin{figure}[h]
\includegraphics[scale=0.45]{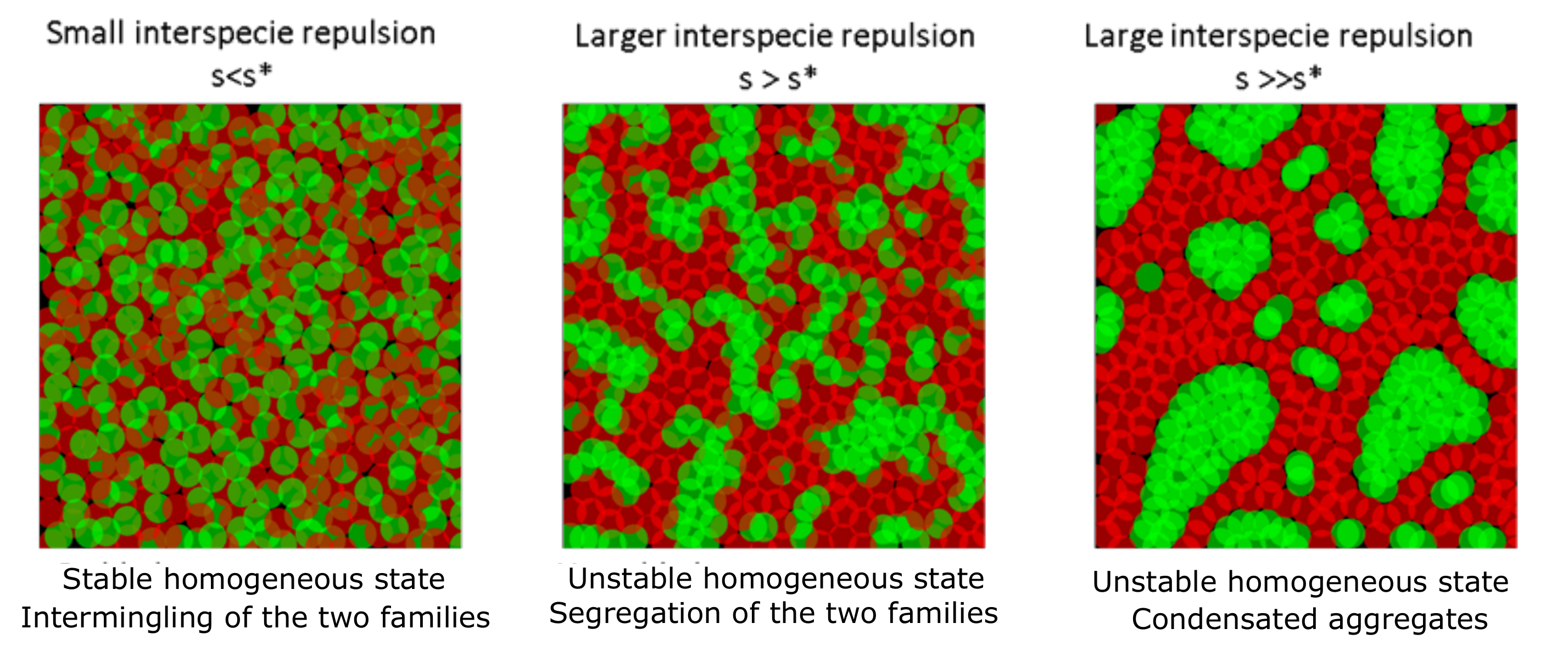}
\caption{Scheme of the predictions of the linear stability analysis by acting on the interspecies repulsion force.  \label{scheme}}
\end{figure}

Our study shows that in a system composed of two-species repelling each other, the interspecies forces must be large enough to compensate both for the diffusion and for the intra-species repulsion, which both tend to homogeneize the system. Aggregation will therefore be ensured if and only if interspecies repulsion wins over diffusion and intraspecies repulsion.  In the case where attractive interactions are considered, we have noted that a necessary condition for the aggregation of the species (or equivalently instability of the homogeneous steady-state) is that the interspecies forces are of the same sign. To observe aggregates, the two families must therefore either repulse or attract each other, but must have the same effect on each other. On the contrary, if one family is attracted by the other and the other repulses it, we will always observe a homogeneous distribution at equilibrium (intermingling of the two families). A third remark concerns the size of the clusters when aggregation occurs. As the interspecies repulsion force increases, the particles of a given family aggregate more together, leading to a decrease of the size of the local aggregates of the compressed family. These conclusions were confirmed by simulations performed for both the microscopic and macroscopic models and are pictured in Fig. \ref{scheme}. Numerical simulations show that both the micro- and macro- models are in excellent agreement with the predictions of the stability analysis performed on the continuous model. The quantitative agreement obtained between both models show that the macroscopic model is a good approximation of the microscopic model as the number of individuals goes to infinity, provided the interspecies repulsion forces are not too large. Indeed for large interspecies repulsion forces, we find some structural discrepancies between the two models, where the microscopic dynamics seems to favor the formation of rounder aggregates compared to the elongated structures obtained with the macroscopic model. We find that the microscopic dynamics is comprised of two time phases: the first phase consists of a fast segregation between the two families and is followed later by a reorganisation of the clusters which get rounder at large times. The macroscopic dynamics does not seem to contain the second phase (restructuring), suggesting that this phenomenon can be due to finite size effects. We postulate that these discrepancies come from the microscopic noise due to thermal fluctuations, which gives rise to instabilities and allows the agent-based system to reach new states which are not available in the deterministic description, or produce spatial correlations which in turn dominate the macroscopic system behavior. Several works have reported these phenomena \cite{Nesic_PhysRev2014,Kessler_Nature98,Armero_PhysRev98,Cardy_JSP98}, particularly at onset for transitions from metastable or unstable phases, in which microscopic noise can be amplified to macroscopic time and length scales. The exploration of these effects will be the subject of future works. Finally, we remark that the segregation process is seen in the numerical simulations to be efficient even close to the instability threshold: as soon as the homogeneous state becomes unstable, the system evolves towards a well segregated configuration. This suggests the presence of a subcritical bifurcation \cite{ChayesPanferov10, Barre2016, CaGvPaSch}, the study of which we also leave for future work.

%\red{I have not put more remarks on the discrepancies micro/macro, possibly related to metastable states, etc. I found the current explanations mentioning only "finite size effects" to be enough. (J)}

The paper is organized as follows. In Section \ref{sec:2}, we first give the main ingredients of the microscopic model, and then we sketch the main steps of the derivation of the macroscopic dynamics. We present two approaches that involve taking limit of large number of individuals and large scale/fast network remodelling limit of the microscopic model in the different orders; the details  of one of them are moved to Appendixes \ref{AppendixA} and \ref{AppendixB}. Section \ref{sec:3} is devoted to the stability analysis around the homogeneous steady-states of the macroscopic model: in Section \ref{sec:3.1}, we give the stability results in the whole space, Section \ref{sec:3.2} explores the case of periodic boundary conditions and Section \ref{sec:separatedphases} is devoted to the case of phase separated initial conditions. Finally, Section \ref{sec:numerical} presents the numerical results, performed on the microscopic and macroscopic models  in different regime of parameters, with a particular focus on the qualitative and quantitative comparison between the two models.

\section{Mathematical modelling}\label{sec:2}
\subsection{Microscopic model}
\label{sec:IBM}

The model features two families of particles referred to as type $A$ and type $B$. Each particle can link/unlink with  neighbors located in a ball of radius $R$ from its center. Each particle can link with a neighbor of its own family as well as with a neighbor of the other family, with no restriction on the number of links per particle and with the same detection radius $R$ no matter the type of  link (intra- or inter- species). In order to model tissue plasticity, the links are not permanent but supposed to be  created and suppressed via random processes. In this way, the model allows for constant remodelling of the link network. Each link between two particles generates a spring-like interaction potential, which depends on the link type (intra- or inter- species link). Finally, particle positions are subject to random positional noise to model the movements of the tissue.

In this paper, we restrict ourselves to a two-dimensional model. We consider a set of $N_A$ particles of type $A$ and $N_B$ particles of type $B$ described by their centers $(X^A_i,X^B_\ell) \in \mathbb{R}^2 \times \mathbb{R}^2$, $i\in [1,N_A], \ell \in [1,N_B]$ respectively. The link creation and suppression are supposed to follow Poisson processes of frequencies $\nu_{c,N,\ep}^{AA}, \nu_{c,N,\ep}^{AB},\nu_{c,N,\ep}^{BB}$ and $\nu_{d,\ep}^{AA}, \nu_{d,\ep}^{AB},\nu_{d,\ep}^{BB}$, where the subscripts $c$ and $d$ refer to 'creation' and 'deletion' respectively, and the superscripts $AA,BB$ and $AB$ denote intraspecies links (AA,  BB) and interspecies links (AB); $\ep$ is a scaling parameter, and the subscripts $\ep$, and $N$ for the $\nu_c$'s make  the dependency of these rates on $N_A$, $N_B$ and $\ep$ explicit, as will be explained below. We suppose that the intraspecies links generate pairwise symmetric potentials $\Phi^{AA}(X^A_i,X^A_j)$ and $\Phi^{BB}(X^B_\ell,X^B_m)$, not necessarily equal, and that the interspecies links generate non-symmetric potentials $\Phi^{AB}(X^A_i,X^B_\ell) \neq \Phi^{BA}(X^B_\ell,X^A_i)$, modelling the fact that the two particle families act differently on each other. For the moment we do not specify  interaction potentials, trying to
keep the derivation at maximal level of generality. Note that $\Phi^{AB}$ refers to the action a type B particle exerts on a type $A$ particle while $\Phi^{BA}$ is the action a type $A$ particle exerts on a particle of type $B$. We define the total "energy" $W^A$ of the $A$-particles  as the sum over all pairwise link potentials acting on particles $A$, and $W^B$ is the sum over all pairwise link potentials acting on particles~$B$:
% \blue{\emph{There is a slight abuse of notation here, $i(k_1)$ and $i(k_3)$ do not represent the same function... I did not change it because I was afraid a fix would make things worse.}}
\begin{align}
W^A(X^A,X^B) = &\sum_{k_1=1}^{K_{AA}} \Phi^{AA}(X^A_{i(k_1)},X^A_{j(k_1)})+ \sum_{k_3=1}^{K_{AB}} \Phi^{AB}(X^A_{i(k_3)},X^B_{\ell(k_3)}) \label{PotA}\\
W^B(X^A,X^B) = & \sum_{k_2=1}^{K_{BB}} \Phi^{BB}(X^B_{\ell(k_2)},X^B_{m(k_2)}) + \sum_{k_3=1}^{K_{AB}} \Phi^{BA}(X^B_{\ell(k_3)},X^A_{i(k_3)}),\label{PotB}
\end{align}
\noindent where $K^{AA},K^{BB},K^{AB}$ denote the (time-dependent) total number of links between particles of type $A$, particles of type $B$, and interspecies links respectively. In the formulas above  $(i(k_1),j(k_1))$ denote the indices of particles of type $A$ connected by the intraspecies link $k_1$, $(\ell(k_2),m(k_2))$ the indices of particles of type $B$ connected by link $k_2$. By a slight abuse of notation we denote by $(i(k_3),\ell(k_3))$ the indices of particles of type A connected to particles of type B by link $k_3$. 

Particle motion during a time interval between two linking/unlinking events is supposed to occur in the so-called overdamped regime. The resulting equations contain a drift term 
in the steepest descent direction of the "energies" $W^A$ and $W^B$ and a noise term:
\begin{empheq}[left=\empheqlbrace]{align} 
{dX^A_i} = - \mu \, \nabla_{X^A_i} W^A(X^A,X^B)dt + \sqrt{2 D^A} dB_i, \quad \forall i \in \{1, \ldots, N_A \}, \label{IBMA}\\
{dX^B_i}= - \mu \, \nabla_{X^B_\ell} W^B(X^A,X^B)dt + \sqrt{2 D^B} dB_\ell, \quad \forall \ell \in \{1, \ldots, N_B \} \label{IBMB},
\end{empheq}
where $\mu>0$ is the mobility coefficient considered to be given, and $B_i$ is a 2-dimensional Brownian motion $B_i = (B_i^1,B_i^2)$ of intensity $D^A>0$ for species $A$ and $D^B>0$ for species $B$. Inserting \eqref{PotA}-\eqref{PotB} into \eqref{IBMA}-\eqref{IBMB}, we obtain
\begin{align}
dX_i^A = &-\mu dt \bigg[ \sum_{k_1=1}^{K_{AA}} \big(\nabla_{x_1} \Phi^{AA}(X^A_{i(k_1)},X^A_{j(k_1)})\delta_{i(k_1)}(i) + \nabla_{x_2} \Phi^{AA}(X^A_{i(k_1)},X^A_{j(k_1)})\delta_{j(k_1)}(i)\big) \nonumber\\
& + \sum_{k_3=1}^{K_{AB}} \nabla_{x_1} \Phi^{AB}(X^A_{i(k_3)},X^B_{\ell(k_3)})\delta_{i(k_3)}(i)\bigg] + \sqrt{2D^A} dB_i\label{MicroA}\\
dX_\ell^B = &-\mu dt \bigg[ \sum_{k_2=1}^{K_{BB}} \big( \nabla_{x_1}\Phi^{BB}(X^B_{\ell(k_2)},X^B_{m(k_2)})\delta_{\ell(k_2)}(\ell) + \nabla_{x_2}\Phi^{BB}(X^B_{\ell(k_2)},X^B_{m(k_2)})\big)\delta_{m(k_2)}(\ell) \nonumber\\
& + \sum_{k_3=1}^{K_{AB}} \nabla_{x_1}\Phi^{BA}(X^B_{\ell(k_3)},X^A_{i(k_3)})\delta_{\ell(k_3)}(\ell)\bigg] + \sqrt{2D^B} dB_\ell,\label{MicroB}
\end{align}
where $\delta_{i}(j)$ stands for the Kronecker delta.\\
\noindent This model bears similarities with the works of \cite{Taylor2011,Taylor2012,Taylor2017}. The main difference lies in the fact that cells are modelled as individual spheres here, while in \cite{Taylor2017} each cell is supposed to be composed of a set of several spheres maintained in a ring of a given radius. Therefore we do not take into account the role of cell deformation in this paper. However, this simpler modelling choice enables us to reduce the complexity of the system and to derive a macroscopic model as performed in the next section. The influence of this modelling simplification on the result will be the subject of future works. In the next section, we present the macroscopic model obtained in the limit of a large number of individuals and in the fast linking-unlinking process.

%%%%%%%%%%%%%%%%%%%%%%%%%%%%%%%%%%%%%%%%%%%%%%%%%%%%%%%%%%%%
%%%%%%%%%%%%%%%%%%%%%%%%%%%%%%%%%%%%%%%%%%%%%%%%%%%%%%%%%%%%
%%%%%%%%%%%%%%%%%%%%%%%%%%%%%%%%%%%%%%%%%%%%%%%%%%%%%%%%%%%%
%%%%%%%%%%%%%%%%%%%%IBM-CINETIQUE%%%%%%%%%%%%%%%%%%%%%%%
%%%%%%%%%%%%%%%%%%%%%%%%%%%%%%%%%%%%%%%%%%%%%%%%%%%%%%%%%%%%
%%%%%%%%%%%%%%%%%%%%%%%%%%%%%%%%%%%%%%%%%%%%%%%%%%%%%%%%%%%%
%%%%%%%%%%%%%%%%%%%%%%%%%%%%%%%%%%%%%%%%%%%%%%%%%%%%%%%%%%%%

\subsection[Macroscopic model]{Macroscopic model}\label{sec:macro}
The derivation of a macroscopic model from the microscopic model defined by Eqs. \eqref{MicroA}-\eqref{MicroB} requires two limits: 
(i) limit of large number of individuals and large number of links, denoted $N,K \to \infty$ and (ii) a large scale or fast network remodelling limit, denoted 
$\varepsilon \to 0$. 

Performing the  limit (i)  first yields, after a mean-field assumption, a kinetic system, from which the macroscopic dynamics can be derived  in the $\varepsilon \to 0$ limit. This approach will be referred to as the {\bf Approach I}.

Performing the limit (ii)  first yields, after averaging over the network configurations, an effective dynamics for the particles, from which the same macroscopic 
dynamics can be derived, in the $N,K \to \infty$ limit. This approach will be referred to as the {\bf Approach II}.

This structure is summarized in the diagram: 
\begin{equation}
\begin{CD} \label{eq:CD}
{\rm Microscopic}~\left\{(X_i^A,X_j^B)\right\}  @> N,K \to \infty >> {\rm Kinetic}~\left\{(f^S,g^{ST},h^{ST})_{S,T=A,B}\right\} \\
@VV \varepsilon \to 0V  @ VV \varepsilon \to 0 V \\
{\rm Averaged~microscopic}~\left\{(\tilde{X}_i^A,\tilde{X}_j^B)\right\}           @> N,K \to \infty >> {\rm Macroscopic}~\left\{(f^A,f^{B})\right\}
\end{CD} 
\end{equation}
The meaning of the $f^S,g^{ST}, h^{ST}, \tilde{X}^S$ is explained below. 
Our final goal is a macroscopic model describing the evolution in time of the particle distributions $f^{A}(x,t)$ and $f^{B}(x,t)$ of the type-$A$ particles and type-$B$ particles respectively.  Let us roughly discuss the two possible approaches mentioned above.

\subsubsection*{Sketch of  Approach I}
For finite $N_A,N_B$ we define: 
\begin{equation}\label{fAB}
f^A_N(x,t) = \frac{1}{N_A} \sum_{i=1}^{N_A} \delta_{X^A_i(t)}(x),\qquad 
f^B_N(x,t) = \frac{1}{N_B} \sum_{\ell=1}^{N_B} \delta_{X^B_\ell(t)}(x),
\end{equation}
\noindent
where $\delta_{X^S_i(t)}(x)$ denotes the Dirac delta located at $X^S_i(t)$ for $S$ being either $A$ or $B$. In the large $N_S$ limit it gives the probability to find a particle of type $S$ at point $x$ at time $t$. 

To write the kinetic model, we need to define the (symmetric) empirical measures $g^{SS}_N(x_1,x_2,t)$ of the intraspecies links ($S$ being either A or B) by:
\begin{equation}\label{gAA}
\begin{split}
g^{AA}_N(x_1,x_2,t) =& {\frac{1}{2N_A}} \sum_{k_1=1}^{K^{AA}} \delta_{X^A_{i(k_1)}, X^A_{j(k_1)}}(x_1,x_2) +  \delta_{X^A_{j(k_1)}, X^A_{i(k_1)}}(x_1,x_2),\\
g^{BB}_N(x_1,x_2,t) =& {\frac{1}{2N_{B}}} \sum_{k_2=1}^{K^{BB}} \delta_{X^B_{\ell(k_2)}, X^B_{m(k_2)}}(x_1,x_2) +  \delta_{X^B_{m(k_2)}, X^B_{\ell(k_2)}}(x_1,x_2),
 \end{split}
\end{equation} 
\noindent 
with a similar definition of the Dirac deltas. Such a $g_N^{{SS}}(x_1,x_2,t)$ gives in the large $N_S$ limit the density of links connecting a particle of a given type and located within a volume $dx_1$ about $x_1$ with a particle of the same type located within a volume $dx_2 $ about $x_2$, normalized by $N_S$; note the integral of $g^{SS}_N$ is not $1$. As will become clear below, we will be interested in a regime where $K_{SS}$ and $N_S$ have the same order of magnitude; hence the chosen normalization ensures that $g^{SS}_N$ is of order $1$. Now, we define a non-symmetric empirical measures for the interspecies links $g^{AB}_N(x_1,x_2,t)$:
\begin{equation}\label{gAB}
g^{AB}_N(x_1,x_2,t) ={\frac{1}{N_{A}}} \sum_{k_3=1}^{K^{AB}} \delta_{X^A_{i(k_3)}, X^B_{\ell(k_3)}}(x_1,x_2), 
\end{equation} 
and analogously 
\begin{equation}\label{gBA}
g^{BA}_N(x_1,x_2,t) = {\frac{1}{N_{B}}} \sum_{k_3=1}^{K^{AB}} \delta_{X^B_{\ell(k_3)}, X^A_{i(k_3)}}(x_1,x_2) = \frac{N_A}{N_B} g^{AB}_N(x_2,x_1,t). 
\end{equation} 
\noindent Note that $g^{AB}_N(x_1,x_2,t)$ gives in the large $N_A$ limit the density of interspecies links between a particle of type A located within a volume $dx_1$ about $x_1$ and a particle of type B located within a volume $dx_2$ about $x_2$, normalized by $N_A$. 

Finally, the derivation of the kinetic model also involves the two-particle distribution functions defined by:
\begin{align}
h^{{AA}}_N(x_1,x_2,t) &=  \frac{1}{2N_A(N_A-1)} \sum_{i=1}^{N_A}\sum_{j=1, j\neq i}^{N_A} \big(\delta_{X^A_i(t),X^A_j(t)}(x_1,x_2) + \delta_{X^A_j(t),X^A_i(t)}(x_1,x_2)\big) \label{TPDFAA}\\
h^{BB}_N(x_1,x_2,t) &=  \frac{1}{2N_B(N_B-1)} \sum_{\ell=1}^{N_B}\sum_{m=1, m\neq \ell}^{N_B} \big(\delta_{X^B_\ell(t),X^B_m(t)}(x_1,x_2) + \delta_{X^B_m(t),X^B_\ell(t)}(x_1,x_2) \big) \label{TPDFBB}\\
h^{AB}_N(x_1,x_2,t) &=  \frac{1}{N_A N_B} \sum_{i=1}^{N_A}\sum_{m=1}^{N_B} \delta_{X^A_i(t),X^B_m(t)}(x_1,x_2).\label{TPDFAB}
\end{align} 
\noindent Here, $h^{AA}_N$ and $h^{BB}_N$ give in the large $N_A,N_B$ limit the probabilities of finding pairs of not necessarily linked particles of the same species around $x_1$ and $x_2$, while $h^{AB}_N(x_1,x_2,t)$ gives the probability of finding a particle of type A around $x_1$ and a particle of type B around $x_2$.

The kinetic system provides evolution equations for $f^A,f^B, g^{AA}, g^{BB}, g^{AB}$, the large $N$-limits of the corresponding empirical densities 
defined above. Since this derivation follows closely the works of \cite{Barre2016, Degond_Peurichard_2016} adapted to a two species system, we leave the details in Appendix \ref{AppendixA}. 
The fast network remodelling limit taken on this kinetic system then formally yields a macroscopic system of evolution equations involving only $f^A$ and $f^B$,
i.e. the macroscopic evolution we are looking for, given by \eqref{Eqfmacro}, below.

\subsubsection*{Sketch of  Approach II}
%In this section, we sketch a formal derivation of this averaged microscopic model \eqref{eq:averagedA}-\eqref{eq:averagedB}.
 We denote by $A_{ij}(t),B_{ij}(t),C_{ij}(t)$ the adjacency matrices of particles $A$, $B$, and
cross-links $A-B$ respectively. In particular, for $i,j\in \{1,\ldots,N_A\}$, $A_{ij}(t)=1$ (resp. $=0$) if particles of type $A$ $i$ and $j$ are connected at time $t$ (resp. not connected). The definition of matrix $B$ is similar. For $i\in  \{1,\ldots,N_A\}$, $j \in  \{1,\ldots,N_B\}$, $C_{ij}(t)=1$ (resp. $=0$) if
particle $i$ of type $A$ and particle $j$ of type $B$ are connected at time $t$ (resp. not connected). $A$ and $B$ are square symmetric matrices, and $C$ is an $N_A\times N_B$ rectangular matrix. 

The derivation of the reduced microscopic model relies on averaging. The diffusions of particles positions $X_i^A, X_j^B(t)$ are slow processes, and the links $A_{ij}(t),B_{ij}(t),C_{ij}(t)$ are fast processes: they quickly converge to stationary measures which depend on $X_i^{A}(t)$ and $X_i^{B}(t)$. We will then compute the evolution of $X_i^{A}(t)$ and $X_i^{B}(t)$ by averaging the basic dynamical equations \eqref{IBMA}-\eqref{IBMB} over these stationary measures of the link processes.

The process for the links is written
\begin{eqnarray}
dA_{ij}(t) &=& -A_{ij}(t)dN^{AA,d}_{ij}(t) +[1-A_{ij}(t)] \chi_{\{|X^A_i(t)-X^A_j(t)|\leq R\}} dN^{AA,c}_{ij}(t) \label{fastAA} \\
dB_{ij}(t) &=& -B_{ij}(t)dN^{BB,d}_{ij}(t) +[1-B_{ij}(t)] \chi_{\{|X^B_i(t)-X^B_j(t)|\leq R\}} dN^{BB,c}_{ij}(t) \label{fastBB} \\
dC_{ij}(t) &=& -C_{ij}(t)dN^{AB,d}_{ij}(t) +[1-C_{ij}(t)] \chi_{\{|X^A_i(t)-X^B_j(t)|\leq R\}} dN^{AB,c}_{ij}(t) \label{fastAB} 
\end{eqnarray}
where the $N^{AA,d}_{ij},N^{AA,c}_{ij},N^{BB,d}_{ij},N^{BB,c}_{ij},N^{AB,d}_{ij},N^{AB,c}_{ij},$ are independent Poisson processes with rates $\nu_{d,\ep}^{ST}$ for destruction of the link connecting particles of type $S$ and $T$, $S,T\in\{A,B\}$, and $\nu^{ST}_{c,N,\ep}$ for creation of a link between particles of type $S$ and $T$. We will moreover consider the following scaling of these rates
\eq{
\nu^{AA}_{d,\ep}&=\nu^{AA}_d\varepsilon^{-2},\quad \nu^{AA}_{c,N,\ep}=\nu^{AA}_{c,\ep} N_A^{-1} = \nu^{AA}_c N_A^{-1}\varepsilon^{-2},\\
\nu^{BB}_{d,\ep}&=\nu^{BB}_d\varepsilon^{-2},\quad\nu^{BB}_{c,N,\ep}=\nu^{BB}_{c,\ep} N_B^{-1}=\nu^{BB}_c N_B^{-1}\varepsilon^{-2},\\
\nu^{AB}_{d,\ep}&=\nu^{AB}_d\varepsilon^{-2},\quad \nu^{AB}_{c,N,\ep}=\nu^{AB}_{c,\ep} N_B^{-1}=\nu^{AB}_c N_B^{-1}\varepsilon^{-2}.  \label{eq:nuscalings}
}
The subscript $\ep$ (resp. $N_{A}, N_B$) signals a dependency on $\ep$ (resp. $N_{A}, N_B$), and the rates without $N_{A}, N_B,\ep$ subscripts are assumed to be independent of $N_{A},N_B,\ep$. Relations \eqref{eq:nuscalings} make  the scaling with $\ep$ and $N_{A}, N_B$ needed to perform the limit procedures in \eqref{eq:CD} explicit. In particular, the scaling with $N_{A}, N_B$ ensures that the connectivity of any particle remains of order~$1$, and the scaling with $\ep$ controls the speed of the linking/unlinking process.

Conditionally on the positions $X^A_i,X^B_j$, all the processes $A_{ij},B_{ij},C_{ij}$ are independent. The stationary measures of \eqref{fastAA}-\eqref{fastBB}-\eqref{fastAB}, for fixed positions $X^A_i,X^B_j$ are then simply product of Bernoulli measures ($\mathbb{P}$ denotes the probability):
\begin{eqnarray}
\mathbb{P}(A_{ij}(t)=1)=
\frac{\frac{\nu^{AA}_c}{N_A} \chi_{\{|X^{A}_i(t)-X^A_j(t)|\leq R\}} }{\frac{\nu^{AA}_c}{N_A} +\nu^{AA}_d} &,&
\mathbb{P}(A_{ij}(t)=0)= 1- \mathbb{P}(A_{ij}(t)=1)\\ %1- \frac{\frac{\tilde{\nu}^{AA}_f}{N} \chi_{|X^A_i(t)-X^A_j(t)|\leq R}}{\frac{\tilde{\nu}^{AA}_f}{N}+\tilde{\nu}^{AA}_d} \nonumber \\
\mathbb{P}(B_{ij}(t)=1)=\frac{\frac{\nu^{BB}_c}{N_B} \chi_{\{|X^{B}_i(t)-X^B_j(t)|\leq R\}} }{\frac{\nu^{BB}_c}{N_B} +\nu^{BB}_d} &,&
\mathbb{P}(B_{ij}(t)=0)= 1-\mathbb{P}(B_{ij}(t)=1) \\%1- \frac{\frac{\tilde{\nu}^{BB}_f}{N} \chi_{|X^B_i(t)-X^B_j(t)|\leq R}}{\frac{\tilde{\nu}^{BB}_f}{N}+\tilde{\nu}^{BB}_d} \nonumber \\
\mathbb{P}(C_{ij}(t)=1)=\frac{\frac{\nu^{AB}_c}{N_B} \chi_{\{|X^{A}_i(t)-X^B_j(t)|\leq R\}} }{\frac{\nu^{AB}_c}{N_B} +\nu^{AB}_d} &,&
\mathbb{P}(C_{ij}(t)=0)=1-\mathbb{P}(C_{ij}(t)=1)% 1- \frac{\frac{\tilde{\nu}^{AB}_f}{N} \chi_{|X^A_i(t)-X^B_j(t)|\leq R}}{\frac{\tilde{\nu}^{AB}_f}{N}+\tilde{\nu}^{AB}_d} \nonumber 
\end{eqnarray}
For $N_A,N_B$ large, the above expressions simplify as the $O(1/N_A, 1/N_B)$ terms in the denominators are negligible. 
One can write the equations for the positions, averaged over the stationary measure for the links; calling $\tilde{X}_i^A,\tilde{X}_j^B$ these new processes, we obtain (neglecting terms of order $1/N_A, 1/N_B$)

\eq{\label{eq:averagedA}
d\tilde{X}_i^A =& - \mu \left( \frac{1}{N_A}\frac{\nu^{AA}_c}{\nu^{AA}_d}\sum_{j=1}^{N_A} \chi_{\{|\tilde{X}_i^A-\tilde{X}_j^A|\leq R\}}\nabla  \Phi^{AA}(\tilde{X}_i^A-\tilde{X}_j^A) \right. \\
&\qquad\qquad+\left. \frac{1}{N_B}\frac{\nu^{AB}_c}{\nu^{AB}_d}\sum_{j=1}^{N_B}  \chi_{\{|\tilde{X}_i^A-\tilde{X}_j^B|\leq R\}}\nabla \Phi^{AB}(\tilde{X}_i^A-\tilde{X}_j^B) \right)dt+\sqrt{2D^A} dB_i^A   
}
\eq{ \label{eq:averagedB}
d\tilde{X}_i^B =& -\mu\left( \frac{N_A}{N_B} \frac{1}{N_A}\frac{\nu^{AB}_c}{\nu^{AB}_d}\sum_{j=1}^{N_A}  \chi_{\{|\tilde{X}_i^B-\tilde{X}_j^A|\leq R\}}\nabla \Phi^{BA}(\tilde{X}_i^B-\tilde{X}_j^A) \right. \\
&\qquad\qquad+\left. \frac{1}{N_B}\frac{\nu^{BB}_c}{\nu^{BB}_d}\sum_{j=1}^{N_B}  \chi_{\{|\tilde{X}_i^B-\tilde{X}_j^B|\leq R\}} \nabla \Phi^{BB}(\tilde{X}_i^B-\tilde{X}_j^B)\right)dt +\sqrt{2D^B} dB_i^B.
}
%When the number of particles and links tends to infinity, we assume that 
% \eqh{\frac{\tilde\nu^{K_{AA}}_{c}} {N_A}\to \nu^{AA}_c,\quad \frac{\tilde\nu^{K_{BB}}_{c}}{ N_B} \to\nu^{BB}_c, \quad \frac{\tilde\nu^{K_{AB}}_{c}}{ N_B} \to\nu^{AB}_c,\\
% \tilde\nu^{K_{AA}}_{d}\to \nu^{AA}_d,\quad \tilde\nu^{K_{BB}}_{d}\to\nu^{BB}_d, \quad \tilde\nu^{K_{AB}}_{d}\to\nu^{AB}_d,
 %}
 % this scaling of the Poisson rates with the number of particles ensures that the mean connectivity of particles remains of order $1$. 
 Here we tacitly assumed the translation invariance of the potential $\Phi^ST(X_i,X_j)=\Phi^ST(X_i-X_j)$, which can be relaxed, see the Appendix.
From Eqs. \eqref{eq:averagedA}-\eqref{eq:averagedB}, one can deduce
the following proposition, describing the the dynamics of the particles density $f^A$ and $f^B$ in the large $N_A,N_B$ limit:
\begin{prop}\label{prop1}
Assume Eqs. \eqref{eq:averagedA}-\eqref{eq:averagedB} and that the potentials are radially symmetric $\Phi^{ST}(X_i-X_j)=\Phi^{ST}(|X_i-X_j|$; then in the limit $N_A,N_B\to \infty, N_A/N_B \to r_{AB}>0$, the one particle distribution functions
$f^A$ and $f^B$ are solution of the system
\eq{\label{Eqfmacro}
\partial_t f^A &= D^A \Delta_x f^A +\nabla_x \cdot \bigg(f^A(x,t) \nabla_x \big(\tilde{\Phi}^{AA} \ast f^A \big) (x,t) \bigg) + \nabla_x \cdot \bigg(f^A(x,t) \nabla_x  \big(\tilde{\Phi}^{AB} \ast f^B\big) (x,t) \bigg), \\
\partial_t f^B &= D^B \Delta_x f^B +\nabla_x \cdot \bigg(f^B(x,t) \nabla_x \big(\tilde{\Phi}^{BB} \ast f^B \big) (x,t) \bigg) +  \nabla_x \cdot \bigg(f^B(x,t) \nabla_x  \big(\tilde{\Phi}^{BA} \ast f^A\big) (x,t) \bigg),
}
where the potentials  $\tilde{\Phi}^{ST}$ are given by:
% obtained from the 
% $\Phi^{ST}(X_i,X_j)=\Phi^{ST}(X_i-X_j)$ by 
% truncating at distance $R$ and reweighing according to the rates of the Poisson processes and $r_{AB}$:
\begin{eqnarray}
\tilde{\Phi}^{AA}(x) &=& \frac{\nu_c^{AA}}{\nu_d^{AA}} \left(\Phi^{AA}(x) \chi_{\{|x|\leq R\}}+ \Phi^{AA}(R)\chi_{\{|x|>R\}}\right), \\
\tilde{\Phi}^{BB}(x) &=& \frac{\nu_c^{BB}}{\nu_d^{BB}} \left(\Phi^{BB}(x) \chi_{\{|x|\leq R\}}+ \Phi^{BB}(R)\chi_{\{|x|>R\}}\right), \\
\tilde{\Phi}^{AB}(x) &=& \frac{\nu_c^{AB}}{\nu_d^{AB}} \left(\Phi^{AB}(x) \chi_{\{|x|\leq R\}}+ \Phi^{AB}(R)\chi_{\{|x|>R\}}\right), \\
\tilde{\Phi}^{BA}(x) &=& r_{AB}\frac{\nu_c^{AB}}{\nu_d^{AB}} \left(\Phi^{BA}(x) \chi_{\{|x|\leq R\}}+ \Phi^{BA}(R)\chi_{\{|x|>R\}}\right), 
\end{eqnarray}
and $\ast$ denotes the convolution.
\end{prop}

An alternate proof of this proposition based on Approach I is presented in the Appendix. The Hookean potential considered in this paper satisfies the radial symmetry assumption. However, for purposes of derivation of the macroscopic model, this assumption can be relaxed at the expense of more complex formulas (see the Appendix). 
The macroscopic model consists of an aggregation-diffusion equations with nonlocal terms, where each particle interacts with its close neighbors of the same family (second term of the right hand side of Eq. \eqref{Eqfmacro}) as well as with the ones of the other family (third term of the right hand side of Eq. \eqref{Eqfmacro}), and where the diffusive term corresponds to the Brownian motion of individual particles. In the next section, we perform the linear stability analysis to identify phase transitions of the homogeneous steady-state.

%\red{I have used potentials that depend only on one variable: $\Phi(x_1,x_2)= \Phi(|x_1-x_2|)$. Is it useful to do something more general, as done at the 
%beginning of this section?}
%\red{I would do everything with $\Phi(|x_1-x_2|)$... (J)}

\section{Analysis of the macroscopic system in the whole space}\label{sec:3}
\subsection{Linear stability in the whole space}\label{sec:3.1}
In this section, we perform a linear stability analysis of the macroscopic model. We recall the macroscopic equations for $f^A, f^B$:
\begin{align*}
\partial_t f^A &= D^A \Delta_x f^A + \nabla_x \cdot \bigg(f^A \nabla_x \big(\tilde{\Phi}^{AA} \ast f^A \big) \bigg) +  \nabla_x \cdot \bigg(f^A \nabla_x  \big(\tilde{\Phi}^{AB} \ast f^B\big) \bigg)\\
\partial_t f^B &= D^B \Delta_x f^B + \nabla_x \cdot \bigg(f^B \nabla_x \big(\tilde{\Phi}^{BB} \ast f^B \big) \bigg) +  \nabla_x \cdot \bigg(f^B \nabla_x  \big(\tilde{\Phi}^{BA} \ast f^B\big) \bigg),
\end{align*} 
\noindent where the factors  $\frac{\nu^{ST}_c}{\nu^{ST}_d}$ and $r_{AB}$ have been included in the potential functions $\tilde{\Phi}^{ST}$. We first linearize around the homogeneous steady states, i.e $f_*^A = const.$, $f_*^B = const.$. Writing
$$
f^A = f_*^A + \tilde{f}^A, \; f^B = f_*^B + \tilde{f}^B,
$$ 
\noindent we have , at the first order:
\eq{\label{Eqfmacrosimp}
\partial_t f^A &= f_*^A \Delta \big[ f^A \ast \tilde\Phi^{AA} +  f^B \ast \tilde\Phi^{AB} \big] + D^A f^A \\
\partial_t f^B &= f_*^B \Delta \big[ f^B \ast \tilde\Phi^{BB} +  f^A \ast \tilde\Phi^{BA} \big] + D^B f^B.
}
\noindent For any integrable function $F$ and a vector $y\in \mathbb{R}^2$ we recall the definition of  the spatial Fourier transform
$$
\hat{F}(y) = \frac{1}{2\pi}\int_{\mathbb{R}^2} \exp^{-ix \cdot y} F(x) dx.
$$
Applying it  to both sides of \eqref{Eqfmacrosimp}, we obtain the following system:
\eq{\label{Fourier}
\partial_t \begin{pmatrix} \hat{f}^A\\ \hat{f}^B\end{pmatrix}(y,t) &= \begin{pmatrix} - f_*^A |y|^2  \big(2\pi \hat{\tilde\Phi}^{AA}(y) + \frac{D^A}{f_*^A}\big) & -f_*^A |y|^2  2\pi \hat{\tilde\Phi}^{AB}(y)\\
-f_*^B |y|^2  2\pi\hat{\tilde\Phi}^{BA}(y) & -f_*^B |y|^2  \big(2\pi\hat{\tilde\Phi}^{BB}(y) + \frac{D^B}{f_*^B}\big) \end{pmatrix} \begin{pmatrix} \hat{f}^A\\ \hat{f}^B\end{pmatrix}(y,t) \\
&:= M(y)\begin{pmatrix} \hat{f}^A\\ \hat{f}^B\end{pmatrix}(y,t).
}
\noindent Therefore, 
$$\begin{pmatrix} \hat{f}^A\\ \hat{f}^B\end{pmatrix}(y,t) = {c}_1(y) \exp^{\lambda_1(y) t} \vec{u}_1(y) + {c}_2(y) \exp^{\lambda_2(y) t} \vec{u}_2(y),$$
\noindent where $\lambda_1(y), \lambda_2(y)$ are the eigenvalues of the matrix $M(y)$ and $\vec{u}_1(y),\vec{u}_2(y)$ are the corresponding eigenvectors. 

\subsubsection*{General case}

In the general case, the homogeneous steady state will be unstable if at least one of the eigenvalues of matrix $M(y)$ is positive. By computing the determinant of M:
$$
\Delta(M) = |y|^4f_*^A f_*^B \bigg[\left(2\pi\hat{\tilde\Phi}^{AA} + \frac{D^A}{f_*^A}\right)\left(2\pi\hat{\tilde\Phi}^{BB} + \frac{D^B}{f_*^B}\right) - (2\pi)^2\hat{\tilde\Phi}^{AB}\hat{\tilde\Phi}^{BA}\bigg],
$$	
\noindent we can see that, for general interaction potentials, the constant steady states will be unstable if one of the two conditions is met:
\begin{itemize}
\item 1. $\Delta (M) < 0$: $(2\pi\hat{\tilde\Phi}^{AA} + \frac{D^A}{f_*^A})(2\pi\hat{\tilde\Phi}^{BB} + \frac{D^B}{f_*^B}) < 4\pi^2 
\hat{\tilde\Phi}^{AB}\hat{\tilde\Phi}^{BA}$
\item 2. $\Delta(M) \geq 0$ and $Tr(M)>0$: $\begin{cases}
(2\pi\hat{\tilde\Phi}^{AA} + \frac{D^A}{f_*^A})(2\pi\hat{\tilde\Phi}^{BB} + \frac{D^B}{f_*^B}) \geq 4\pi^2\hat{\tilde\Phi}^{AB}\hat{\tilde\Phi}^{BA}\\
f_*^A(2\pi\hat{\tilde\Phi}^{AA} + \frac{D^A}{f_*^A}) + f_*^B(2\pi\hat{\tilde\Phi}^{BB} + \frac{D^B}{f_*^B})<0.
\end{cases}$
\end{itemize}

\subsubsection*{Hookean interaction potentials}

\begin{figure}[h]
\includegraphics[scale=0.4]{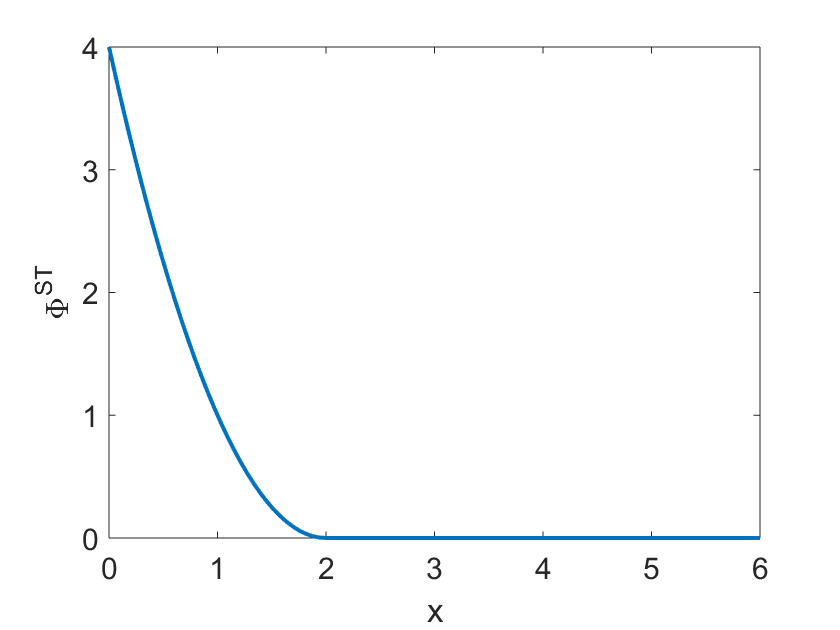}
\caption{Form of the potential $\tilde{\Phi}^{ST}(x)$ for $R=2$ and $\frac{\nu^{ST}_c}{\nu^{ST}_d} \frac{\kappa^{ST}}{2} = 1$. The potential is repulsive on its support. \label{PhiST}}
\end{figure}

\begin{figure}[h]
\includegraphics[scale=0.3]{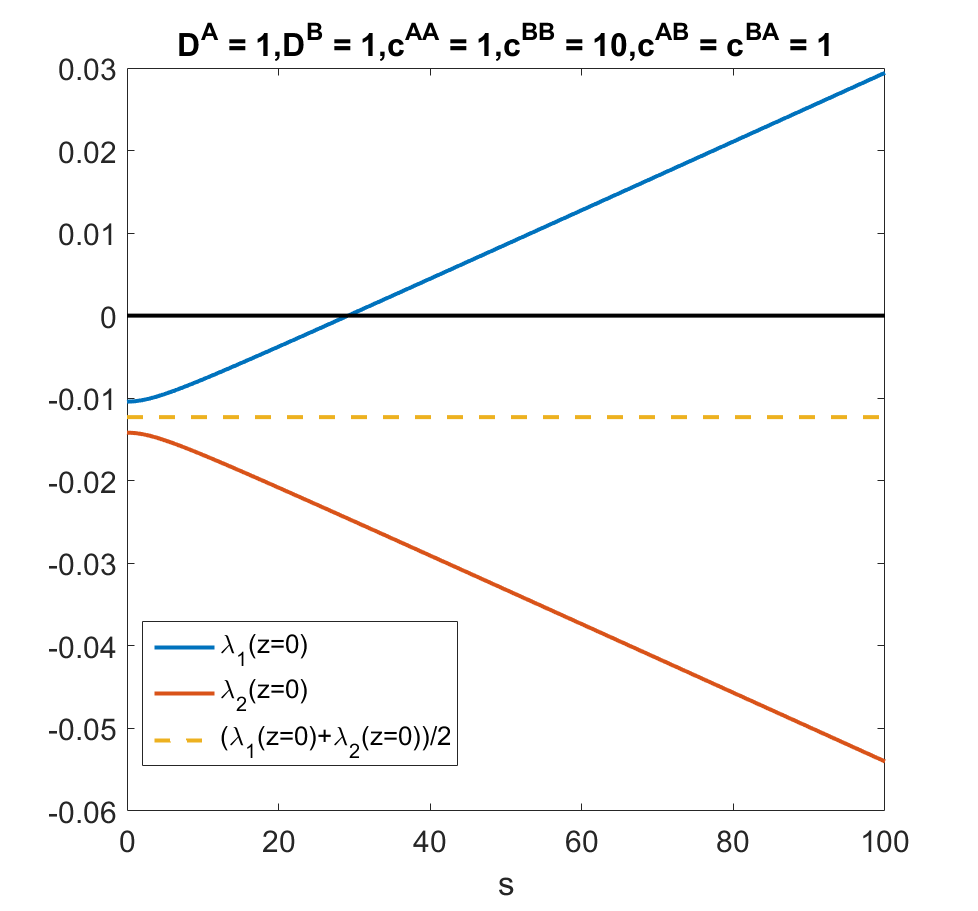}
\caption{Values of $\lambda_1$ (blue curve), $\lambda_2$ (orange curve) and their mean (yellow dotted line) near $z=0$ (z=0.1), plotted as functions of the scaling parameter $s$ for $R=1$, $D^A=D^B = 1, c'^{AA}=c'^{AB}=c'^{BA}=1, c'^{BB}=10$.\label{PloteigV0}}
\end{figure}

\begin{figure}[h]
\includegraphics[scale = 0.4] {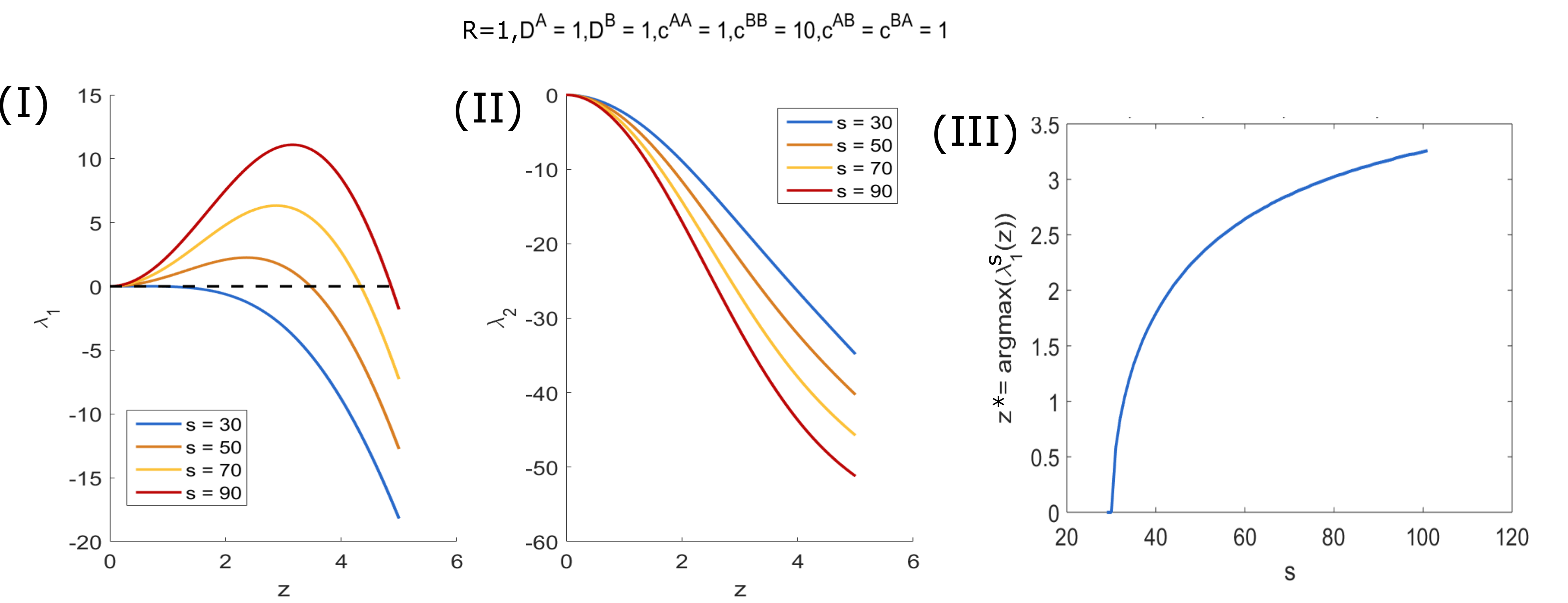}
\caption{(I): Values of $\lambda_1(z)$ as functions of $z$ for $R=1$, $D^A=D^B = 1,\, c^{AA}=c^{AB}=c^{BA}=1,\, c^{BB}=10$ and for different values of $s$ in the instability regime: $s=30$ (blue curve), $s=50$ (orange curve), $s=70$ (yellow curve), and $s=90$ (red curve). (II): same plots for $\lambda_2(z)$. (III) Plot of $z^*$ defined in \eqref{z_star} as a function of parameter $s$. \label{PloteigV}}
\end{figure}

In order to explicit these conditions as functions of the model parameters, we now suppose that the intra- and inter- species links act as springs of rest length $R$ between the particles. As the detection radius for the interaction is also $R$, this amounts to consider that particles only repulse each other up until distance~$R$. To keep  enough generality, we consider different interaction intensities between the type A and type B intra- and inter- species springs:
\begin{align*}
\Phi^{AA}(x_1,x_2) = \frac{\kappa^{AA}}{2}(|x_1-x_2| - R)^2, \quad \Phi^{BB}(x_1,x_2) = \frac{\kappa^{BB}}{2}(|x_1-x_2| - R)^2\\
\Phi^{AB}(x_1,x_2) = \frac{\kappa^{AB}}{2}(|x_1-x_2| - R)^2,\quad \Phi^{AB}(x_1,x_2) = \frac{\kappa^{BA}}{2}(|x_1-x_2| - R)^2, 
\end{align*}  
% with $\kappa^{AA}, \kappa^{BB}, \kappa^{AB},\kappa^{BA}\geq 0$ 
\noindent and not necessarily mutually equal. We have:
$$
\int f^T(x') \chi_{\{|x-x'|\leq R\}} \nabla_x {\Phi}^{ST}(x,x')dx' =  \kappa^{ST} \int f^T(x')(|x-x'|-R) \frac{x-x'}{|x-x'|}\chi_{\{|x-x'|\leq R\}}dx'.
$$
\noindent Remembering that the factors $\frac{\nu^{ST}_c}{\nu^{ST}_d}$ were included in the potential functions $\tilde{\Phi}^{ST}$ and following derivation \eqref{fg_app0}-\eqref{fg_app3} from the Appendix \ref{scaling}, we may introduce:
\eq{\label{Hookean}
\tilde{\Phi}^{ST}(x) = \frac{\nu^{ST}_c}{\nu^{ST}_d} \frac{\kappa^{ST}}{2}
\begin{cases}
(|x|-R)^2 \qquad \text{for $|x|\leq R$} \\
0 \qquad \qquad \qquad \text{for $|x|> R$}
\end{cases},
}
\noindent which ensures that the equations for $f^A,\, f^B$ are of the form \eqref{Eqfmacro} and the linearized system is of the form \eqref{Eqfmacrosimp}.  The form of the new potential is plotted in Fig. \ref{PhiST}.
 We can now compute the Fourier transform of $\tilde{\Phi}^{ST}$.
% $$
% \hat{\tilde{\Phi}}^{ST}(y) = \frac{1}{2\pi}\int_{\mathbb{R}^2} \exp^{-ix \cdot y} \tilde{\Phi}^{ST}(x) dx.
% $$
Using the radial symmetry of $\tilde{\Phi}^{ST}$ and denoting $ \Upsilon=|y|$ we can show that
$$
\hat{\tilde{\Phi}}^{ST}(\Upsilon) =  \frac{\nu^{ST}_c}{\nu^{ST}_d} \kappa^{ST} \bigg(J_0(\Upsilon R) \frac{R^2}{\Upsilon^2} - J_1( \Upsilon R)\frac{2R}{\Upsilon^3}+ \frac{\pi R^2}{2\Upsilon^2} \big[ J_1(\Upsilon R) H_0(\Upsilon R) - J_0(\Upsilon R) H_1(\Upsilon R)\big] \bigg),
$$
\noindent where $J_0,J_1$ are the Bessel functions of the first kind of order 0 and 1:
$$
J_i(x) = \sum_{m=0}^\infty \frac{(-1)^m}{m!\Gamma(m+1+i)}\lr{\frac{x}{2}}^{2m+i}
$$
\noindent and $H_0,H_1$ are the Struve functions defined by:
$$
H_i(x) = \sum_{m=0}^\infty \frac{(-1)^m}{\Gamma(m+3/2)\Gamma(m+3/2+i)}\lr{\frac{x}{2}}^{2m+i+1}.
$$
\noindent We refer the reader to \cite{Barre2016} for the computation of these terms. Replacing the Fourier transforms of the potentials by their expressions in $M(y)$, we can write
$$
M(y) = -\begin{pmatrix}  c^{AA} H(y) + D^A |y|^2  & c^{AB} H(y)\\
c^{BA} H(y) & c^{BB} H(y) + D^B |y|^2 \end{pmatrix},
$$
\noindent with $c^{ST} = \frac{2\pi \kappa^{ST} f_*^S \nu^{ST}_c R^2}{\nu^{ST}_d}$ for all $S, T\in\{A,B\}$ and 
\begin{align*}
H(y) &= J_0(|y|R) - J_1(|y|R)\frac{2}{|y|R}+ \frac{\pi}{2} \big[ J_1(|y|R) H_0(|y|R) - J_0(|y|R) H_1(|y|R)\big] \\
&= \frac{\pi}{2} \big[J_1(|y|R) H_0(|y|R) - J_0(|y|R) H_1(|y|R)\big] - J_2(R|y|).
\end{align*}
\noindent Writing $z = |y|R$, the determinant of $M$ can now be written:
$$
\Delta(M) = \frac{1}{R^4}(D^A z^2 +c'^{AA}\tilde{H}(z))(D^B z^2 +c'^{BB}\tilde{H}(z)) - c'^{AB}c'^{BA}\tilde{H}(z)^2,
$$
\noindent where $c'^{ST} = R^2 c^{ST} =  \frac{2\pi \kappa^{ST} f_*^S \nu^{ST}_c R^4}{\nu^{ST}_d}$ and 
\begin{align*}
\tilde H(z) &= \frac{\pi}{2} \big[J_1(z) H_0(z) - J_0(z) H_1(z)\big] - J_2(z).
\end{align*}
\noindent Now, lengthy but straightforward computations show that for $z$ close to the origin we have:
\eq{\label{approx_space}
\tilde H(z) = \frac{1}{24} z^2 + O(z^4),
}
\noindent and so, close to the origin $z=0$, we have
\eq{\label{deltaM}
\Delta(M) &= \frac{z^4}{R^4}\bigg((D^A +\frac{c'^{AA}}{24})(D^B +\frac{c'^{BB}}{24}) - \frac{c'^{AB}c'^{BA}}{24^2}\bigg)\\
&=\frac{z^4}{R^4}\bigg(D^AD^B +\frac{1}{24}(D^A c'^{BB} +D^B c'^{AA}) + \frac{c'^{AA}c'^{BB}-c'^{AB}c'^{BA}}{24^2}\bigg).
}
\noindent In order to simplify the analysis, we suppose the following hypothesis:
\begin{hypot}\label{hyp:1}
The intraspecies links generate repulsive potentials, i.e $\kappa^{AA}, \kappa^{BB}>0$.
\end{hypot}
\noindent We first note that under Hypothesis \ref{hyp:1}, the trace of $M$
$$
Tr(M) = -\frac{z^2}{R^2} \bigg[D^A+D^B +\frac{c'^{AA}+c'^{BB}}{24}\bigg] +o(z^4),
$$
\noindent is negative for small $z$. 
% \blue{\emph{
% Maybe recall the argument to exclude instability for large $z$? (useful at least for me!)}} 
Therefore, the homogeneous steady states will be unstable for small $z$ only if $\Delta(M)\leq0$. Note also that the parameters $c'^{AB}$ and $c'^{BA}$ should have the same sign to allow the determinant of $M$ to be negative, otherwise the homogeneous steady state will be a stable case. Now, we scale the interspecies link potential intensities with a parameter $s \in \mathbb{R}$ such that $\kappa^{AB} = s \tilde{\kappa}^{AB},\kappa^{BA} = s \tilde{\kappa}^{BA}$. It corresponds to the relevant scaling of the parameters $c'^{AB}$ and $c'^{BA}$, to simplify the notation we denote the corresponding reference values by the same symbols $c'^{AB}$ and $c'^{BA}$, then
$$
\Delta(M) = \frac{z^4}{R^4}\bigg(D^AD^B +\frac{1}{24}(D^A c'^{BB} +D^B c'^{AA}) + \frac{c'^{AA}c'^{BB}-s^2{c}'^{AB}{c}'^{BA}}{24^2}\bigg),
$$
\noindent and we immediately note that $s$ must be large enough to allow $\Delta(M)$ to be negative. More precisely, the two eigenvalues of $M$ are written:
\eq{\label{l1l2}
\lambda_1 &= \frac{1}{2}\bigg(Tr(M) + \sqrt{Tr(M)^2 - 4\Delta(M)}\bigg) = \frac{z^2}{2 R^2} \bigg(-C_+ + \sqrt{C_-^2 + s^2 \frac{c'^{AB}c'^{BA}}{144}} \bigg)\\
\lambda_2 &= \frac{1}{2}\bigg(Tr(M) - \sqrt{Tr(M)^2 - 4\Delta(M)}\bigg)= \frac{z^2}{2R^2} \bigg(-C_+ - \sqrt{C_-^2 + s^2 \frac{c'^{AB}c'^{BA}}{144}} \bigg),
}
\noindent where $C_+ = D^A+D^B+\frac{c'^{AA}+c'^{BB}}{24}$ and $C_- = D^A-D^B+\frac{c'^{AA}-c'^{BB}}{24}$. We can therefore plot their values as functions of $s$ near $z=0$ (see Fig. \ref{PloteigV0} for $z=0.1$ and parameter values $D^A=D^B = 1, c'^{AA}=c'^{AB}=c'^{BA}=1, c'^{BB}=10$).
As shown by Fig. \ref{PloteigV0}, there exists a critical value $s^*$ of $s$ such that for $s>s^*$ the homogeneous state will be unstable (i.e we will observe cell aggregates). From the definition of $s$, this means that the interspecies repulsion force must be large enough to compensate the intraspecie repulsion and diffusion and enable the two species to separate into clusters. These first results are in accordance with the observations of \cite{Taylor2017}. Note that, by  equating \eqref{deltaM} to $0$, we can directly compute the value of $s^*$ as a function of the model parameters:
\eq{\label{s_star}
s^* = \sqrt{\frac{576}{c'^{AB}c'^{BA}} \big(D^A+\frac{c'^{AA}}{24}\big)\big(D^B+\frac{c'^{BB}}{24}\big) }.
}
In Fig. \ref{PloteigV}, we aim to plot the values of $\lambda_1(z),\lambda_2(z)$ in the unstable regime $s>s^*$, $s^*$ being determined on Fig.\ref{PloteigV0} (critical value of $s$ such that $\lambda_1(z\approx 0)$ becomes positive). 
We select four values of $s>s^*$ and for each of them we compute the functions $\lambda_1(z),\lambda_2(z)$ using matrix $M(z)$ before Taylor expanding it near 0. The Bessel and Struve functions are approximated numerically. 

As one can see in Fig. \ref{PloteigV} (I), $\lambda_1$ is an increasing function of $s$ at fixed $z$. Furthermore, the most unstable parameter $z$, i.e. the value $z^*$ for which $\lambda_1(z)$ reaches its maximal value, defined by:
\eq{\label{z_star}
z^*(s) = \underset{z\in \mathbb{R}^+}{\text{argmax}} \; \lambda_1^s(z),}
increases with $s$ (see Fig. \ref{PloteigV} (III)). Hence, at small times, the instability should be dominated by Fourier modes with parameter around $z^\ast(s)$. Assuming this remains qualitatively true in the time asymptotic regime, one then expects that larger values of $s$ will lead to more aggregated (smaller) patterns.
%$$
%z^*(s) = \frac{C}{L(s)}.
%$$
%In Fig. \ref{PloteigV}, we plot the values of $z^*$ as a function of $s$, and observe that $z^*$ is an increasing function of $s$ (see the discussion section for more details). 
% \begin{figure}[h]
%\includegraphics[scale = 0.3]{LinearStab_eigV_ystar.png}
%\caption{Plot of $z^*$ defined in \eqref{z_star} as a function of parameter $s$. \label{PloteigV}}
%\end{figure}

\subsection{Linear stability in the periodic box}\label{sec:3.2}
For the sake of numerical simulations, we now interpret the above results in the case of a space-periodic domain. In practice, instead of the whole plane we consider the two-dimensional square periodic domain $[-L,L]\times[-L,L]$, however, for $L$ sufficiently large the results change quantitatively but not qualitatively. 
Introducing the shorthand notation for the Fourier modes
\eq{ e_{k_1,k_2}=\exp{\left[\frac{i\pi}{L} \lr{k_1x_1+k_2x_2}\right]},\label{not_e}}
we may write for $S=A,B$:
\eqh{
f^S(x_1,x_2)=\sum_{k_1,k_2\in \mathbb{Z}}\hat f^S_{k_1,k_2}e_{k_1,k_2},
}
where the Fourier coefficients $\hat f^S_{k_1,k_2}$ are given by
\eqh{
\hat f^S_{k_1,k_2}=\frac{1}{4L^2}\int_{-L}^{L}\int_{-L}^{L} f^S(x_1,x_2)e_{-k_1,-k_2}\, dx_1\, dx_2.}
The periodic analog of the macroscopic system after Fourier transform is (see \eqref{Fourier} for comparison):
\eq{\label{FourierP}
\partial_t \begin{pmatrix} \hat{f}_{k_1,k_2}^A\\ \hat{f}_{k_1,k_2}^B\end{pmatrix}&= \begin{pmatrix} - f_*^A \frac{\pi^2(k_1^2+k_2^2)}{L^2}  \big(4L^2 \hat{\tilde\Phi}_{k_1,k_2}^{AA} + \frac{D^A}{f_*^A}\big) & -f_*^A \frac{\pi^2(k_1^2+k_2^2)}{L^2}  4L^2 \hat{\tilde\Phi}_{k_1,k_2}^{AB}\\
-f_*^B \frac{\pi^2(k_1^2+k_2^2)}{L^2}  4L^2\hat{\tilde\Phi}_{k_1,k_2}^{BA} & -f_*^B \frac{\pi^2(k_1^2+k_2^2)}{L^2} \big( 4L^2\hat{\tilde\Phi}_{k_1,k_2}^{BB} + \frac{D^B}{f_*^B}\big) \end{pmatrix} \begin{pmatrix} \hat{f}_{k_1,k_2}^A\\ \hat{f}_{l_1,l_2}^B\end{pmatrix} \\
&:= M_{k_1,k_2}\begin{pmatrix} \hat{f}_{k_1,k_2}^A\\ \hat{f}_{k_1,k_2}^B\end{pmatrix}.
}
Since $f^A, f^B$ are both probability measures, we take
$$f_*^A= f_*^B=\frac{1}{4L^2},$$
which means that the matrix $M$ from \eqref{FourierP} has now the form
\eq{
M_{k_1,k_2}=\begin{pmatrix} -  \frac{\pi^2(k_1^2+k_2^2)}{L^2}  \big( \hat{\tilde\Phi}_{k_1,k_2}^{AA} + D^A\big) & - \frac{\pi^2(k_1^2+k_2^2)}{L^2}  \hat{\tilde\Phi}_{k_1,k_2}^{AB}\\
- \frac{\pi^2(k_1^2+k_2^2)}{L^2}  \hat{\tilde\Phi}_{k_1,k_2}^{BA} & -\frac{\pi^2(k_1^2+k_2^2)}{L^2} \big(\hat{\tilde\Phi}_{k_1,k_2}^{BB} + D^B\big) \end{pmatrix} .}
This form of the system can be again studied for general potentials. Here we immediately focus on the Hookean case \eqref{Hookean}, for which the Fourier transform equals
\eq{
\hat{\tilde\Phi}^{ST}_{k_1,k_2}=\frac{\nu_c^{ST}}{\nu_d^{ST}}\frac{\pi\kappa^{ST}}{2L^2 }\frac{R^4}{z_{k_1,k_2}^2}
\lr{\frac{\pi}{2} \big[J_1(z_{k_1,k_2}) H_0(z_{k_1,k_2}) - J_0(z_{k_1,k_2}) H_1(z_{k_1,k_2})\big] - J_2(z_{k_1,k_2})},}
where we denoted
\eq{\label{zd}
z_{k_1,k_2}=\frac{\pi R}{L}\sqrt{k_1^2+k_2^2}.}
As explained in \cite{Barre2016}, due to specific shape of the potential, we only need to check the stability of the first mode $(k_1,k_2)=\pm(1,0)$ or $(k_1,k_2)=\pm(0,1)$ in order to find  whether the whole system is stable.
This is the point where the analysis in the space-periodic domain differs from the whole space case. Note that in \eqref{approx_space} we took $z\to 0$, which cannnot be assumed for the discrete values \eqref{zd}. 
Exactly as before, we compute:
$$Tr(M_{1,0})=-\frac{\pi^2}{L^2}\lr{D^A+D^B+\frac{R^2}{2\pi}\lr{\frac{\nu_c^{AA}}{\nu_d^{AA}}\kappa^{AA}+\frac{\nu_c^{BB}}{\nu_d^{BB}}\kappa^{BB}}\tilde H\lr{\frac{\pi R}{L}}},$$
\eq{\label{deltaMd}
\Delta(M_{1,0})=\frac{\pi^4}{L^4}\bigg[ & D^AD^B+\frac{R^2}{2\pi}\lr{D^B\frac{\nu_c^{AA}}{\nu_d^{AA}}\kappa^{AA}+D^A\frac{\nu_c^{BB}}{\nu_d^{BB}}\kappa^{BB}}\tilde H\lr{\frac{\pi R}{L}}\\
& +\frac{R^4}{4\pi^2}\lr{\frac{\nu_c^{AA}}{\nu_d^{AA}}\frac{\nu_c^{BB}}{\nu_d^{BB}}\kappa^{AA}\kappa^{BB}
-s^2\frac{\nu_c^{AB}}{\nu_d^{AB}}\frac{\nu_c^{BA}}{\nu_d^{BA}}\kappa^{AB}\kappa^{BA}}\tilde H\lr{\frac{\pi R}{L}}^2
\bigg].}
Recall however, that assuming Hypothesis \ref{hyp:1} we still have that  $Tr(M_{1,0})$ is always negative.  Based on  numerical simulations close to $z=0$ for the whole space case (see Fig. \ref{PloteigV} I) we observe that  the constant steady state is unstable again only for sufficiently small $z_{1,0}$ and iff $\Delta(M_{1,0})<0$.  This leads to two  eigenvalues $\lambda_1$, $\lambda_2$ of different signs, that can be computed using the formulas \eqref{l1l2}. The value of parameter $s$ that corresponds to the phase transition can be computed by equating \eqref{deltaMd} to $0$ and therefore
\eq{\label{s_star_per}
s^*=\frac{1}{\frac{R^2}{2\pi}\tilde H\lr{\frac{\pi R}{L}}}
\sqrt{\frac{
\lr{D^A+\frac{R^2}{2\pi}\frac{\nu_c^{AA}}{\nu_d^{AA}}\kappa^{AA}\tilde H\lr{\frac{\pi R}{L}}}\lr{D^B+\frac{R^2}{2\pi}\frac{\nu_c^{BB}}{\nu_d^{BB}}\kappa^{BB}\tilde H\lr{\frac{\pi R}{L}}}}{\frac{\nu_c^{AB}}{\nu_d^{AB}}\frac{\nu_c^{BA}}{\nu_d^{BA}}\kappa^{AB}\kappa^{BA}}}.
}

\subsection{On phase separated initial conditions}\label{sec:separatedphases}
The authors of \cite{Taylor2017} report segregation experiments between two types of cells. In particular, in some experiments (see Fig.1 i)-k) in \cite{Taylor2017}) the initial condition is perfectly segregated: at $t=0$ there is a sharp 
interface between a subdomain with only  cells of type $A$ and a subdomain with only  cells of type $B$. Depending on various parameters, the experiments show different possible evolutions  of the interface (see Fig.1 i)-k) in \cite{Taylor2017}):
\begin{itemize}
\item the interface remains sharp and does not move;
\item the interface remains sharp, does not deform, and moves in one direction;
\item the interface becomes blurred: the two species start to mix.
\end{itemize}
Our goal here is to try to relate these different outcomes to different parameter ranges in our model. The heuristic remarks of this paragraph will be backed by numerical simulations in Section~\ref{sec:numerical}.
We consider here a bounded domain $\Omega\subset \mathbb{R}^2$, and make some further hypotheses in order to simplify the system:
\begin{itemize}
 \item[i)] the inter-species interaction is symmetric $\Phi^{AB}=\Phi^{BA}$;
 \item[ii)] the interaction radius $R$ is much smaller than the scale of the experiments;
 \item[iii)] the diffusion is neglected.
 \end{itemize}
 Under hypothesis i), the macroscopic system of equation admits a free energy, or Lyapunov functional
\eqh{
 \mathcal{F}[f^A,f^B] =& \int_\Omega \left( \frac12 \tilde{\Phi}^{AA}(x-y) f^A(x)f^A(y)+\frac12\tilde{\Phi}^{BB}(x-y) f^B(x)f^B(y)+\tilde{\Phi}_{AB}(x-y) f^A(x)f^B(y)\right)dxdy \\
 & + \int_\Omega \left(D^A f^A(x) \ln f^A(x) +D^B f^B(x) \ln f^B(x)\right)dx.
}
The macroscopic dynamics is a gradient flow of the functional $\mathcal{F}$, with respect to the Wasserstein distance $W_2$ \cite{JKO98}. 
Owing to hypothesis ii), we will replace the potentials $\Phi^{\alpha\beta}$ by Dirac delta functions with  weights equal to their integrals, which we call 
$\gamma^{\alpha\beta}$, respectively.
Owing to hypothesis iii), we  neglect diffusion, which all together leads to the simplified functional $\mathcal{E}$:
 \begin{equation}
 \label{eq:Esimpl}
 \mathcal{E}[f^A,f^B] =  \int_\Omega \left( \frac12 \gamma^{AA} (f^A)^2(x) +\gamma^{AB}f^A(x)f^B(x)+\frac12 \gamma^{BB}(f^B)^2(x)\right)dx. 
 \end{equation}
Under hypotheses i)-iii) we expect that the system evolves in order to minimize \eqref{eq:Esimpl}, with constraints of constant total mass and positivity (we assume that at $t=0$, there is the same amount of $A$ and $B$ cells):
\begin{equation}
\label{eq:constraints}
\int_\Omega f^A(x)dx=\int_\Omega f^B(x) dx=\frac12~,~f^A(x)\geq 0~,~f^B(x)\geq 0.
\end{equation}
This problem is easily solved if $\gamma^{AA}<0$ or $\gamma^{BB}<0$. Then the minimum energy is $-\infty$ and it corresponds to one or both species infinitely concentrated. From now on we assume that  
$\gamma^{AA}>0$ and $\gamma^{BB}>0$, i.e. the intraspecies interactions are repulsive.
Another singular case corresponds to the range $\gamma^{AB}<-\sqrt{\gamma^{AA}\gamma^{BB}}$; then there is a pair of constants $(n^A>0,n^B>0)$ such that 
\[
\gamma^{AA}(n^A)^2 + 2\gamma^{AB}n^An^B +\gamma^{BB}(n^B)^2 <0.
\]
Now, building a sequence of configurations where $A$ and $B$ cells have concentration $n^A/\ep$, $n^B/\ep$ on the same domain of measure 
$\ep$, we see that the associated energy tends to $-\infty$ when $\ep$ tends to $0$.
Therefore, we also assume from now on that
$\gamma^{AB}> -\sqrt{\gamma^{AA}\gamma^{BB}}$. 

\begin{lemma}\label{lemma:opt1}
Assume that $\gamma^{AA}>0$, $\gamma^{BB}>0$, and  $-\sqrt{\gamma^{AA}\gamma^{BB}}<\gamma^{AB}<\sqrt{\gamma^{AA}\gamma^{BB}}$.
Then minimum energy is reached by the homogeneous system with constant densities $f^A(x)=f^B(x)=1/(2|\Omega|)$.
\end{lemma}
{\it Proof:} Since $\gamma^{AA}\gamma^{BB}-(\gamma^{AB})^2>0$, the function 
\[
(u,v)\mapsto \frac12 \gamma^{AA} u^2 +\gamma^{AB}uv+\frac12 \gamma^{BB}v^2
\]
is convex. Hence, for any functions $f^A(x),f^B(x)$:
\[
\frac{1}{|\Omega|}\mathcal{E}[f^A,f^B] \geq 
\frac12 \gamma^{AA} \left(\frac{\bar{f}^A}{|\Omega|}\right)^2 +\gamma^{AB}\left(\frac{\bar{f}^A}{|\Omega|}\right)\left(\frac{\bar{f}^B}{|\Omega|}\right)+\frac12 \gamma^{BB}\left(\frac{\bar{f}^B}{|\Omega|}\right)^2,
\]
where $|\Omega|$ is the volume of the domain and $\bar{f}^{S}=\int_\Omega f^{S}(x) dx=1/2$, $S=A,B$. 
This shows immediately that for any function $f^A(x),f^B(x)$
\[
\mathcal{E}[f^A,f^B] \geq \mathcal{E}[\bar{f}^A,\bar{f}^B]
=\frac{1}{8|\Omega|}(\gamma^{AA}+2\gamma^{AB}+\gamma^{BB}).
\]
The minimum energy is then reached for a homogeneous system, with constant densities, $f^A(x)=f^B(x)=1/(2|\Omega|)$.\hfill $\Box$
\medskip

We will now assume that $\gamma^{AB}>\sqrt{\gamma^{AA}\gamma^{BB}}$.
We want to show that the minimal configuration is perfectly phase separated, ie $\Omega$ divided in two subdomains $\Omega=\Omega^A \cup \Omega^B$ with  $f^{A}=0$ on $\Omega^{B}$ and $f^{B}=0$ on $\Omega^{A}$. Note that the geometry of the subdomains is of no relevance within this simplified model, only their measure matters.
\begin{lemma}\label{lemma:opt2}
Among all perfectly phase separated
configurations, the following is optimal:
\begin{eqnarray}
f^A(x)=n^A=const.~on~\Omega^A~&,&~f^B(x)=n^B=const.~on~\Omega^B \label{eq:opt} \\
with~|\Omega^A| = l_A^{\rm opt} &=& |\Omega|\left( \frac{\gamma_{AA}^{1/2}}{\gamma_{AA}^{1/2}+\gamma_{BB}^{1/2}}\right)  \label{eq:Linhom} 
\end{eqnarray}
The associated energy is
\begin{equation}
\mathcal{E}_{\rm inhom}^{\rm opt} =\frac{1}{8|\Omega|}\left(\sqrt{\gamma_{AA}}+\sqrt{\gamma_{BB}}\right)^2 \label{eq:Einhom}
\end{equation}
In other words: for any $(f^A,f^B)$ perfectly phase separated configuration,
\[
\mathcal{E}[f^A,f^B]\geq \mathcal{E}_{\rm inhom}^{\rm opt}
\]
with equality if and only if \eqref{eq:opt}-\eqref{eq:Linhom} hold. 
\end{lemma}
{\it Proof:} First note that the functions 
\[
u\mapsto \frac12 \gamma^{AA}u^2~,~ u\mapsto \frac12 \gamma^{BB}u^2~
\]
are convex; hence any perfectly phase separated configuration must be piecewise homogeneous in order to be optimal.
Let us now consider a general piecewise homogeneous phase separated configuration:
\begin{equation}
\forall x \in \Omega^A~f^A(x)=n^A~,~f^B(x)=0  ~,~\forall x \in \Omega^B~f^B(x)=n^B~,~f^A(x)=0 
\end{equation}
with
\[
|\Omega^A|=l_A~,~ |\Omega^B|=l_B~,~l_A+l_B = |\Omega|.
\]
Then $n^A=1/(2l_A)$, $n^B=1/(2|\Omega|-2l_A)$, and the configuration is characterized by the parameter $l_A$, which we can optimize. The associated energy is
\[
\mathcal{E} = \frac18 \lr{\frac{\gamma_{AA}}{l_A}+\frac{\gamma_{BB}}{|\Omega|-l_A}}.
\]
Minimizing over $l_A$, one finds 
\begin{eqnarray}
l_A^{\rm opt}&=&|\Omega| \left( \frac{\gamma_{AA}^{1/2}}{\gamma_{AA}^{1/2}+\gamma_{BB}^{1/2}}\right)  
\end{eqnarray}
with associated energy given by \eqref{eq:Einhom}. \hfill $\Box$
\medskip

\begin{lemma}\label{lemma:opt3}
Assume $\gamma^{AA}>0$ and $\gamma^{BB}>0$, and $\gamma^{AB}>\sqrt{\gamma^{AA}\gamma^{BB}}$. Then the optimal configuration is the perfectly phase separated one described in the previous lemma, ie
for any configuration $(f^A,f^B)$, we have
\[
\mathcal{E}[f^A,f^B]\geq \mathcal{E}_{\rm inhom}^{\rm opt}
\]
\end{lemma}
{\it Proof:}
Take any configuration $f^A,f^B$. 
Using $\int_\Omega f^{A,B} =1/2$ and \eqref{eq:Einhom}, we have
\eqh{
\mathcal{E}[f^A,f^B]-\mathcal{E}_{\rm inhom}^{\rm opt}
=& \frac12 \gamma^{AA}\int (f^A)^2 + \frac12 \gamma^{BB}\int (f^B)^2 +\gamma^{AB} \int f^Af^B  \\
& -\frac12{\frac{\gamma^{AA}}{|\Omega|} }\left(\int f^A\right)^2
-\frac12 {\frac{\gamma^{BB}}{|\Omega|} }\left(\int f^B\right)^2
-{\frac{\sqrt{\gamma^{AA}\gamma^{BB}}}{|\Omega|}}\left(\int f^A\right)\left(\int f^B\right) \\
=& \frac12 \int \left(\sqrt{\gamma^{AA}}f^A -{\frac{\sqrt{\gamma^{AA}}}{|\Omega|}}\int f^A \right)^2 + \frac12 \int \left(\sqrt{\gamma^{BB}}f^B -{\frac{\sqrt{\gamma^{BB}}}{|\Omega|}}\int f^B \right)^2 \\
&- \int \left( \sqrt{\gamma^{AA}}f^A-{\frac{\sqrt{\gamma^{AA}}}{|\Omega|}}\int f^A \right)\left(\sqrt{\gamma^{BB}}f^B -{\frac{\sqrt{\gamma^{BB}}}{|\Omega|}}\int f^B \right) \\ 
&+(\gamma^{AB}-\sqrt{\gamma^{AA}\gamma^{BB}}){\int f^Af^B}  \\
\geq& \frac12 \int \left[\left(\sqrt{\gamma^{AA}}f^A-\sqrt{\gamma^{AA}}\int f^A\right)-\left(\sqrt{\gamma^{BB}}f^B -\sqrt{\gamma^{BB}}\int f^B\right)\right]^2 
}
where we have used $(\gamma^{AB})> [\gamma^{AA}\gamma^{BB}]^{1/2}$ for the last inequality. Hence the inequality is strict unless $f^Af^B=0$, ie unless the configuration is phase separated.
We have proved that the optimal configuration is 
phase separated, and given by \eqref{eq:Linhom}. \hfill $\Box$
\medskip

To summarize, assuming $\gamma^{AA}>0$, $\gamma^{BB}>0$, $\gamma^{AB}>-\sqrt{\gamma^{AA}\gamma^{BB}}$  and putting together Lemmas  \ref{lemma:opt1},\ref{lemma:opt2},\ref{lemma:opt3}, one expects the following phenomenology:
\begin{itemize}
\item If $|\gamma^{AB}|<\gamma_{AA}\gamma_{BB}$, then the homogeneous solution is favoured. The two types of cells should then start to mix and the interface should be blurred. 
\item If $\gamma^{AB}>\sqrt{\gamma_{AA}\gamma_{BB}}$, then the phase separated solution is favoured; in this case the initially sharp interface should stay sharp. There are two subcases:
\begin{itemize}
\item $\gamma_{AA} \simeq \gamma_{BB}$: then the optimal $l$ is close to $1/2$, which is the initial condition. Hence the interface should not move.
\item $\gamma_{AA}$ and $\gamma_{BB}$ significantly different: then the optimal $l$ is not close to the initial $1/2$, and one expects the sharp interface to move, as the system tries to approach the energy minimum.
\end{itemize}
\end{itemize}
These three scenarios are qualitatively similar to the ones reported in \cite{Taylor2017} (see their Figure~1), and are seen in numerical simulations, as shown on Fig.\ref{fig:step}; a quantitative comparison is difficult, since hypothesis i)-iii) are not necessarily satisfied in numerical simulations.

\section{Numerical results} \label{sec:numerical}
% \subsection{Averaged microscopic model}\label{sec:averaged-micro}
The microscopic model described in Section \ref{sec:IBM} is very demanding numerically; it is one reason to introduce the macroscopic model  \eqref{Eqfmacro}, which relies on the double limit $(N_A,N_B) \to \infty$, $\varepsilon \to 0$. The averaged microscopic model \eqref{eq:averagedA}-\eqref{eq:averagedB}, obtained in the limit $\varepsilon \to 0$ with $N_A,N_B$ fixed, can be simulated at a reasonable numerical cost, which makes comparisons with the macroscopic model possible. 

\medskip

Numerical simulations for the averaged microscopic   model \eqref{eq:averagedA}-\eqref{eq:averagedB} as well as macroscopic model \eqref{Eqfmacro} are performed on a 2D periodic domain $[-7.5, \; 7.5] \times [ -7.5, \; 7.5]$. All simulations are performed  with cell detection radii $R_A=R_B = 1$. We fix the diffusion constants $D_A=D_B = 10^{-4}$ and explore different values of  the inter- and intra- species repulsion intensities $\kappa^{AA}, \kappa^{BB}, \kappa^{AB} = s \tilde{\kappa}^{AB}, \kappa^{BA} = s \tilde{\kappa}^{BA}$. We explore four different cases: 
\begin{itemize}
\item case 1: For the same intraspecies repulsion  $\kappa^{AA} = \kappa^{BB}$ and symmetric inter-species repulsion $\tilde{\kappa}^{AB} = \tilde{\kappa}^{BA}$
\item case 2: For the same intraspecies repulsion  $\kappa^{AA} = \kappa^{BB}$ and non-symmetric inter-species repulsion $\tilde{\kappa}^{AB} < \tilde{\kappa}^{BA}$ ($A$-cells repulse  $B$-cells more strongly than the reverse)
\item case 3: For different intraspecies repulsion $\kappa^{AA} > \kappa^{BB}$ ($A$-cells repulse  each other more strongly  than $B$-cells), and symmetric inter-species repulsion $\tilde{\kappa}^{AB} = \tilde{\kappa}^{BA}$
\item case 4: For different intraspecies repulsion $\kappa^{AA} > \kappa^{BB}$ ($A$-cells repulse  each other more strongly  than $B$-cells) and non-symmetric inter-species repulsion $\tilde{\kappa}^{AB} < \tilde{\kappa}^{BA}$ ($A$-cells repulse  $B$-cells more strongly than the reverse)
\end{itemize}
For each case, we consider two types of initial conditions, (i) when cells are initially randomly distributed, which approaches the homogeneous stationary state, and (ii) when $B$-cells are initially randomly placed on the left-half of the domain and $A$-cells randomly distributed on the right part, which corresponds to the phase separated initial conditions considered in Section \ref{sec:separatedphases}, and we explore two different regimes (stable regime, when $s<s^*$ for the interspecies repulsion and unstable regime, when $s>s^*$). Table \ref{Table1} sums up the model parameters used for inter- and intra- species repulsion forces. 

\begin{figure}[H]
\includegraphics[scale=0.4]{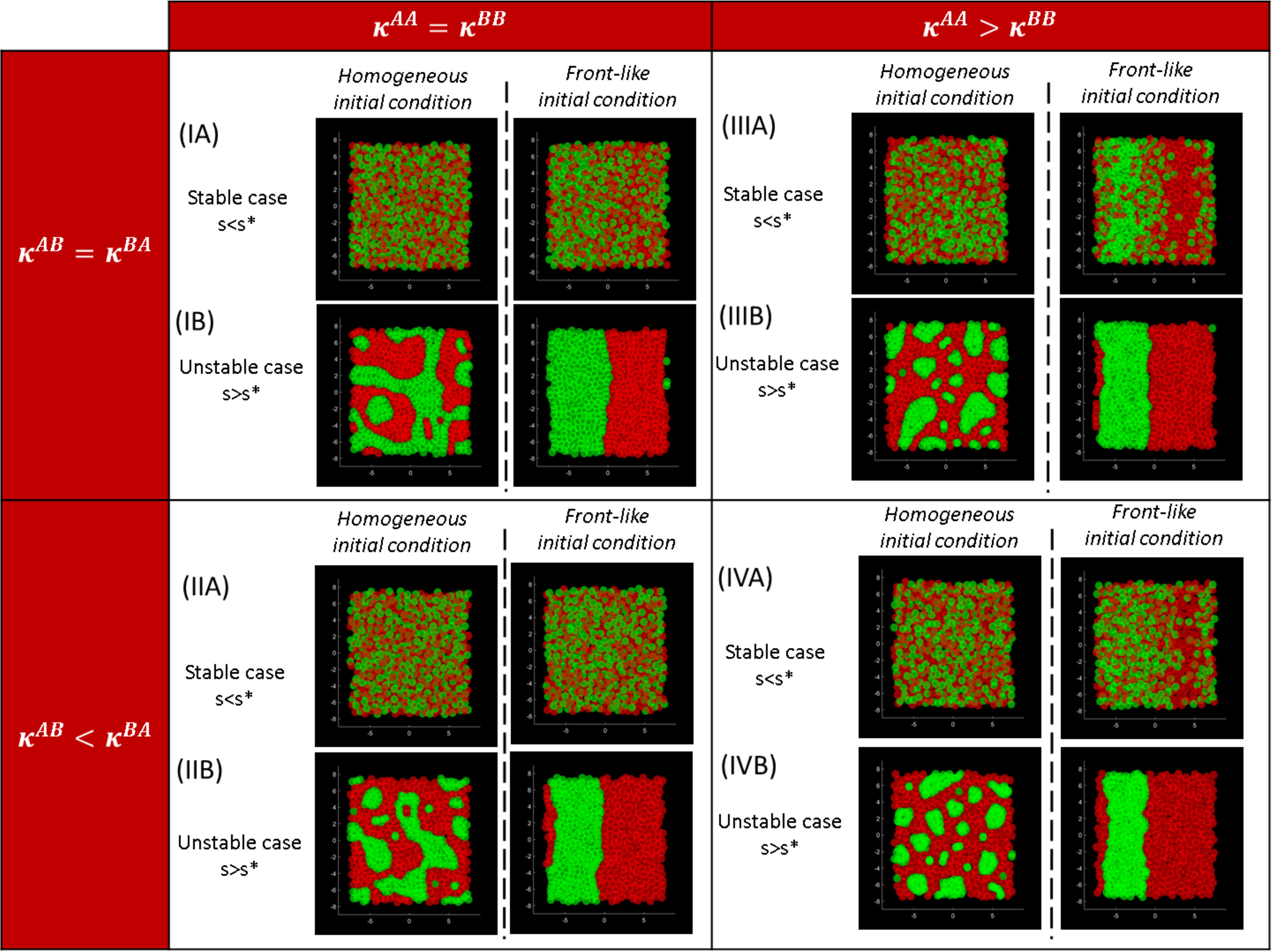}
\caption{Microscopic simulations for Cases 1-4 for parameters described in Table \ref{Table1}. $A$-cells are represented as red disks, $B$-cells as green disks. For each subsection, the left figure is obtained starting from a homogeneous distribution of particles, the right one from a segregated initial distribution ($B$-cells on the half-left of the domain, $A$-cells on the right).\label{SimusFig1}}
\end{figure}

\begin{table}[H]
\begin{tabular}{ccc}
 &value of $s$ & comment\\
\hline
\multicolumn{3}{c}{Case I: $\kappa^{AA}=\kappa^{BB}= {2}$, $\tilde{\kappa}^{AB} = \tilde{\kappa}^{BA}= {2}$.}\\
\hline
IA & 0.5 & Stable regime ($s < s^* \approx 1.01$)\\
IB & 4 & Unstable regime ($s>s^* \approx 1.01$)\\
\hline
\multicolumn{3}{c}{Case II: $\kappa^{AA}=\kappa^{BB}= 2$, $1=\tilde{\kappa}^{AB} < \tilde{\kappa}^{BA}= 2$.}\\
\hline
IIA & 0.5 & Stable regime ($s < s^* \approx 1.43$)\\
IIB & 4 & Unstable regime ($s>s^* \approx 1.43$)\\
\hline
\multicolumn{3}{c}{Case III: $2=\kappa^{AA}>\kappa^{BB}= 1$, $\tilde{\kappa}^{AB} = \tilde{\kappa}^{BA}=2$.}\\
\hline
IIIA & 0.5 & Stable regime ($s < s^* \approx 0.72$)\\
IIIB & 4 & Unstable regime ($s>s^* \approx 0.72$)\\
\hline
\multicolumn{3}{c}{Case IV: $2=\kappa^{AA}>\kappa^{BB}= 1$, $1=\tilde{\kappa}^{AB} < \tilde{\kappa}^{BA}= 2$.}\\
\hline
IVA & 0.5 & Stable regime ($s < s^* \approx 1.02$)\\
IVB & 4 & Unstable regime ($s>s^* \approx 1.02$)
\end{tabular}
\caption{Model parameters for the inter- and intra- species forces. The value of parameter $s^*$ has been computed numerically from the formula  \eqref{s_star_per}. \label{Table1}}
\end{table}
In Fig. \ref{SimusFig1}, we show the final states of the simulations for each case and each regime previously described for the microscopic model. $A$-cells are represented as red disks, $B$-cells as green ones and we use $N_A=N_B=250$ particles for each family of cells. The visualisation of the macroscopic results also uses disks to resemble the microscopic ones, for more details we refer to the Appendix \ref{AppendixC}.

 As one can note in the stable regime ($s<s^*$) and for initially randomly distributed particles, we observe a homogeneous distribution of particles as expected, with no aggregation. When starting from a front-like intial distribution in the stable regime, $B$- and $A$- cells intermingle at the front. However in the unstable regime (cases $s>s^*$), one can observe a segregation of cells by type. In case 1 and 2 (when intraspecies repulsion is the same), we observe the formation of mazes of $B$-cells when starting from an initial homogeneous distribution, and the maintenance  of a sharp front when starting from a non-homogeneous initial distribution (Fig. \ref{SimusFig1} (IB, IIB)). Note that the $B$-cell clusters are smaller in case 2 than in case 1 (i.e when $A$-cells act more strongly on $B$-cells than the reverse compared with the case where inter-species repulsion is symmetric). When starting from a segregated initial condition we observe a slight compression of the green cells by the red ones. This suggests that nonsymmetric interspecies repulsion can favor cell aggregation and domination of a population over the other.
 
The same observations as in case 2 can be done for case 3 (i.e when intra-species repulsion is stronger in  cell $A$ type than in  $B$ type and with symmetric interspecies forces, (IIIB)): smaller clusters than in case 1 and slight compression of $B$-family in the case of segregated initial condition. This suggests that inter- and intra- specie repulsion act in the same manner: decreasing the intra-species force for one family has the same impact on the final structures as decreasing the repulsion force of one family onto the other one.

Finally for case 4 (intra-species repulsion is stronger in  cell $A$ type than in  $B$ type and stronger repulsion of $B$-cells by $A$-cells than the reverse), we observe the emergence of small clusters of $B$-cells in a medium composed of $A$-cells, and a large compression of the $B$-cells with maintenance of sharp borders when starting from a non-homogeneous distribution (IV B).

It is noteworthy that the numerical simulations of the microscopic model in the limit $\epsilon \rightarrow 0$ are in good accordance with the predictions of the macroscopic model. We indeed observe  homogeneous or non-homogeneous distributions of particles depending on the model parameters, for values of $s$ in the range predicted by the linear stability analysis. 

\begin{figure}[H]
\includegraphics[trim=50 0 0 0, scale=0.6]{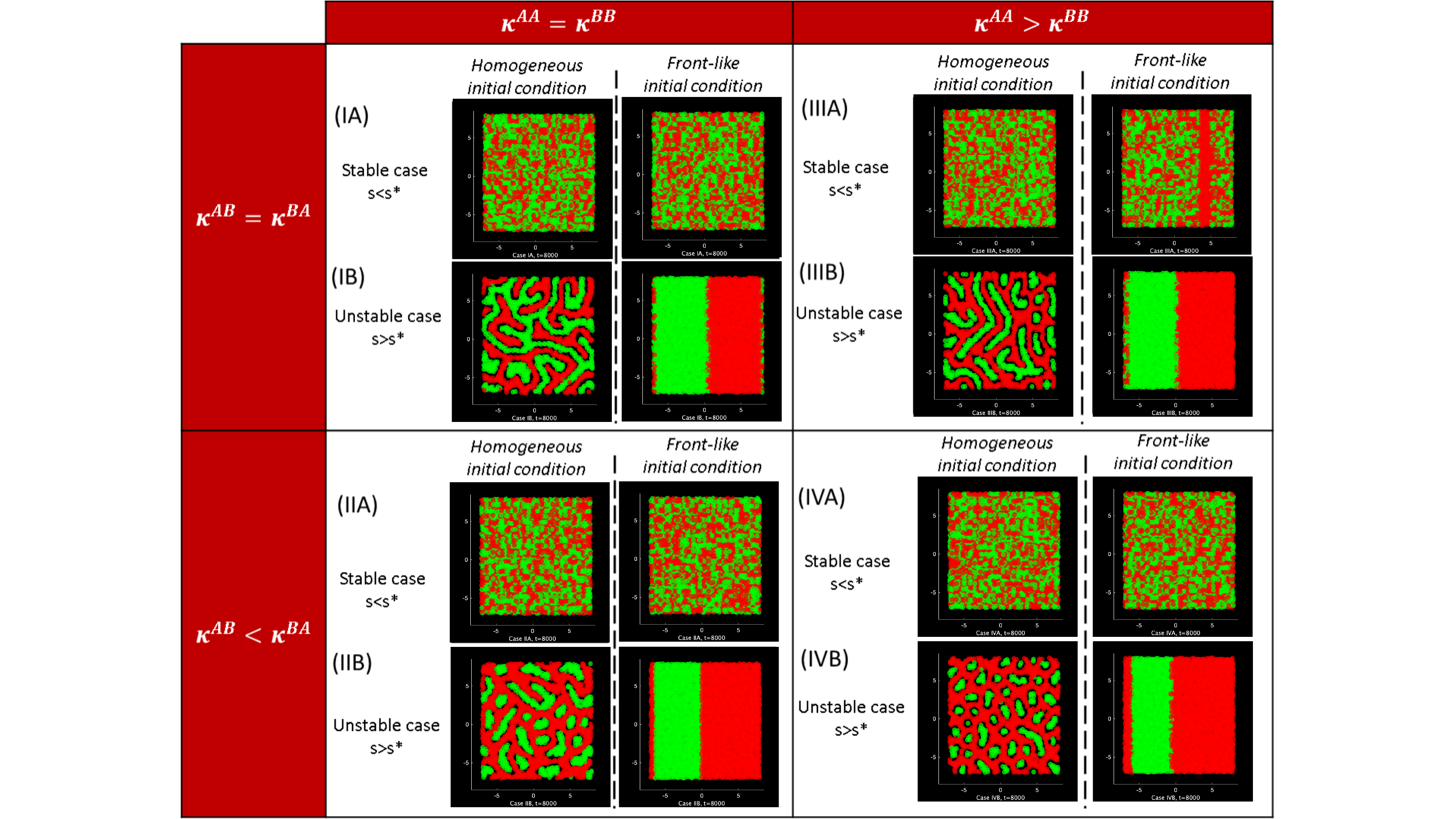}
\caption{Macroscopic simulations for Cases 1-4 for parameters described in table \ref{Table1} for the final time of simulations equal to $T=8000$.\label{fig:step}}
\end{figure}

\begin{table}[]
    \centering
    \begin{tabular}{c|c|c|c|c}
    Case number  &  \multicolumn{2}{c|}{VF green cells (microscopic model)} & \multicolumn{2}{c}{VF green cells  (macroscopic model)}\\
    \hline
    & Homogeneous IC & Front-like IC & Homogeneous IC & Front-like IC\\
    case (IB) & $48.2 \%$ & $50.0 \%$ & $49.7\%$& $50.0\%$\\
    case (IIB) & $40.2 \%$ & $42.3 \%$ & $38.5\%$ &$42.0\%$\\
    case (IIIB) & $40.6 \%$ & $42.4 \%$ & $44.0\%$ &$46.0\%$\\
    case (IVB) & $34.8 \%$ & $35.0 \%$ &$35.9\%$ &$38.0\%$\\
    \end{tabular}
    \caption{Volume fraction of the green family computed on the simulations of FigS. \ref{SimusFig1}-\ref{fig:step} at equilibrium for the microscopic model (left column) and for the macroscopic model (right column).}
    \label{table2}
\end{table}

In Table \ref{table2}, we show the Volume Fraction (VF) of type B cells (green family) computed on the simulation images of Figs. \ref{SimusFig1} for the microscopic model and of Fig. \ref{fig:step} for the macroscopic model, starting from homogeneous Initial Conditions (IC) or front-like Initial Conditions. Given a numerical image such as depicted in Fig. \ref{SimusFig1}, the volume fraction corresponds to the number of green pixels over the total amount of pixels in the image. As one can see in Table \ref{table2}, the volume fraction at equilibrium does not depend on the type of initial conditions, as suggested by the analysis performed in Section \ref{sec:separatedphases}. As expected, we obtain a volume fraction of $50 \%$ when the two cell types have the same inter- and intra- species forces, and the volume fraction occupied by specie B decreases as type $B$ cells' inter- and/or intra- species forces decrease (maintaining the type $A$ cells inter- and intra- species forces constant), due to the compression exerted by the stronger family (type $A$) on the weaker cells (type $B$). These results are in accordance with the theoretical predictions of the macroscopic model, showing that the microscopic and macroscopic model have the same properties. The next section is devoted to deeper numerical comparisons between the two models.

% \subsection{Comparison of the Microscopic and Macroscopic models}\label{sec:4.2}

% \red{Should we move the derivation of image processing to the appendix?}\blue{\emph{I vote yes.}}

\begin{figure}
    \centering
    \includegraphics[scale = 0.25]{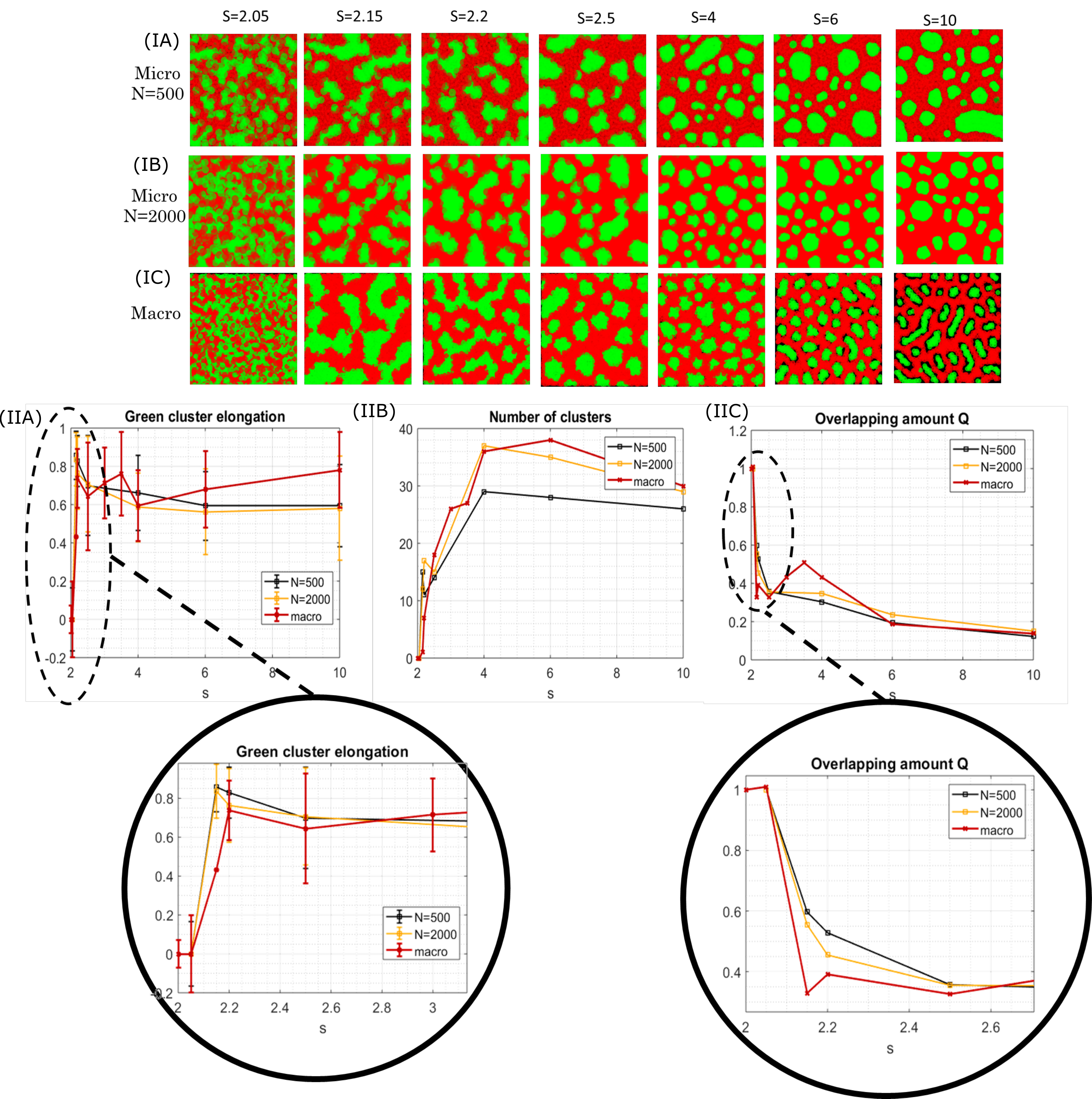}
    \caption{(I) Simulations of the microscopic and macroscopic models for $\kappa^{AA} = 4, \kappa^{BB} = \tilde{\kappa}^{AB} = \tilde{\kappa}^{BA}=1$ for which  $s^* \approx 2.1$. Simulations of the microscopic model are performed with $N_A = N_B =500$ (IA) and $N_A=N_B=2000$ (IB) particles. Simulations of the macroscopic model correspond to (IC). We consider 7 values of the interspecies repulsion intensity, from left to right: for $2.05=s<s^*$, for $s^*<s = \{2.15, 2.2, 2.5, 4, 6, 10 \}$. Type B cells are represented in green, type A cells in red. (II) In (IIA-C), we show the values of the quantifiers at time equilibrium as functions of $s$. Fig. (IA) shows the mean elongation of the green clusters, (IIB) shows the number of green clusters and (IIC) shows the overlapping amount $Q$ described by \eqref{Q}. Black curves are obtained with the microscopic model for $N_A=N_B=500$ (corresponding to Figures (IA)), yellow curves are for $N_A=N_B=2000$ and correspond to Figures (IB) and red curves are obtained with the macroscopic model (Figures (IC)). The two bottom figures correspond to a zoom of the corresponding curves close to the transition region for $s$. }
    \label{fig:kBB1}
\end{figure}

In Figs. \ref{fig:kBB1} (I), we show simulations of the microscopic and macroscopic models for $\kappa^{AA} = 4, \kappa^{BB} = \tilde{\kappa}^{AB} = \tilde{\kappa}^{BA}=1$ for which $s^* \approx 2.1$ can be computed using the formula \eqref{s_star_per}. Simulations of the microscopic model are performed with $N_A = N_B =500$ (IA) and $N_A=N_B=2000$ (IB) particles. Simulations of the macroscopic model correspond to (IC). We consider 7 values of the interspecies repulsion intensity, from left to right: for $2.05=s<s^*$, for $s^*<s = \{2.15, 2.2, 2.5, 4, 6, 10 \}$. (II) In (IIA-C), we show the values of the quantifiers at time equilibrium as functions of $s$. Fig. (IA) shows the mean elongation of the green clusters, (IIB) shows the number of green clusters and (IIC) shows the Overlapping amount $Q$ described by \eqref{Q}. Black curves are obtained with the microscopic model for $N_A=N_B=500$ (corresponding to Figures (IA)), yellow curves are for $N_A=N_B=2000$ and correspond to Figures (IB) and red curves are obtained with the macroscopic model (Figures (IC)).
As one can observe in Fig. \ref{fig:kBB1} (IA-C), we obtain a very good agreement between the microscopic model and the macroscopic simulations. Before the transition from mixed to segregated states (for $s<s^*$, first column), the system at equilibrium corresponds to a perfectly mixed state both for the micro- and for the macro- models, while right after the transition (for $s=2.15$), segregation of the two families is observed for both models. As $s$ increases, clusters get more numerous, smaller and rounder. This qualitative observation is supported by the values of the quantifiers plotted in Figs. (IIA-C). For the different values of $s$, we obtain a very good quantitative agreement between the two models, even more so when the number of particles of the micro- model is increased from $N_A=N_B=500$ to $N_A=N_B=2000$ (compare black and yellow curves to red ones in Figs. (IIA-C)). These results tend to show that the macroscopic model is 
a good approximation of the microscopic dynamics as the number of individuals becomes large. However, it is noteworthy that some discrepancy is observed for very large values of $s$. Indeed for $s=10$ (last column of (IA-C)), one can note that the clusters obtained by the macro- model are significantly more elongated than those obtained with the micro- model (compare black/yellow curves to the red one in Fig.\ref{fig:kBB1} (IIA)). For the details of image processing used to 
prepare the figures we refer to the Appendix \ref{AppendixC2}.

\begin{figure}
    \centering
    \includegraphics[scale = 0.45]{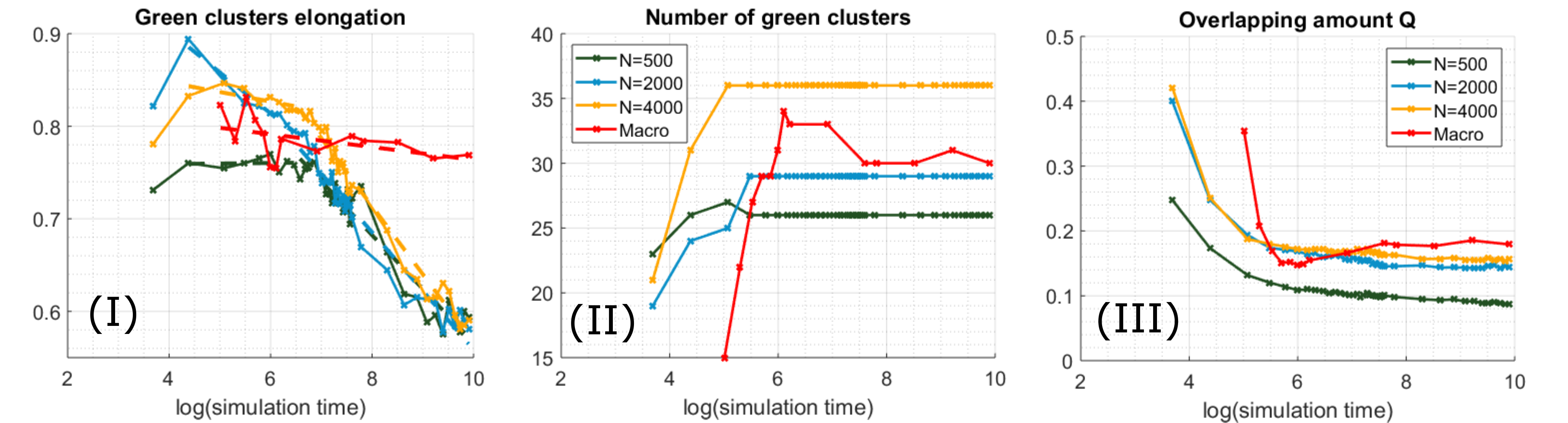}
    \caption{Quantifiers computed on the simulation images as functions of the logarithm of the simulation time for $\kappa^{AA}=4, \kappa^{BB} = \tilde{\kappa}^{AB} = \tilde{\kappa}^{BA} = 1$ for the microscopic model with $N_A=N_B=500$ (green curves), $N_A=N_B=2000$ (blue curves), $N_A=N_B = 4000$ (yellow curves) and for the macro model (red curve). (I) Green cluster elongation, (II) Number of cell clusters and (III) Overlapping amount $Q$. On Figure (I), we superimpose linear fits (dotted lines) for short times and large times, showing the two timescales (two slopes) of the micro model compared to the unique timescale of the macro dynamics (single slope).}
    \label{fig:kBB1_time}
\end{figure}

In order to document the discrepancies between the micro- and macro- structures in the case $s = 10$, we plot in Fig. \ref{fig:kBB1_time} the values of the quantifiers computed on the simulation images as functions of time for $\kappa^{AA}=4, \kappa^{BB} = \tilde{\kappa}^{AB} = \tilde{\kappa}^{BA} = 1$ for the microscopic model with $N_A=N_B=500$ (green curves), $N_A=N_B=2000$ (blue curves), $N_A=N_B = 4000$ (yellow curves) and for the macro model (red curves). Fig. \ref{fig:kBB1_time} (I) shows the evolution of green cluster elongation, (II) gives the number of cell clusters and (III) shows the overlapping amount $Q$ as function of the simulation time. 
Figs. \ref{fig:kBB1_time} reveal that the structures observed with the microscopic model undergo two transitions as function of the time, enabling us to conclude that the dynamics of the microscopic model exhibits two phases:
\begin{itemize}
    \item Very quickly after initialisation ($t \in [0,100]$), the system segregates the two cell types and reaches an equilibrium value for the number of clusters and the overlapping amount $Q$ (see Figs. \ref{fig:kBB1_time} (II) and (III). At the time of the segregation, numerous and well-separated elongated clusters are created and then maintained for a long time. 
    \item A second transition occurs later in time (around $t\approx 3.10^{3}$), where the clusters change shape to attain a new equilibrium composed of rounder clusters (see the drop in the value of the elongation in Fig. \ref{fig:kBB1_time} (I)). In this second phase, the number of clusters and their border properties are maintained, but the shape gradually changes to produce very round clusters.
\end{itemize}

It is noteworthy that this two-phase process is not observed with the macroscopic dynamics (see the red curves of Figs. \ref{fig:kBB1_time} (I-III). On the contrary for the macroscopic model, the segregation between the two families and the production of clusters appear later than with the microscopic model (around $t \in [100, 300]$), and the shape of the clusters (elongated) seems to be at equilibrium. These results tend to show that the macroscopic model fails to capture the second time phase (reorganisation of the clusters) exhibited by the microscopic model.

The  good qualitative and quantitative agreement between the micro- and the macro- models has been also confirmed for another set of parameters $\kappa^{AA}=4, \kappa^{BB} = 3, \tilde{\kappa}^{AB} = \tilde{\kappa}^{BA} = 1$ for which the critical value of $s$ is $s^* \approx 3.52$. We do not include the corresponding figures here, as they are very similar to Figs. \ref{fig:kBB1_time}.

\section{Conclusions and perspectives}

Along the recent biological studies \cite{Taylor2012} our paper demonstrates that contact cell repulsion on its own can generate pattern formation and cell-sorting in tissues composed of different categories of cells. The present paper provides evidence of this by means of both a microscopic and a consistently derived macroscopic model. The advantage of the macroscopic approach is that it provides a mathematical way to investigate the stability of the equilibria and consequently quantitative criteria for the appearance of these patterns. The validity of this analysis is assessed by numerical comparison between microscopic and macroscopic models. We also show that the model is able to capture the border sharpening observed in the biological study \cite{Taylor2012}. In the future better quantitative comparison with experiments will allow for more systematic choice of model parameters. 

Clearly the macroscopic model captures the behaviour of the mean, and it would be interesting to investigate what fluctuations around these means are induced by the finiteness of the number of cells.  Our analysis provides only local information around the threshold for instability, further analysis as in \cite{Barre2016} could inform on the type of bifurcation involved, for instance if it is supercritical or subcritical. The consistent derivation of the macroscopic model from the microscopic dynamics is still only formal and it would be desirable to have a rigorous mathematical proof of convergence when the number of particles goes to infinity. Finally, the biological context could be enriched and applied to clinically relevant situations such as cancer for which the type of bifurcation involved could be critically important.

\section*{Acknowledgements} The authors wish to thank Anais Khuong for introducing  to the biological aspects of the paper. PD acknowledges support from the Engineering and Physical Sciences Research Council (EPSRC) under grant ref. EP/M006883/1 and EP/N014529/1, from the Royal Society and the Wolfson foundation through a Royal Society Wolfson Research Merit Award no. WM130048 and from the National Science Foundation (NSF) under grant RNMS11-07444 (KI-Net). PD is on leave from CNRS, Institut de Math\'ematiques, Toulouse, France.
EZ was supported by the Polish Government MNiSW research grant 2016-2019 "Iuventus Plus"  No.  0888/IP3/2016/74.

\section*{Data statement}
No new data were generated in the course of this research

% \red{Add conclusions and perspectives: major take-away message, future directions}

% \blue{Conclusions: more or less everything is in the introduction?}

% \blue{Some ideas for future perspectives:
% \begin{itemize}
%     \item Quantitative comparison with experiments;
%     \item Understand better the complex pattern formation, and the finite $N$ effects we have identified;
%     \item Prove the convergence towards the macroscopic models in some cases;
%     \item Study the bifurcation of the homogeneous steady state (subcritical?), possibly also using a variational
%     approach (see Chayes-Panferov in the one-specie case).
%     \end{itemize}
% }

\newpage
\appendix
%----------------------------APPENDIX I-----------------
%---------------------derivation micro-macro------------
%\setcounter{equation}{0}
\section{Derivation of a kinetic model from the microscopic model}
\label{AppendixA}

In this section, the derivation of a kinetic model from the Individual Based Model of Section~\ref{sec:IBM} is performed, following the Approach I. Using the expressions for the individual particle distribution, link distribution and two-particle distribution  function defined by Eqs. \eqref{fAB},\eqref{gAA},\eqref{gAB}, \eqref{TPDFAA}-\eqref{TPDFAB}
 and in the limit of a large number of individuals, we have the following formal derivation:\\

\begin{prop}\label{thm1}
Assume that in the limit 
$$N=(N_A, N_B) \rightarrow \infty, \;N_A/N_B\to r_{AB}>0,$$ 
the following convergences hold:
$$f^{A}_{N},f^{B}_{N}  \rightarrow f^{A}, f^B, \quad g^{AA}_{N},g^{BB}_N,g^{AB}_N,g^{BA}_N \rightarrow g^{AA},g^{BB},g^{AB},g^{BA},$$ 
$$h^{AA}{N},h^{BB}_N,h^{AB}_N \rightarrow h^{AA},h^{BB},h^{AB}.$$
If we assume the scalings \eqref{eq:nuscalings} for the rates, and if we assume additional assumptions (\eqref{eq:MF1},\eqref{eq:MF2}  and \eqref{eq:MF1b},\eqref{eq:MF2b} below), then $f^{A}, f^B$ formally solve:
\begin{equation}\label{Eqfkinet}
\begin{cases}
 \partial_t f^A(x,t) = 2\mu  \nabla_x \cdot F^{AA}[g^{AA}](x,t) + \mu  \nabla_x \cdot F^{AB}[g^{AB}](x,t) + D^A\Delta f^A,\\
 \partial_t  f^B (x,t) = 2\mu  \nabla_x \cdot F^{BB}[g^{BB}](x,t) + \mu  \nabla_x \cdot F^{BA}[g^{BA}](x,t) + D^B\Delta f^B,
 \end{cases}
 \end{equation}
 \noindent where:
 \begin{align*}
  F^{AA}[g](x,t) &= \int g(x_1,x_2,t) \nabla_{x_1} \Phi^{AA}(x_1,x_2) dx_2,\quad
      F^{BB}[g](x,t) &= \int g(x_1,x_2,t) \nabla_{x_1} \Phi^{BB}(x_1,x_2) dx_2,\\
    F^{AB}[g](x,t) &= \int g(x_1,x_2,t) \nabla_{x_1} \Phi^{AB}(x_1,x_2) dx_2,\quad
    F^{BA}[g](x,t) &= \int g(x_1,x_2,t) \nabla_{x_1} \Phi^{BA}(x_1,x_2) dx_2,
 \end{align*}
 \noindent and $g^{AA}, g^{BB}$ and $g^{AB}$ formally solve:
 \begin{align}
\partial_t g^{AA}&(x_1,x_2,t) = D^A \big(\Delta_{x_1} g^{AA}(x_1,x_2,t) +\Delta_{x_2} g^{AA}(x_1,x_2,t) \big) \nonumber\\
&+ 2\mu \nabla_{x_1} \cdot \bigg(\frac{g^{AA}(x_1,x_2,t)}{f^A(x_1,t)}F^{AA}[g^{AA}](x_1,t)\bigg)+ \mu  \nabla_{x_1} \cdot \bigg(\frac{g^{AA}(x_1,x_2,t)}{f^A(x_1,t)}F^{AB}[g^{AB}](x_1,t)\bigg)\label{EqgAAkinet}\\
&+ 2\mu  \nabla_{x_2} \cdot \bigg(\frac{g^{AA}(x_1,x_2,t)}{f^A(x_2,t)}F^{AA}[g^{AA}](x_2,t)\bigg)+ \mu  \nabla_{x_2} \cdot \bigg(\frac{g^{AA}(x_1,x_2,t)}{f^A(x_2,t)}F^{AB}[g^{AB}](x_2,t)\bigg) \nonumber \\
&+ \frac{\nu^{AA}_{c,\ep} }{2} h^{AA}(x_1, x_2, t)\chi(|x_1-x_2|\leq R)
- \nu^{AA}_{d,\ep} g^{AA}(x_1,x_2,t),\nonumber
\end{align}
\begin{align}
\partial_t g^{BB}&(x_1,x_2,t) =  D^B\big(\Delta_{x_1} g^{BB}(x_1,x_2,t) +\Delta_{x_2} g^{BB}(x_1,x_2,t) \big) \nonumber\\
&+ 2\mu \nabla_{x_1} \cdot \bigg(\frac{g^{BB}(x_1,x_2,t)}{f^B(x_1,t)}F^{BB}[g^{BB}](x_1,t)\bigg)+ \mu \nabla_{x_1} \cdot \bigg(\frac{g^{BB}(x_1,x_2,t)}{f^B(x_1,t)}F^{BA}[g^{BA}](x_1,t)\bigg)\nonumber\\
&+ 2\mu \nabla_{x_2} \cdot \bigg(\frac{g^{BB}(x_1,x_2,t)}{f^B(x_2,t)}F^{BB}[g^{BB}](x_2,t)\bigg)+ \mu \nabla_{x_2} \cdot \bigg(\frac{g^{BB}(x_1,x_2,t)}{f^B(x_2,t)}F^{BA}[g^{BA}](x_2,t)\bigg) \nonumber \\
&+ \frac{\nu^{BB}_{c,\ep} }{2} h^{BB}(x_1, x_2, t)\chi(|x_1-x_2|\leq R)
- \nu^{BB}_{d,\ep} g^{BB}(x_1,x_2,t),\nonumber
\end{align}
\begin{align}
\partial_t g^{AB}&(x_1,x_2,t) =  \big(D^A \Delta_{x_1} g^{AB}(x_1,x_2,t) +D^B\Delta_{x_2} g^{AB}(x_1,x_2,t) \big) \nonumber\\
&+ 2\mu  \nabla_{x_1} \cdot \bigg(\frac{g^{AB}(x_1,x_2,t)}{f^A(x_1,t)}F^{AA}[g^{AA}](x_1,t)\bigg)+ \mu  \nabla_{x_1} \cdot \bigg(\frac{g^{AB}(x_1,x_2,t)}{f^A(x_1,t)}F^{AB}[g^{AB}](x_1,t)\bigg)\nonumber\\
&+ 2\mu \nabla_{x_2} \cdot \bigg(\frac{g^{AB}(x_1,x_2,t)}{f^B(x_2,t)}F^{BB}[g^{BB}](x_2,t)\bigg)+ \mu \nabla_{x_2} \cdot \bigg(\frac{g^{AB}(x_1,x_2,t)}{f^B(x_2,t)}F^{BA}[g^{BA}](x_2,t)\bigg) \nonumber \\
&+ \nu^{AB}_{c,\ep}  h^{AB}(x_1, x_2, t)\chi(|x_1-x_2|\leq R)
- \nu^{AB}_{d,\ep} g^{AB}(x_1,x_2,t),\nonumber \\
 g^{BA}&(x_1,x_2,t) = r_{AB} g^{AB}(x_2,x_1,t).\nonumber
\end{align}
 \end{prop}

\begin{pf}
 \subsection{Evolution equation for the individual particles} 
 For all observable functions $\phi(x)$, we define:
\begin{align*}
 \langle f^A_N,\phi \rangle  &= \int\limits \phi(x) f^A_N(t,x)  dx = \frac{1}{N_A} \sum_{i=1}^{N_A} \phi(X^A_i(t)),\\
  \langle f^B_N,\phi \rangle  &= \int\limits \phi(x) f^B_N(t,x)  dx = \frac{1}{N_B} \sum_{\ell=1}^{N_B} \phi(X^B_\ell(t)).   
\end{align*}
\noindent Similarly, for all two-particle observable functions $\psi (x_1,x_2)$, we define: 
\begin{align*}
\langle \hspace{-0.8mm} \langle g^{AA}_N,\psi \rangle \hspace{-0.8mm} \rangle &= \int \psi(x_1,x_2) g^{AA}_N(x_1,x_2) dx_1dx_2
= \frac{1}{2N_A} \sum_{k_1=1}^{K_{AA}} \bigg(\psi(X^A_{i(k_1)}, X^A_{j(k_1)}) + \psi(X^A_{j(k)}, X^A_{i(k)})  \bigg),\\
\langle \hspace{-0.8mm} \langle g^{BB}_N,\psi \rangle \hspace{-0.8mm} \rangle &= \int \psi(x_1,x_2) g^{BB}_N(x_1,x_2) dx_1dx_2
= \frac{1}{2N_{B}} \sum_{k_2=1}^{K_{BB}} \bigg(\psi(X^B_{\ell(k_2)}, X^B_{m(k_2)}) + \psi(X^B_{\ell(k)}, X^B_{m(k)})  \bigg),\\
\langle \hspace{-0.8mm} \langle g^{AB}_N,\psi \rangle \hspace{-0.8mm} \rangle &= \int \psi(x_1,x_2) g^{AB}_N(x_1,x_2) dx_1dx_2
= \frac{1}{N_{A}} \sum_{k_3=1}^{K_{AB}} \psi(X^A_{i(k_3)}, X^B_{\ell(k_3)}), \\
\langle \hspace{-0.8mm} \langle g^{BA}_N,\psi \rangle \hspace{-0.8mm} \rangle &= \int \psi(x_1,x_2) g^{BA}_N(x_1,x_2) dx_1dx_2
= \frac{1}{N_{B}} \sum_{k_3=1}^{K_{AB}} \psi(X^B_{\ell(k_3)},X^A_{i(k_3)}),%, \\
%\langle \hspace{-0.8mm} \langle g_B^{K_{AB}},\psi \rangle \hspace{-0.8mm} \rangle &= \int \psi(x_1,x_2) g_B^{K_{AB}}(x_1,x_2) dx_1dx_2
%= \frac{1}{K_{BA}} \sum_{k_3=1}^{K_{AB}} \psi(X^B_{\ell(k_3)}, X^A_{i(k_3)}).
\end{align*}
\noindent where integrals over $x$ are carried over ${\mathbb R}^2$. Then: 
\begin{align*}
 \frac{d}{dt}\langle f^A_N,\phi \rangle  =  \frac{1}{N_A} \sum_{i=1}^{N_A} \nabla_x \phi(X_i(t)) \cdot \frac{d X^A_i(t)}{dt} , \qquad \frac{d}{dt}\langle f^B_N,\phi \rangle  =  \frac{1}{N_B} \sum_{i=1}^{N_B} \nabla_x \phi(X_i(t)) \cdot \frac{d X^B_i(t)}{dt}.
\end{align*}
\noindent For the sake of simplicity, the computations for $f^A_N$ only are developed here. Using~\eqref{IBMA} and It\^o's formula, we obtain formally: 
\begin{align}\label{two_diff}
 \frac{d}{dt}\langle f^A_N,\phi \rangle
 =  -\frac{1}{N_A} \sum_{i=1}^{N_A}  \mu &\nabla_x \phi(X^A_i)  \cdot \nabla_{X^A_i}  W^A(X^A,X^B) + \frac{D^A}{N_A} \sum_{i=1}^{N_A} \Delta \phi(X^A_i) + \frac{\sqrt{2D^A}}{N_A} \sum_{i=1}^{N_A} \nabla_x \phi(X_i^A) \cdot \frac{dB_i^A}{dt}.
  \end{align}
\noindent As the $dB_i^A$ are independent and independent of $\nabla_x \phi(X^A_i)$, it can be shown that in the limit of a large number of particles this term can be neglected \cite{Barre2016}. We get:
\begin{equation*}
\begin{split}
 \frac{d}{dt}\langle f^A_N,\phi \rangle= -  \frac{\mu}{N_A} \sum_{i=1}^{N_A} \nabla_x \phi(X^A_i) \cdot &\bigg(\sum_{k_1=1}^{K_{AA}} (\nabla_{x_1}\Phi^{AA}  \delta_{i(k_1)}(i) + \nabla_{x_2}\Phi^{AA} \delta_{j(k_1)}(i))(X^A_{i(k_1)},X^A_{j(k_1)}) \\
 & + \sum_{k_3=1}^{K_{AB}} \nabla_{x_1}\Phi^{AB} \delta_{i(k_3)}(i)(X^A_{i(k_3)},X^B_{\ell(k_3)})\bigg) + \frac{D^A}{N_A} \sum_{i=1}^{N_A} \Delta \phi(X^A_i).
\end{split}
\end{equation*}
\noindent Now, exchanging the sums in $i$ and $k_1$ and $i$ and $k_3$ in the previous equation, one obtains:
\begin{equation*}
\begin{split}
 \frac{d}{dt}\langle f^A_N,\phi \rangle = -  \frac{\mu}{N_A} & \sum_{k_1=1}^{K_{AA}} \bigg(\nabla_x \phi(X^A_{i(k_1)}) \cdot  \nabla_{x_1}\Phi^{AA}(X^A_{i(k_1)},X^A_{j(k_1)}) + \nabla_x \phi(X^A_{j(k_1)}) \cdot  \nabla_{x_2}\Phi^{AA}(X^A_{i(k_1)},X^A_{j(k_1)}) \bigg)\\
 & -  \frac{\mu}{N_A} \sum_{k_3=1}^{K_{AB}} \nabla_x \phi(X^A_{i(k_3)}) \cdot \nabla_{x_1} \Phi^{AB}(X^A_{i(k_3)},X^B_{\ell(k_3)}) + \frac{D^A}{N_A} \sum_{i=1}^{N_A} \Delta \phi(X^A_i).
\end{split}
\end{equation*}
\noindent From the symmetry of $\Phi^{AA}$, we have:
$$ 
\nabla_{x_2}\Phi^{AA}(X^A_{i(k_1)},X^A_{j(k_1)}) = \nabla_{x_1} \Phi^{AA}(X^A_{j(k_1)},X^A_{i(k_1)}),$$
 \noindent leading to:
\begin{equation*}
\begin{split}
 \frac{d}{dt}\langle f^A_N,\phi \rangle = -  \frac{\mu}{N_A} & \sum_{k_1=1}^{K_{AA}} \bigg(\nabla_x \phi(X^A_{i(k_1)}) \cdot  \nabla_{x_1}\Phi^{AA}(X^A_{i(k_1)},X^A_{j(k_1)}) + \nabla_x \phi(X^A_{j(k_1)}) \cdot  \nabla_{x_1}\Phi^{AA}(X^A_{j(k_1)},X^A_{i(k_1)}) \bigg)\\
 & -  \frac{\mu}{N_A} \sum_{k_3=1}^{K_{AB}} \nabla_x \phi(X^A_{i(k_3)}) \cdot \nabla_{x_1} \Phi^{AB}(X^A_{i(k_3)},X^B_{\ell(k_3)}) + \frac{D^A}{N_A} \sum_{i=1}^{N_A} \Delta \phi(X^A_i).
\end{split}
\end{equation*}
\noindent or again:
\begin{equation*}
\begin{split}
 \frac{d}{dt}\langle f^A_N,\phi \rangle = &  - 2\mu  \langle  \hspace{-0.8mm} \langle g^{AA}_N, \nabla_{x_1} \Phi^{AA}(x_1,x_2) \cdot \nabla_x \phi(x_1) \rangle \hspace{-0.8mm} \rangle \\
 &- \mu \langle  \hspace{-0.8mm} \langle g^{AB}_N, \nabla_{x_1} \Phi^{AB}(x_1,x_2) \cdot \nabla_x \phi(x_1) \rangle \hspace{-0.8mm} \rangle + D^A \langle f^A_N, \Delta \phi \rangle\\
  = &2\mu   \langle  \hspace{-0.8mm} \langle  \nabla_{x_1}  \cdot \big( g^{AA}_N(x_1,x_2) \nabla_{x_1} \Phi^{AA}(x_1,x_2)\big) , \phi(x_1) \rangle \hspace{-0.8mm} \rangle\\
 & + \mu   \langle  \hspace{-0.8mm} \langle\nabla_{x_1} \cdot \big(g^{AB}_N(x_1,x_2) \nabla_{x_1} \Phi^{AB}(x_1,x_2) \big),  \phi(x_1) \rangle \hspace{-0.8mm} \rangle + D^A \langle \Delta f^A_N,  \phi \rangle,
\end{split}
\end{equation*}
\noindent where we have formally integrated by parts to obtain the second equality. We then exchange the order of integration 
% to get:
% \begin{align*}
% \frac{d}{dt} \int f^A_N(x_1,t) \phi(x_1) dx_1 = & 2\mu \int \int \nabla_{x_1} \cdot \big(g^{AA}_N (x_1,x_2) \nabla_{x_1} \Phi^{AA}(x_1,x_2)\big) \phi(x_1) dx_1dx_2 \\
% &+ \mu  \int \int \nabla_{x_1} \cdot \big(g^{AB}_N (x_1,x_2) \nabla_{x_1} \Phi^{AB}(x_1,x_2) \big) \phi(x_1) dx_1dx_2 \\
% &+ D \int \Delta f^A_N(x_1,t)\phi(x_1)dx_1 \\
% = & 2\mu  \int \nabla_{x_1} \cdot \bigg(\int g^{AA}_N (x_1,x_2) \nabla_{x_1} \Phi^{AA}(x_1,x_2) dx_2\bigg)  \phi(x_1) dx_1 \\
% &+ \mu  \int \nabla_{x_1} \cdot \bigg( \int g^{AB}_N (x_1,x_2) \nabla_{x_1} \Phi^{AB}(x_1,x_2) dx_2 \bigg) \phi(x_1) dx_1 \\
% &+ D^A \int \Delta f^A_N(x_1,t)\phi(x_1)dx_1. \\
% \end{align*}
% \noindent Now, in
and pass to the limit $N_A, N_B \rightarrow \infty$. If  $f^A_N \rightarrow f^{A}$, $g^{AA}_N \rightarrow g^{AA}, g^{AB}_N \rightarrow g^{AB}$, then:
\begin{equation}\label{systfA}
 \frac{df^A}{dt} = 2\mu \nabla_x \cdot F^{AA}[g^{AA}](x,t) + \mu  \nabla_x \cdot F^{AB}[g^{AB}](x,t) + D\Delta f^A,
 \end{equation}
 \noindent where:
 \begin{align*}
  F^{AA}[g](x,t) &= \int g^{AA}(x_1,x_2,t) \nabla_{x_1} \Phi^{AA}(x_1,x_2) dx_2\\
    F^{AB}[g](x,t) &= \int g^{AB}(x_1,x_2,t) \nabla_{x_1} \Phi^{AB}(x_1,x_2) dx_2.
 \end{align*}
\noindent Similarly, we can show that if $f^B_N \rightarrow f^{B}$, $g^{BB}_N \rightarrow g^{BB}, g_N^{BA} \rightarrow g^{BA}$ as 
$N_A,N_B \rightarrow \infty$ with $N_A/N_B\to r_{AB}$, we get:
\begin{equation}\label{systf}
 \frac{df^B}{dt} = 2\mu \nabla_x \cdot F^{BB}[g^{BB}](x,t) + \mu \nabla_x \cdot F^{BA}[g^{BA}](x,t) + D\Delta f^B,
 \end{equation}
 \noindent with:
 \begin{align*}
  F^{BB}[g](x,t) &= \int g(x_1,x_2,t) \nabla_{x_1} \Phi^{BB}(x_1,x_2) dx_2\\
    F^{BA}[g](x,t) &= \int g(x_1,x_2,t) \nabla_{x_1} \Phi^{BA}(x_1,x_2) dx_2.
 \end{align*}
\noindent We now turn towards the computation of the inter- and intra- species link distributions.

\subsection{Evolution equation for the particle links}
Here we develop only the computations for one intraspecies link distribution, namely $g^{AA}_N$, the computation of $g^{BB}_N$ being similar. From its asymmetry, the computation of the interspecies link distribution needs special treatment and we will develop  the computation of $g^{AB}_N$ further. 
%and leave the computation of $g_B^{K_{AB}}$ to the reader. 
We remark that the noise in \eqref{IBMA} transforms directly into a linear
diffusion term for $f^A$, all other contributions, analogously as in \ref{two_diff}, vanish in the large $N_A$ limit. It is not difficult to see that the same simplification takes place for $g^{AA}_N$ in the $N_{A} \rightarrow \infty$ limit. Thus, to reduce the
computations we will first use \eqref{IBMA} without noise, and reintroduce the diffusion term in the end.

\subsubsection*{Intraspecies link distribution}

Following the same principle as for $f^A_N$, one can write:
\begin{equation}\label{dg1}
\begin{split}
\frac{d}{dt} \langle  \hspace{-0.8mm} \langle  g^{AA}_N,\Psi  \rangle  \hspace{-0.8mm} \rangle  =\frac{1}{2N_{A}} \sum_{k=1}^{K_{AA}}& \bigg[ \nabla_{x_1} \Psi(X^A_{i(k)},X^A_{j(k)})  \cdot \frac{dX^A_{i(k)}}{dt} + \nabla_{x_1} \Psi(X^A_{j(k)},X^A_{i(k)}) \cdot \frac{dX^A_{j(k)}}{dt}\\
 & + \nabla_{x_2} \Psi(X^A_{i(k)},X^A_{j(k)}) \cdot \frac{dX^A_{j(k)}}{dt} + \nabla_{x_2} \Psi(X^A_{j(k)},X^A_{i(k)}) \cdot \frac{dX^A_{i(k)}}{dt} \bigg]\\
  = E_1 + E_2, 
 \end{split}
 \end{equation}
 \noindent where $E_k$ corresponds to the $k$-th line of~\eqref{dg1}. For the sake of simplicity, the computation of $E_1$ only is developed here. The computation of the other ones are similar and omitted. 
 
From Eq.~\eqref{IBMA}, one obtains: 
 \begin{equation*}
 \begin{split}
 E_1  = &\frac{1}{2N_{A}} \sum_{k=1}^{K_{AA}} \bigg[\nabla_{x_1} \Psi(X^A_{i(k)},X^A_{j(k)}) \cdot \frac{dX^A_{i(k)}}{dt} + \nabla_{x_1} \Psi(X^A_{j(k)},X^A_{i(k)}) \cdot \frac{dX^A_{j(k)}}{dt}\bigg]\\
  = &- \frac{\mu}{2N_{A}} \sum_{k=1}^{K_{AA}} \bigg[
   \nabla_{x_1} \Psi(X^A_{i(k)},X^A_{j(k)}) \cdot \sum_{k_1=1}^{K_{AA}}  \bigg(\nabla_{x_1} \Phi^{AA} \delta_{(i(k_1),i(k))} + \nabla_{x_2} \Phi^{AA} \delta_{(j(k_1),i(k))}\bigg)(X^A_{i(k_1)},X^A_{j(k_1)})\\
   &- \frac{\mu}{2N_{A}} \sum_{k=1}^{K_{AA}} \bigg[
   \nabla_{x_1} \Psi(X^A_{i(k)},X^A_{j(k)}) \cdot \sum_{k_3=1}^{K_{AB}}  \nabla_{x_1} \Phi^{AB} \delta_{(i(k_3),i(k))} (X^A_{i(k_3)},X^B_{j(k_3)})\bigg]\\
&  - \frac{\mu}{2N_{A}} \sum_{k=1}^{K_{AA}} \bigg[
   \nabla_{x_1} \Psi(X^A_{j(k)},X^A_{i(k)}) \cdot \sum_{k_1=1}^{K_{AA}}  \bigg(\nabla_{x_1} \Phi^{AA} \delta_{(i(k_1),j(k))} + \nabla_{x_2} \Phi^{AA} \delta_{(j(k_1),j(k))}\bigg)(X^A_{i(k_1)},X^A_{j(k_1)})\\
      &- \frac{\mu}{2N_{A}} \sum_{k=1}^{K_{AA}} \bigg[
   \nabla_{x_1} \Psi(X^A_{j(k)},X^A_{i(k)}) \cdot \sum_{k_3=1}^{K_{AB}}  \nabla_{x_1} \Phi^{AB} \delta_{(i(k_3),j(k))} (X^A_{i(k_3)},X^B_{j(k_3)})\bigg].
 \end{split}
 \end{equation*}
 \noindent  Now, exchanging the sums in $k$ and $k_1$ and $k$ and $k_3$ and using the symmetry of $\Phi^{AA}$, one obtains:
  \begin{equation}\label{dg}
 \begin{split}
  E_1= 
 - \frac{\mu}{2N_{A}} \sum_{k_1=1}^{K_{AA}} \nabla_{x_1} \Phi^{AA}(X^A_{i(k_1)},X^A_{j(k_1)}) \cdot \sum_{k=1}^{K_{AA}} & \big(\nabla_{x_1} \Psi(X^A_{i(k)},X^A_{j(k)}) \delta_{(i(k),i(k_1))}
 + \nabla_{x_1} \Psi(X^A_{j(k)},X^A_{i(k)}) \delta_{(j(k),i(k_1))} \big) \\
- \frac{\mu}{2N_{A}} \sum_{k_1=1}^{K_{AA}} \nabla_{x_1} \Phi^{AA}(X^A_{j(k_1)},X^A_{i(k_1)}) \cdot \sum_{k=1}^{K_{AA}} & \big(\nabla_{x_1} \Psi(X^A_{i(k)},X^A_{j(k)}) \delta_{(i(k),j(k_1))}
 + \nabla_{x_1} \Psi(X^A_{j(k)},X^A_{i(k)}) \delta_{(j(k),j(k_1))} \big) \\
 - \frac{\mu}{2N_{A}} \sum_{k_3=1}^{K_{AB}} \nabla_{x_1} \Phi^{AB}(X^A_{i(k_3)},X^B_{j(k_3)}) \cdot \sum_{k=1}^{K_{AA}} & \big(\nabla_{x_1} \Psi(X^A_{i(k)},X^A_{j(k)}) \delta_{(i(k),i(k_3))}
 + \nabla_{x_1} \Psi(X^A_{j(k)},X^A_{i(k)}) \delta_{(j(k),i(k_3))} \big)
  \end{split}
 \end{equation}
 In the first two lines of the above expression, $k=k_1$ in the internal sums play a special role: indeed we know that the $k=k_1$ term, at variance with the other terms, always contribute to the sum. On the contrary, $k=k_3$ in the last line above should not be distinguished. So, we get:
 %\red{the link $k_3$ should not be distinguished as it does not say anything about the links between A and another A, so I still approximate it using the old %technique}. We know that they always contribute to the sum, and so they should be distinguished. For example, the first term is equal to
 \eqh{
&\sum_{k=1}^{K_{AA}}  \big(\nabla_{x_1} \Psi(X^A_{i(k)},X^A_{j(k)}) \delta_{(i(k),i(k_1))}
 + \nabla_{x_1} \Psi(X^A_{j(k)},X^A_{i(k)}) \delta_{(j(k),i(k_1))} \big) \\
 &=
 \Grad_{x_1}  \Psi(X^A_{i(k_1)},X^A_{j(k_1)}) +
 \sum_{k\neq k_1} \big(\nabla_{x_1} \Psi(X^A_{i(k)},X^A_{j(k)}) \delta_{(i(k),i(k_1))}
 + \nabla_{x_1} \Psi(X^A_{j(k)},X^A_{i(k)}) \delta_{(j(k),i(k_1))} \big) 
 } 
 
\noindent Because there is no restriction on the number of links per particle, the sums over $k$ cannot be simplified in this case. In order to express the terms in the $k\neq k_1$ sum, we define the number of intra-species links connected to a particle:
$$
C_i^{k_1} = \mbox{Card}(\{ k \; | \; i(k)=i(k_1) \; or \; j(k) = i(k_1) \}), 
$$
\noindent where Card denotes the cardinal of a set. 
Then, as $N_{A} \rightarrow \infty$, we assume that the following mean-field approximation holds for any chosen link $k_1$: 
 \begin{equation}
 \begin{split}
 \frac{1}{2C_i^{k_1}} \sum_{k\neq k_1} \big( \Psi(X^A_{i(k)},X^A_{j(k)})&\delta_{i(k),i(k_1)} + \Psi(X^A_{j(k)},X^A_{i(k)}) \delta_{j(k),i(k_1)})\big)  \underset{N_{A} \rightarrow \infty}{\rightarrow} \int (\Psi P^{AA})(X^A_{i(k_1)},x_2) dx_2,
\end{split}
\label{eq:MF1}
 \end{equation}
 \noindent where 
 \begin{equation}
P^{AA}(X^A_{i(k_1)},x_2)= \frac{g^{AA}(X^A_{i(k_1)},x_2)}{\int g^{AA}(X^A_{i(k_1)},x_2) dx_2}  ,
\label{eq:MF2}
 \end{equation}
 \noindent is the conditional probability of finding an intraspecies link conditioned on the fact that one of the particles of this link has the same location as $i(k_1)$. Then, as $N_A \rightarrow \infty$ , $C_i^{k_1}$ is the mean number of intraspecies links per particle. The mean number of intraspecies links AA in the volume $dX^A_{i(k_1)}$ is $N_{A}\int g^{AA}(X^A_{i(k_1)},x_2) dx_2$ and the mean number of particles of type A in $dX^A_{i(k_1)}$ is $N_Af^A(X^A_{i(k_1)})$. Thus:
 \begin{equation*}
  C_i^{k_1} \underset{N_A \rightarrow \infty }{\longrightarrow} \frac{\int g^{AA}(X^A_{i(k_1)},x_2) dx_2 }{f^A(X^A_{i(k_1)})}.
 \end{equation*}
 \noindent So, we get:
 \begin{equation*}
 \begin{split}
\sum_{k\neq k_1}\big( \Psi(X^A_{i(k)},X^A_{j(k)})&\delta_{i(k),i(k_1)} + \Psi(X^A_{j(k)},X^A_{i(k)}) \delta_{j(k),i(k_1)})\big)\underset{ N_A \rightarrow \infty }{\rightarrow}
 \frac{2}{f^A(X^A_{i(k_1)})} \int (\Psi g^{AA})(X_{i(k_1)}, x_2) dx_2.
 \end{split}
 \end{equation*}
 \noindent Inserting these expressions in Eq.~\eqref{dg}, one obtains, when $N_A, N_{B}$ are large:
  \begin{equation*}
 \begin{split}
&\underset{\underset{N_A/N_B\rightarrow r_{AB}}{N_A,N_B \rightarrow \infty }}{\lim}  E_1 =
-
\underset{\underset{N_{A}/N_B \rightarrow r_{AB}>0}{ N_A,N_B \rightarrow \infty }}{\lim}   \frac{\mu}{2N_{A}} \sum_{k_1=1}^{K_{AA}} \nabla_{x_1} \Phi^{AA}(X^A_{i(k_1)},X^A_{j(k_1)})\cdot
\Grad_{x_1}  \Psi(X^A_{i(k_1)},X^A_{j(k_1)})\\
&-
\underset{\underset{N_{A}/N_B \rightarrow r_{AB}>0}{ N_A,N_B \rightarrow \infty }}{\lim}   \frac{\mu}{2N_{A}} \sum_{k_1=1}^{K_{AA}} \nabla_{x_1} \Phi^{AA}(X^A_{j(k_1)},X^A_{i(k_1)})\cdot
\Grad_{x_1}  \Psi(X^A_{j(k_1)},X^A_{i(k_1)})\\
  &- \underset{\underset{N_{A}/N_B \rightarrow r_{AB}>0}{ N_A,N_B \rightarrow \infty }}{\lim} \frac{\mu}{N_{A}} \sum_{k_1=1}^{K_{AA}} \nabla_{x_1} \Phi^{AA}(X^A_{i(k_1)},X^A_{j(k_1)}) \cdot \psi^1_{AA}(X^A_{i(k_1)}) + \nabla_{x_1} \Phi^{AA}(X^A_{j(k_1)},X^A_{i(k_1)}) \cdot \psi^1_{AA}(X^A_{j(k_1)}) \\
&- \underset{\underset{N_{A}/N_B \rightarrow r_{AB}>0} {N_A,N_B \rightarrow \infty }}{\lim} \frac{\mu}{N_{A}} \sum_{k_3=1}^{K_{AB}} \nabla_{x_1} \Phi^{AB}(X^A_{i(k_3)},X^B_{j(k_3)}) \cdot \psi^1_{AA}(X^A_{i(k_3)}) 
%&-\underset{\underset{\frac{K_{AA}}{N_A} \rightarrow \xi^{AA}>0}{\underset{K_{AA} \rightarrow \infty}{ N_A \rightarrow \infty }}}{\lim} 
%- \frac{\mu}{2K_{AA}} \sum_{k_3=1}^{K_{AB}} \nabla_{x_1} \Phi^{AB}(X^A_{i(k_3)},X^B_{j(k_3)}) \cdot \nabla_{x_1} \Psi(X^A_{i(k_3)},X^A_{j(k_3)}) 
  \end{split}
 \end{equation*}
\noindent where 
 \begin{equation*}
\psi^1_{AA}(x_1) = \frac{1}{f^A(x_1)}\int \big( g^{AA} \nabla_{x_1} \Psi \big)(x_1,x_2)  dx_2.
 \end{equation*}
 \noindent Finally, we find:
  \begin{equation*}
 \begin{split}
E_1 \underset{\underset{\frac{N_{A}}{N_B}, \rightarrow r_{AB}>0}{ N_A,N_B \rightarrow \infty }}{\rightarrow}  &-2\mu \langle  \hspace{-0.8mm} \langle  g^{AA}, \nabla_{x_1} \Phi^{AA}(x_1,x_2) \cdot \psi^1_{AA}(x_1)  \rangle  \hspace{-0.8mm} \rangle -\mu \langle  \hspace{-0.8mm} \langle  g^{AB}, \nabla_{x_1} \Phi^{AB}(x_1,x_2) \cdot \psi^1_{AA}(x_1)  \rangle  \hspace{-0.8mm} \rangle \\
&-\mu\langle  \hspace{-0.8mm}\langle g^{AA},\Grad_{x_1}\Phi^{AA}(x_1,x_2)\cdot\Grad_{x_1}\Psi(x_1,x_2)  \rangle  \hspace{-0.8mm}\rangle
%-\frac{\mu}{2}\langle  \hspace{-0.8mm}\langle g^{AB},\Grad_{x^1}\Phi^{AB}(x_1,x_2)\Grad_{x^1}\Psi(x_1,x_2)  \rangle  \hspace{-0.8mm}\rangle
  \end{split}
 \end{equation*}
 
 \noindent After the same treatment for $E_2$ of Eq.~\eqref{dg1} and in the limit $N_A,N_B \rightarrow \infty, \frac{N_{A}}{N_B} \rightarrow r_{AB}>0$, one obtains the final equation for $g^{AA}$:
   \begin{align}
&\frac{d}{dt} \langle  \hspace{-0.8mm} \langle  g^{AA}(x_1,x_2),\Psi(x_1,x_2)  \rangle  \hspace{-0.8mm} \rangle \nonumber\\
&= -2\mu \langle  \hspace{-0.8mm} \langle  g^{AA}, \nabla_{x_1} \Phi^{AA}(x_1,x_2) \cdot \psi^1_{AA}(x_1)  \rangle  \hspace{-0.8mm} \rangle -\mu \langle  \hspace{-0.8mm} \langle  g^{AB}, \nabla_{x_1} \Phi^{AB}(x_1,x_2) \cdot \psi^1_{AA}(x_1)  \rangle  \hspace{-0.8mm} \rangle \\
& -2\mu \langle  \hspace{-0.8mm} \langle  g^{AA}, \nabla_{x_1} \Phi^{AA}(x_1,x_2) \cdot \psi^2_{AA}(x_1)  \rangle  \hspace{-0.8mm} \rangle -\mu \langle  \hspace{-0.8mm} \langle  g^{AB}, \nabla_{x_1} \Phi^{AB}(x_1,x_2) \cdot \psi^2_{AA}(x_1)  \rangle  \hspace{-0.8mm} \rangle  \\
&-{\mu\langle  \hspace{-0.8mm}\langle g^{AA},\Grad_{x_1}\Phi^{AA}(x_1,x_2)\cdot\Grad_{x_1}\Psi(x_1,x_2)  \rangle  \hspace{-0.8mm}\rangle
-\mu\langle  \hspace{-0.8mm}\langle g^{AA},\Grad_{x_2}\Phi^{AA}(x_1,x_2)\cdot\Grad_{x_2}\Psi(x_1,x_2)  \rangle  \hspace{-0.8mm}\rangle}
, \label{simpg}
  \end{align}
\noindent where,
 \begin{equation*}
\psi^2_{AA}(x_1) = \frac{1}{f^A(x_1)}\int \big( g^{AA} \nabla_{x_2} \Psi \big)(x_2,x_1)  dx_2.
 \end{equation*}
\noindent Integrating by parts, changing the variables and order of integrals we easily obtain:
   \begin{align*}
&\frac{d}{dt} \langle  \hspace{-0.8mm} \langle  g^{AA}(x_1,x_2),\Psi(x_1,x_2)  \rangle  \hspace{-0.8mm} \rangle \nonumber\\
&= 2\mu  \langle  \hspace{-0.8mm} \langle \nabla_{x_1} \cdot \bigg( \frac{g^{AA}(x_1,x_2)}{f^A(x_1)} \int \big(g^{AA}\nabla_{x_1} \Phi^{AA}\big)(x_1,x_2) dx_2 \bigg), \Psi(x_1,x_2)  \rangle  \hspace{-0.8mm} \rangle\\
& +\mu \langle  \hspace{-0.8mm} \langle \nabla_{x_1} \cdot \bigg( \frac{g^{AA}(x_1,x_2)}{f^A(x_1)} \int \big(g^{AB}\nabla_{x_1} \Phi^{AB}\big)(x_1,x_2) dx_2 \bigg), \Psi(x_1,x_2)  \rangle  \hspace{-0.8mm} \rangle \\
&+2\mu  \langle  \hspace{-0.8mm} \langle \nabla_{x_2} \cdot \bigg( \frac{g^{AA}(x_1,x_2)}{f^A(x_2)} \int \big(g^{AA}\nabla_{x_1} \Phi^{AA}\big)(x_2,x_1) dx_1 \bigg), \Psi(x_1,x_2)  \rangle  \hspace{-0.8mm} \rangle\\
& +\mu  \langle  \hspace{-0.8mm} \langle \nabla_{x_2} \cdot \bigg( \frac{g^{AA}(x_1,x_2)}{f^A(x_2)} \int \big(g^{AB}\nabla_{x_1} \Phi^{AB}\big)(x_2,x_1) dx_1 \bigg), \Psi(x_1,x_2)  \rangle  \hspace{-0.8mm} \rangle\\
  &{+\mu \langle  \hspace{-0.8mm} \langle \nabla_{x_1} \cdot \bigg( g^{AA}(x_1,x_2) \nabla_{x_1} \Phi^{AA}(x_1,x_2)  \bigg), \Psi(x_1,x_2)  \rangle  \hspace{-0.8mm} \rangle}\\
  &{+\mu \langle  \hspace{-0.8mm} \langle \nabla_{x_2} \cdot \bigg( g^{AA}(x_1,x_2) \nabla_{x_2} \Phi^{AA}(x_1,x_2)  \bigg), \Psi(x_1,x_2)  \rangle  \hspace{-0.8mm} \rangle}.
 \end{align*}
\noindent Finally, restoring the noise, we obtain the final equation for $g^{AA}$:
\begin{align}
\partial_t g^{AA}&(x_1,x_2,t) = D^A \big(\Delta_{x_1} g^{AA}(x_1,x_2,t) +\Delta_{x_2} g^{AA}(x_1,x_2,t) \big) \nonumber\\
&+ 2\mu  \nabla_{x_1} \cdot \bigg(\frac{g^{AA}(x_1,x_2,t)}{f^A(x_1,t)}F^{AA}[g^{AA}](x_1,t)\bigg)+ \mu \nabla_{x_1} \cdot \bigg(\frac{g^{AA}(x_1,x_2,t)}{f^A(x_1,t)}F^{AB}[g^{AB}](x_1,t)\bigg)\label{EqgAA}\\
&+ 2\mu \nabla_{x_2} \cdot \bigg(\frac{g^{AA}(x_1,x_2,t)}{f^A(x_2,t)}F^{AA}[g^{AA}](x_2,t)\bigg)+ \mu \nabla_{x_2} \cdot \bigg(\frac{g^{AA}(x_1,x_2,t)}{f^A(x_2,t)}F^{AB}[g^{AB}](x_2,t)\bigg) \nonumber \\
&{+\mu \nabla_{x_1} \cdot \bigg( g^{AA}(x_1,x_2) \nabla_{x_1} \Phi^{AA}(x_1,x_2)  \bigg)+\mu \nabla_{x_2} \cdot \bigg( g^{AA}(x_1,x_2) \nabla_{x_2} \Phi^{AA}(x_1,x_2)  \bigg)},\nonumber
\end{align}
\noindent where
\begin{align*}
F^{AA}[g^{AA}](x_1,t) &= \int g^{AA}(x_1,x_2) \nabla_{x_1} \Phi^{AA}(x_1,x_2) dx_2, \\
%\; F^{AA}[g^{AA}](x_2,t) = \int g^{AA}(x_2,x_1) \nabla_{x_1} \Phi^{AA}(x_2,x_1) dx_1\\
F^{AB}[g^{AB}](x_1,t) &= \int g^{AB}(x_1,x_2) \nabla_{x_1} \Phi^{AB}(x_1,x_2) dx_2. %\; F^{AB}[g^{AB}](x_2,t) = \int g^{AB}(x_2,x_1) \nabla_{x_1} \Phi^{AB}(x_2,x_1) dx_1.
\end{align*}
Eq. \eqref{EqgAA} does not take into account the phenomena of creation and destruction
of intraspecies links. According to the description at the beginning of this paper, our model describes a process of creation of links with rate $\nu^{AA}_{c,N,\ep}$, provided the two type $A$ particles are sufficiently close to each others. Hence, the number of new intraspecies links will be proportional to the number of pairs of particles such that one of them is close to $x_1$ and the other one is close to $x_2$,
whose distance is less than $R$:
$$
\frac{N_A(N_A - 1)}{2} h^{AA}(x_1, x_2, t)\chi_{\{|x_1-x_2|\leq R\}} dx_1 dx_2 dt,
$$
\noindent where $h^{AA}(x_1, x_2, t) = \lim_{N_A \rightarrow \infty} h^{AA}_N(x_1,x_2,t)$ and where $h^{AA}_N(x_1,x_2,t)$ is the two particle distribution function for particles of type A defined by Eq.\eqref{TPDFAA}. This number has to be decreased by the number of pairs of particles of the same type that are already connected by existing links
$$
N_{A} g^{AA}(x_1,x_2,t)dx_1dx_2.
$$
Therefore, the average number of new intraspecies links created in the interval $[t,t+\Delta t[$ is equal to 
$$
\nu^{AA}_{c,N,\ep} \bigg(\frac{N_A(N_A - 1)}{2} h^{AA}(x_1, x_2, t)\chi_{\{|x_1-x_2|\leq R\}} - N_{A} g^{AA}(x_1,x_2,t) \bigg) dx_1 dx_2 \Delta t 
$$
\noindent Dividing this expression by $N_{A}$ and using \eqref{eq:nuscalings}, the rate of creation of new  intraspecies link at $(x_1,x_2)$ is in the limit $N_A\to \infty$:
$$
\frac{\nu^{AA}_{c,\ep}  }{2} h^{AA}(x_1, x_2, t)\chi_{\{|x_1-x_2|\leq R\}}.
$$
\noindent Notice the scaling of $\nu^{AA}_{c,N,\ep}$ in \eqref{eq:nuscalings}: it ensures that among the $O(N_A^2)$ possible links, only $O(N_A)$ are effectively present, and the 
total number of AA links $K^{AA}=O(N_A)$. The rate of destruction of existing intraspecies link at $(x_1,x_2)$ is:
$$
\nu^{AA}_{d,\ep}  g^{AA}(x_1,x_2,t).
$$
\noindent Including these source terms in \eqref{EqgAA}, we obtain
\begin{align}
\partial_t g^{AA}&(x_1,x_2,t) = D^A \big(\Delta_{x_1} g^{AA}(x_1,x_2,t) +\Delta_{x_2} g^{AA}(x_1,x_2,t) \big) \nonumber\\
&+ 2\mu  \nabla_{x_1} \cdot \bigg(\frac{g^{AA}(x_1,x_2,t)}{f^A(x_1,t)}F^{AA}[g^{AA}](x_1,t)\bigg)+ \mu  \nabla_{x_1} \cdot \bigg(\frac{g^{AA}(x_1,x_2,t)}{f^A(x_1,t)}F^{AB}[g^{AB}](x_1,t)\bigg)\label{EqgAAb}\\
&+ 2\mu \nabla_{x_2} \cdot \bigg(\frac{g^{AA}(x_1,x_2,t)}{f^A(x_2,t)}F^{AA}[g^{AA}](x_2,t)\bigg)+ \mu \nabla_{x_2} \cdot \bigg(\frac{g^{AA}(x_1,x_2,t)}{f^A(x_2,t)}F^{AB}[g^{AB}](x_2,t)\bigg) \nonumber \\
&{+\mu \nabla_{x_1} \cdot \bigg( g^{AA}(x_1,x_2) \nabla_{x_1} \Phi^{AA}(x_1,x_2)  \bigg)+\mu \nabla_{x_2} \cdot \bigg( g^{AA}(x_1,x_2) \nabla_{x_2} \Phi^{AA}(x_1,x_2)  \bigg)}\nonumber\\
&+ \frac{\nu^{AA}_{c,\ep}  }{2} h^{AA}(x_1, x_2, t)\chi_{\{|x_1-x_2|\leq R\}}
- \nu^{AA}_{d,\ep}  g^{AA}(x_1,x_2,t).\nonumber
\end{align}

\noindent Quite straightforwardly, we can show that in the limit $N_A,N_B \rightarrow \infty, \frac{N_{A}}{N_B} \rightarrow r_{AB}>0$, $g^{BB}(x_1,x_2,t)$ solves:
\begin{align}
\partial_t g^{BB}&(x_1,x_2,t) = D^B \big(\Delta_{x_1} g^{BB}(x_1,x_2,t) +\Delta_{x_2} g^{BB}(x_1,x_2,t) \big) \nonumber\\
&+ 2\mu \nabla_{x_1} \cdot \bigg(\frac{g^{BB}(x_1,x_2,t)}{f^B(x_1,t)}F^{BB}[g^{BB}](x_1,t)\bigg)+ \mu  \nabla_{x_1} \cdot \bigg(\frac{g^{BB}(x_1,x_2,t)}{f^B(x_1,t)}F^{BA}[g^{BA}](x_1,t)\bigg)\nonumber\\
&+ 2\mu \nabla_{x_2} \cdot \bigg(\frac{g^{BB}(x_1,x_2,t)}{f^B(x_2,t)}F^{BB}[g^{BB}](x_2,t)\bigg)+ \mu \nabla_{x_2} \cdot \bigg(\frac{g^{BB}(x_1,x_2,t)}{f^B(x_2,t)}F^{BA}[g^{BA}](x_2,t)\bigg) \nonumber \\
&{+\mu \nabla_{x_1} \cdot \bigg( g^{BB}(x_1,x_2) \nabla_{x_1} \Phi^{BB}(x_1,x_2)  \bigg)+\mu \nabla_{x_2} \cdot \bigg( g^{BB}(x_1,x_2) \nabla_{x_2} \Phi^{BB}(x_1,x_2)  \bigg)}\nonumber\\
&+ \frac{\nu^{BB}_{c,\ep}  }{2} h^{BB}(x_1, x_2, t)\chi_{\{|x_1-x_2|\leq R\}}
- \nu^{BB}_{d,\ep}  g^{BB}(x_1,x_2,t),\nonumber
\end{align}
\noindent where  {we have assumed the scaling \eqref{eq:nuscalings} for $\nu^{BB}_{c,N,\ep} N_B$}, and 
\begin{align*}
F^{BB}[g^{BB}](x_1,t) &= \int g^{BB}(x_1,x_2) \nabla_{x_1} \Phi^{BB}(x_1,x_2) dx_2, \\
% \; F^{BB}[g^{BB}](x_2,t) = \int g^{BB}(x_2,x_1) \nabla_{x_1} \Phi^{BB}(x_2,x_1) dx_1\\
F^{BA}[g^{BA}](x_1,t) &= \int g^{BA}(x_1,x_2) \nabla_{x_1} \Phi^{BA}(x_1,x_2) dx_2.% \; F^{BA}[g^{BA}](x_2,t) = \int g^{BA}(x_2,x_1) \nabla_{x_1} \Phi^{BA}(x_2,x_1) dx_1.
\end{align*}
We stress the fact that $\Phi^{AB}(x_1,x_2)$ (force of a particle of type $B$ close to $x_2$ exerted on a particle $A$ close to $x_1$) is not necessarily equal to $\Phi^{BA}(x_2,x_1)$ (force of a particle of type $A$ close to $x_1$ exerted on a particle $B$ close to $x_2$).

\subsubsection*{Computations of the interspecies link distribution}
Here, we develop the computations for the interspecies link distribution $g^{AB}(x_1,x_2,t)$. Proceeding as before, we write:
\begin{equation}\label{dg1a}
\begin{split}
\frac{d}{dt} \langle  \hspace{-0.8mm} \langle  g^{AB}_N,\Psi  \rangle  \hspace{-0.8mm} \rangle  =\frac{1}{N_{A}} \sum_{k=1}^{K_{AB}}& \bigg[\nabla_{x_1} \Psi(X^A_{i(k)},X^B_{\ell(k)})  \cdot \frac{dX^A_{i(k)}}{dt}  + \nabla_{x_2} \Psi(X^A_{i(k)},X^B_{\ell(k)}) \cdot \frac{dX^B_{\ell(k)}}{dt}\bigg]\\
  = E_1 + E_2, 
 \end{split}
 \end{equation}
\noindent where, ignoring the noise and using Eq.\eqref{IBMA}:
\begin{align*}
E_1 =& \frac{1}{N_{A}} \sum_{k=1}^{K_{AB}} \nabla_{x_1} \Psi(X^A_{i(k)},X^B_{\ell(k)})  \cdot \frac{dX^A_{i(k)}}{dt} \\
=&  \frac{-\mu}{N_{A}} \sum_{k=1}^{K_{AB}} \nabla_{x_1} \Psi(X^A_{i(k)},X^B_{\ell(k)})  \cdot \bigg[ \sum_{k_1=1}^{K_{AA}} \big(\nabla_{x_1} \Phi^{AA}(X^A_{i(k_1)},X^A_{j(k_1)})\delta_{i(k_1), i(k)} + \nabla_{x_2} \Phi^{AA}(X^A_{i(k_1)},X^A_{j(k_1)})\delta_{j(k_1), i(k)}\big) \bigg]\\
& \frac{-\mu}{N_{A}} \sum_{k=1}^{K_{AB}} \nabla_{x_1} \Psi(X^A_{i(k)},X^B_{\ell(k)})  \cdot \sum_{k_3=1}^{K_{AB}} \nabla_{x_1} \Phi^{AB}(X^A_{i(k_3)},X^B_{\ell(k_3)})\delta_{i(k_3),i(k)}
\end{align*}
\noindent Exchanging the sums in the first two terms and using the symmetry of $\Phi^{AA}$, we have
\begin{align*}
-\frac{\mu}{N_{A}}\sum_{k=1}^{K_{AB}}\nabla_{x_1} \Psi(X^A_{i(k)},X^B_{\ell(k)}) \cdot &\bigg[\sum_{k_1=1}^{K_{AA}} \nabla_{x_1} \Phi^{AA}(X^A_{i(k_1)},X^A_{j(k_1)})\delta_{i(k_1),i(k)} + \nabla_{x_1} \Phi^{AA}(X^A_{j(k_1)},X^A_{i(k_1)})\delta_{j(k_1),i(k)}\big] \\
& =  -\frac{\mu}{N_{A}} \bigg[\sum_{k_1=1}^{K_{AA}} \nabla_{x_1} \Phi^{AA}(X^A_{i(k_1)},X^A_{j(k_1)}) \cdot \sum_{k=1}^{K_{AB}} \nabla_{x_1} \Psi(X^A_{i(k)},X^B_{\ell(k)}) \delta_{i(k),i(k_1)}\\
& + \sum_{k_1=1}^{K_{AA}} \nabla_{x_1} \Phi^{AA}(X^A_{j(k_1)},X^A_{i(k_1)}) \cdot \sum_{k=1}^{K_{AB}} \nabla_{x_1} \Psi(X^A_{i(k)},X^B_{\ell(k)}) \delta_{i(k),j(k_1)}.\bigg]
\end{align*}
\noindent Now, in the same spirit as before we define the number of interspecies links linked to a particle of type $A$:
$$
C_{i,A}^{k_3} = \mbox{Card}(\{k | i(k) = i(k_3)\}),
$$
\noindent and we make the following mean-field assumption
\begin{align}
\frac{1}{C_{i,A}^{k_3}}  \sum_{k=1}^{K_{AB}} \nabla_{x_1} \Psi(X^A_{i(k)},X^B_{\ell(k)}) \delta_{i(k),i(k_3)}  \underset{N_A, N_{B} \rightarrow \infty, \frac{N_{A}}{N_B}\rightarrow r_{AB}}{\rightarrow} \int \big(\nabla_{x_1} \Psi P^{AB}\big)(X^A_{i(k_3)},x_2)dx_2, \label{eq:MF1b}
\end{align}
\noindent where $P^{AB}(X^A_{i(k_3)},x_2)$ is the conditional probability of finding an interspecies link conditioned on the fact that the type $A$ particle of the link has the same location as $i(k_3)$:
\begin{equation}
P^{AB}(X^A_{i(k_3)},x_2) = \frac{g^{AB}(X^A_{i(k_3)},x_2)}{\int g^{AB}(X^A_{i(k_3)},x_2) dx_2}. \label{eq:MF2b}
\end{equation}
\noindent Now as $N_A,N_B \rightarrow \infty$, $C_{i,A}^{k_3}$ is the mean number of interspecies links per particle of type A. The mean number of interspecies links the type A particle of which belonging to the volume $dX^A_{i(k_3)}$ is $N_{A} \int g^{AB}(X^A_{i(k_3)},x_2)dx_2$, and the mean number of particles of type A is $N_A f^A(X^A_{i(k_3)})$. Therefore, 
$$
C_{i,A}^{k_3} \underset{N_A, N_{B} \rightarrow \infty, \frac{N_{A}}{N_B}\rightarrow r_{AB}}{\rightarrow} \frac{\int g^{AB}(X^A_{i(k_3)},x_2) dx_2}{f^A(X^A_{i(k_3)})},
$$
\noindent leading, when $N_A,N_{B}$ are large, to
\begin{align*}
-\frac{\mu}{N_{A}}& \sum_{k_1=1}^{K_{AA}} \nabla_{x_1} \Phi^{AA}(X^A_{i(k_1)},X^A_{j(k_1)}) \cdot \sum_{k=1}^{K_{AB}} \nabla_{x_1} \Psi(X^A_{i(k)},X^B_{\ell(k)}) \delta_{i(k),i(k_1)} \\
&= -\frac{\mu}{N_{A}} \sum_{k_1=1}^{K_{AA}} \nabla_{x_1} \Phi^{AA}(X^A_{i(k_1)},X^A_{j(k_1)}) \cdot \bigg( \frac{N_{A} \int g^{AB}(X^A_{i(k_1)},x_2) \nabla_{x_1}\Psi(X^A_{i(k_1)},x_2) dx_2}{N_A f^A(X^A_{i(k_1)})}\bigg).
\end{align*}
\noindent and 
\begin{align*}
-\frac{\mu}{N_{A}}& \sum_{k_1=1}^{K_{AA}} \nabla_{x_1} \Phi^{AA}(X^A_{j(k_1)},X^A_{i(k_1)}) \cdot \sum_{k=1}^{K_{AB}} \nabla_{x_1} \Psi(X^A_{i(k)},X^B_{\ell(k)}) \delta_{i(k),j(k_1)}. \\
&= -\frac{\mu}{N_{A}} \sum_{k_1=1}^{K_{AA}} \nabla_{x_1} \Phi^{AA}(X^A_{j(k_1)},X^A_{i(k_1)}) \cdot \bigg( \frac{N_{A} \int g^{AB}(X^A_{j(k_1)},x_2) \nabla_{x_1} \Psi(X^A_{j(k_1)},x_2) dx_2}{N_A f^A(X^A_{j(k_1)})}\bigg).
\end{align*}
\noindent and altogether:
\begin{align*}
-\frac{\mu}{N_{A}}&\sum_{k=1}^{K_{AB}}\nabla_{x_1} \Psi(X^A_{i(k)},X^B_{\ell(k)}) \cdot \bigg[\sum_{k_1=1}^{K_{AA}} \nabla_{x_1} \Phi^{AA}(X^A_{i(k_1)},X^A_{j(k_1)})\delta_{i(k_1),i(k)} + \nabla_{x_1} \Phi^{AA}(X^A_{j(k_1)},X^A_{i(k_1)})\delta_{j(k_1),i(k)}\big] \\
&\underset{N_A, N_{B} \rightarrow \infty, \frac{N_{A}}{N_B}\rightarrow r_{AB}}{\rightarrow} - 2\mu \langle  \hspace{-0.8mm} \langle g^{AA}, \nabla_{x_1} \Phi^{AA}(x_1,x_2) \cdot \bigg( \frac{\int (\nabla_{x_1} \Psi g^{AB} )(x_1,x_2) dx_2}{ f^A(x_1)}\bigg) \rangle  \hspace{-0.8mm} \rangle.
\end{align*}
\noindent Now, exchanging the sums in the last term of $E_1$, we obtain
\begin{align*}
 &\frac{-\mu}{N_{A}} \sum_{k=1}^{K_{AB}} \nabla_{x_1} \Psi(X^A_{i(k)},X^B_{\ell(k)})  \cdot \sum_{k_3=1}^{K_{AB}} \nabla_{x_1} \Phi^{AB}(X^A_{i(k_3)},X^B_{\ell(k_3)})\delta_{i(k_3),i(k)} \\
 = &  {-\frac{\mu}{N_{A}} \sum_{k_3=1}^{K_{AB}} \nabla_{x_1} \Phi^{AB}(X^A_{i(k_3)},X^B_{\ell(k_3)})\cdot
 \nabla_{x_1} \Psi(X^A_{i(k_3)},X^B_{\ell(k_3)})}\\
&- \frac{\mu}{N_{A}} \sum_{k_3=1}^{K_{AB}} \nabla_{x_1} \Phi^{AB}(X^A_{i(k_3)},X^B_{\ell(k_3)}) \cdot \sum_{k\neq k_3} \nabla_{x_1} \Psi(X^A_{i(k)},X^B_{\ell(k)}) \delta_{i(k),i(k_3)},
\end{align*}
where, as previously, we distinguish the link $k_3$, that always contributes to the inner sum. Using the theory previously developed, we can write, when $N_A,N_{B}$ are large:
\begin{align*}
 \frac{-\mu}{N_{A}}& \sum_{k_3=1}^{K_{AB}} \nabla_{x_1} \Phi^{AB}(X^A_{i(k_3)},X^B_{\ell(k_3)}) \cdot \sum_{k\neq k_3} \nabla_{x_1} \Psi(X^A_{i(k)},X^B_{\ell(k)}) \delta_{i(k),i(k_3)} \\
 & =  \frac{-\mu}{N_{A}} \sum_{k_3=1}^{K_{AB}} \nabla_{x_1} \Phi^{AB}(X^A_{i(k_3)},X^B_{\ell(k_3)}) \cdot  \bigg(\frac{1}{f^A(X^A_{i(k_3)})}\int \big(\nabla_{x_1} \Psi g^{AB}\big)(X^A_{i(k_3)},x_2)dx_2\bigg)\\
 &\underset{N_A, N_{B} \rightarrow \infty, \frac{N_{A}}{N_B}\rightarrow r_{AB}}{\rightarrow} -\mu \langle \hspace{-0.8mm}\langle g^{AB},\nabla_{x_1} \Phi^{AB}(x_1,x_2) \cdot \bigg(\frac{1}{f^A(x_1)}\int \big(\nabla_{x_1} \Psi g^{AB}\big) (x_1,x_2) dx_2 \bigg) \rangle \hspace{-0.8mm}\rangle,
\end{align*}
while
\eq{
{-\frac{\mu}{K_{AB}} \sum_{k_3=1}^{K_{AB}} \nabla_{x_1} \Phi^{AB}(X^A_{i(k_3)},X^B_{\ell(k_3)})
 \nabla_{x_1} \Psi(X^A_{i(k_3)},X^B_{\ell(k_3)})}\\
\underset{N_A, K_{AB} \rightarrow \infty}{\rightarrow} 
 {-\mu\langle  \hspace{-0.8mm}\langle g^{AB},\Grad_{x^1}\Phi^{AB}(x_1,x_2)\Grad_{x^1}\Psi(x_1,x_2)  \rangle  \hspace{-0.8mm}\rangle}.
}
\noindent Altogether, performing the same computations for $E_2$ in \eqref{dg1a}, we obtain:
\begin{align*}
\frac{d}{dt} \langle  \hspace{-0.8mm} \langle  g^{AB},\Psi  \rangle  \hspace{-0.8mm} \rangle =& - 2\mu \langle \hspace{-0.8mm}\langle g^{AA},\nabla_{x_1} \Phi^{AA}(x_1,x_2) \cdot \bigg(\frac{1}{f^A(x_1)}\int \big( \nabla_{x_1} \Psi g^{AB}\big) (x_1,x_2) dx_2 \bigg) \rangle \hspace{-0.8mm}\rangle \\
&-\mu \langle \hspace{-0.8mm}\langle g^{AB},\nabla_{x_1} \Phi^{AB}(x_1,x_2) \cdot \bigg(\frac{1}{f^A(x_1)}\int \big( \nabla_{x_1} \Psi g^{AB}\big) (x_1,x_2) dx_2 \bigg) \rangle \hspace{-0.8mm}\rangle\\
&- 2\mu \langle \hspace{-0.8mm}\langle g^{BB},\nabla_{x_1} \Phi^{BB}(x_1,x_2) \cdot \bigg(\frac{1}{f^B(x_1)}\int \big( {\nabla_{x_2}} \Psi g^{AB}\big) (x_2,x_1) dx_2 \bigg) \rangle \hspace{-0.8mm}\rangle \\
&-\mu \langle \hspace{-0.8mm}\langle g^{BA},\nabla_{x_1} \Phi^{BA}(x_1,x_2) \cdot \bigg(\frac{1}{f^B(x_1)}\int  \big(\nabla_{x_2} \Psi g^{AB}\big) (x_2,x_1) dx_2 \bigg) \rangle \hspace{-0.8mm}\rangle\\
& {-\mu\langle  \hspace{-0.8mm}\langle g^{AB},\Grad_{x_1}\Phi^{AB}(x_1,x_2)\cdot\Grad_{x_1}\Psi(x_1,x_2)  \rangle  \hspace{-0.8mm}\rangle}\\
& {-\mu\langle  \hspace{-0.8mm}\langle g^{AB},\Grad_{x_2}\Phi^{BA}(x_1,x_2)\cdot\Grad_{x_2}\Psi(x_1,x_2)  \rangle  \hspace{-0.8mm}\rangle}.
\end{align*}
\noindent By carefully performing integration by parts and change of order of integrals, we can obtain:
\begin{align*}
\frac{d}{dt} \langle  \hspace{-0.8mm} \langle  g^{AB},\Psi  \rangle  \hspace{-0.8mm} \rangle =&2\mu \langle \hspace{-0.8mm}\langle \nabla_{x_1} \cdot \bigg( \frac{g^{AB}(x_1,x_2)}{f^A(x_1)} \int \big(g^{AA} \nabla_{x_1} \Phi^{AA} \big)(x_1,x_2) dx_2 \bigg), \Psi(x_1,x_2) \rangle \hspace{-0.8mm}\rangle\\
&+\mu \langle \hspace{-0.8mm}\langle \nabla_{x_1} \cdot \bigg( \frac{g^{AB}(x_1,x_2)}{f^A(x_1)} \int \big(g^{AB} \nabla_{x_1} \Phi^{AB} \big)(x_1,x_2) dx_2 \bigg), \Psi(x_1,x_2) \rangle \hspace{-0.8mm}\rangle\\
&+2\mu \langle \hspace{-0.8mm}\langle \nabla_{x_2} \cdot \bigg( \frac{g^{AB}(x_1,x_2)}{f^B(x_2)} \int \big(g^{BB} \nabla_{x_1} \Phi^{BB}\big)(x_2,x_1) dx_1 \bigg), \Psi(x_1,x_2) \rangle \hspace{-0.8mm}\rangle \\
&+\mu \langle \hspace{-0.8mm}\langle \nabla_{x_2} \cdot \bigg( \frac{g^{AB}(x_1,x_2)}{f^B(x_2)} \int \big(g^{BA} \nabla_{x_1} \Phi^{BA}\big)(x_2,x_1) dx_1 \bigg), \Psi(x_1,x_2) \rangle \hspace{-0.8mm}\rangle\\
  &{+\mu \langle  \hspace{-0.8mm} \langle \nabla_{x_1} \cdot \bigg( g^{AB}(x_1,x_2) \nabla_{x_1} \Phi^{AB}(x_1,x_2)  \bigg), \Psi(x_1,x_2)  \rangle  \hspace{-0.8mm} \rangle}\\
  &{+\mu \langle  \hspace{-0.8mm} \langle \nabla_{x_2} \cdot \bigg( g^{AB}(x_1,x_2) \nabla_{x_2} \Phi^{BA}(x_1,x_2)  \bigg), \Psi(x_1,x_2)  \rangle  \hspace{-0.8mm} \rangle}.
\end{align*}
\noindent Finally, note that in the case of interspecies links, the average number of new links created in the interval $[t,t+\Delta t[$ is equal to 
$$
\nu^{AB}_{c,N,\ep} \bigg(N_AN_B h^{AB}(x_1, x_2, t)\chi(|x_1-x_2|\leq R) - N_{A} g^{AB}(x_1,x_2,t) \bigg) dx_1 dx_2 \Delta t.
$$
\noindent Assuming for $\nu_{c,N,\ep}^{AB}$ 
the scaling \eqref{eq:nuscalings}
as the number of particles tends to infinity and adding the noise, we obtain the final equation for $g^{AB}$:
\begin{align}
\partial_t g^{AB}&(x_1,x_2,t) = \big(D^A\Delta_{x_1} g^{AB}(x_1,x_2,t) + D^B\Delta_{x_2} g^{AB}(x_1,x_2,t) \big) \nonumber\\
&+ 2\mu \nabla_{x_1} \cdot \bigg(\frac{g^{AB}(x_1,x_2,t)}{f^A(x_1,t)}F^{AA}[g^{AA}](x_1,t)\bigg)+ \mu \nabla_{x_1} \cdot \bigg(\frac{g^{AB}(x_1,x_2,t)}{f^A(x_1,t)}F^{AB}[g^{AB}](x_1,t)\bigg)\nonumber\\
&+ 2\mu \nabla_{x_2} \cdot \bigg(\frac{g^{AB}(x_1,x_2,t)}{f^B(x_2,t)}F^{BB}[g^{BB}](x_2,t)\bigg)+ \mu \nabla_{x_2} \cdot \bigg(\frac{g^{AB}(x_1,x_2,t)}{f^B(x_2,t)}F^{BA}[g^{BA}](x_2,t)\bigg)\nonumber \\
&{+\mu \nabla_{x_1} \cdot \bigg( g^{AB}(x_1,x_2) \nabla_{x_1} \Phi^{AB}(x_1,x_2)  \bigg)+\mu \nabla_{x_2} \cdot \bigg( g^{AB}(x_1,x_2) \nabla_{x_2} \Phi^{BA}(x_1,x_2)  \bigg)}\nonumber\\
&+ \nu^{AB}_{c,\ep}  h^{AB}(x_1, x_2, t)\chi(|x_1-x_2|\leq R)
- \nu^{AB}_{d,\ep}  g^{AB}(x_1,x_2,t),\nonumber
\end{align}
\noindent where $h^{AB}(x_1, x_2, t) = \lim_{N_A,N_B \rightarrow \infty} h^{N_{AB}}(x_1, x_2, t)$ with $h^{N_{AB}}(x_1, x_2, t)$ defined by Eq. \eqref{TPDFAB}. %Note that $g^{AB}(x_1,x_2) = g^{BA}(x_2,x_1)$.
\end{pf}
\section{Scaling of the kinetic model}
\label{AppendixB}

\subsection{Dimensionless Equations}
In order to express the problem in dimensionless variables, we denote by $t_0$  the unit of time and $x_0$, $f^S_0=\frac{1}{x_0^2}$, $g^{ST}_0 = \frac{1}{x_0^4}$, $g^{ST}_0 =  \frac{1}{x_0^4}$ the units of space and distribution functions, where $S$ and $T$ can be either $A$ or $B$ and refer to the particle type. The scaling of $f^S(x,\theta)$, $g^{ST}(x_1,x_2)$ and $h^{ST}(x_1,x_2)$ comes from the fact that they are probability distribution functions on a 2D domain. The following dimensionless variables are defined: 
\begin{equation*}
\bar{t} = \frac{t}{t_0}, \; \bar{x} = \frac{x}{x_0},\; \bar{f^S}=\frac{f^S}{f_0} = f^S x_0^2, \;  \bar{g}^{ST} = \frac{g^{ST}}{g_0} = g^{ST} x_0^4, \;\bar{h}^{ST} = \frac{h^{ST}}{h_0} = h^{ST} x_0^4.
\end{equation*}
\noindent and the following dimensionless parameters are introduced: 
\begin{equation*}
\mu' = \frac{\mu}{t_0}, \; \nu'^{ST}_{c,\infty}  = t_0 \nu^{ST}_{c,\infty} ,\;  \nu'^{ST}_{d,\infty}  = t_0 \nu^{ST}_{d,\infty} ,\; R'=\frac{R}{x_0}, \; D'=\frac{{D t_0}}{x_0^2}, \Phi'^{ST}=\frac{\Phi^{ST} t_0^2}{x_0^2},
\end{equation*}
\noindent where we assumed that the potential scales as the potential energy $\frac{x_0^2}{t_0^2}$. We first have:
\eqh{
\partial_t f^S(x,t)& = \frac{1}{t_0 x_0^2} \partial_{\bar{t}}\bar{f}^S(\bar{x},\bar{t}), \; \partial_t g^{ST}(x_1,x_2,t) \\
&= \frac{1}{t_0 x_0^4} \partial_{\bar{t}}\bar{g}^{ST}(\bar{x}_1,\bar{x}_2,\bar{t}), \; \partial_t h^{ST}(x_1,x_2,t) = \frac{1}{t_0 x_0^4} \partial_{\bar{t}}\bar{h}^{ST}(\bar{x}_1,\bar{x}_2,\bar{t})
}
\noindent and:
\begin{equation*}
\begin{split}
\nabla_{x} \cdot \big(F^{ST}[g^{ST}](x)\big) = \frac{1}{x_0} \nabla_{x'} \cdot \big(\int g^{ST}\nabla_{x_1} \Phi^{ST} dx_2 \big) &= \frac{1}{x_0} \nabla_{x'} \cdot \big(\int \frac{\bar{g}^{ST}}{x_0^4} \frac{x_0}{t_0^2} \nabla_{x'_1} \Phi'^{ST} x_0^2 d\bar{x}_2 \big)\\
&= \frac{1}{x_0^2t_0^2} \nabla_{x'} \cdot \big(F'^{ST}[\bar{g}^{ST}](x)\big).
\end{split}
\end{equation*}
\noindent In this new set of variables, choosing $S=A, T=A$ (the same scaling apply for the other equations) and omitting the primes and bar for clarity, Eqs.~\eqref{Eqfkinet}-\eqref{EqgAAkinet} become :
\begin{align*}
 \partial_t f^A(x,t) = 2\mu \nabla_x \cdot F^{AA}[g^{AA}](x,t) + \mu \nabla_x \cdot F^{AB}[g^{AB}](x,t) + D^A\Delta f^A,
 \end{align*}
\noindent and 
 \begin{align*}
\partial_t g^{AA}&(x_1,x_2,t) = D^A \big(\Delta_{x_1} g^{AA}(x_1,x_2,t) +\Delta_{x_2} g^{AA}(x_1,x_2,t) \big) \nonumber\\
&+ 2\mu \nabla_{x_1} \cdot \bigg(\frac{g^{AA}(x_1,x_2,t)}{f^A(x_1,t)}F^{AA}[g^{AA}](x_1,t)\bigg)+ \mu \nabla_{x_1} \cdot \bigg(\frac{g^{AA}(x_1,x_2,t)}{f^A(x_1,t)}F^{AB}[g^{AB}](x_1,t)\bigg)\\
&+ 2\mu \nabla_{x_2} \cdot \bigg(\frac{g^{AA}(x_1,x_2,t)}{f^A(x_2,t)}F^{AA}[g^{AA}](x_2,t)\bigg)+ \mu \nabla_{x_2} \cdot \bigg(\frac{g^{AA}(x_1,x_2,t)}{f^A(x_2,t)}F^{AB}[g^{AB}](x_2,t)\bigg) \nonumber \\
&{+\mu \nabla_{x_1} \cdot \bigg( g^{AA}(x_1,x_2,t) \nabla_{x_1} \Phi^{AA}(x_1,x_2)  \bigg)+\mu \nabla_{x_2} \cdot \bigg( g^{AA}(x_1,x_2,t) \nabla_{x_2} \Phi^{AA}(x_1,x_2)  \bigg)}\nonumber\\
&+ \frac{\nu^{AA}_{c,\ep}  }{2} h^{AA}(x_1, x_2, t)\chi(|x_1-x_2|\leq R)
- \nu^{AA}_{d,\ep}  g^{AA}(x_1,x_2,t).\nonumber
\end{align*}
\noindent Finally, we choose the space and time scales $x_0$, $t_0$ such that $\mu = 1$.

\subsection{Scaled equations}\label{scaling}
In order to describe the system at a macroscopic scale, a small parameter $\varepsilon \ll 1$ is introduced and the space and time units are set to $\tilde{x_0} = \varepsilon^{-1/2} x_0$, $\tilde{t_0} = \varepsilon^{-1} t_0$. The variables $x$, $t$, $R$ and unknowns $f$ and $g$ are then correspondingly changed to $\tilde{x} = \sqrt{\varepsilon} x$, $\tilde{t} = \varepsilon t$, $\tilde{R} = \sqrt{\varepsilon} R$. Therefore,  $\tilde{f}^S(\bar{x}) = \varepsilon^{-1} f^S(x)$, $\tilde{g}^{ST}(\tilde{x}_1,\tilde{x}_2,\tilde{t}) = \varepsilon^{-2} g(x_1,x_2,t)$ and $\tilde{h}^{ST}(\tilde{x}_1,\tilde{x}_2,\tilde{t}) = \varepsilon^{-2} h(x_1,x_2,t)$. The diffusion constant is supposed to be of order 1, $D^S = \tilde{D}^S$, and we suppose that the interaction potentials scale as $\Phi^{ST}(x_1,x_2) = \tilde{\Phi}^{ST}(\tilde{x}_1,\tilde{x}_2)$. Then, 
\eqh{
\nabla_{x} \cdot \big(F^{ST}[g^{ST}](x)\big)& = \frac{1}{\sqrt{\varepsilon}} \nabla_{\tilde{x}} \cdot \big(\int \varepsilon^2 \tilde{g}^{ST} \frac{1}{\sqrt{\varepsilon}} \nabla_{\tilde{x}_1} \tilde{\Phi}^{ST} \varepsilon d\tilde{x}_2\big)\\
 &= \varepsilon^2 \nabla_{\tilde{x}} \cdot \big(\int \tilde{g}^{ST}  \nabla_{\tilde{x}_1} \tilde{\Phi}^{ST}  d\tilde{x}_2\big) = \varepsilon^2 \nabla_{\tilde{x}} \cdot \big( \tilde{F}^{ST}[\tilde{g}^{ST}](\tilde{x})\big),
}
\noindent and with {$\mu = 1$}, we obtain the same equation for $\tilde{f}^S$ (for $S$=A for instance):
\begin{align*}
 \partial_{\tilde{t}} \tilde{f}^A(\tilde{x},\tilde{t}) = 2 \nabla_{\tilde{x}} \cdot \tilde{F}^{AA}[\tilde{g}^{AA}](\tilde{x},\tilde{t}) + \nabla_{\tilde{x}} \cdot \tilde{F}^{AB}[\tilde{g}^{AB}](\tilde{x},\tilde{t}) + D^A\Delta \tilde{f}^A.
 \end{align*}
\noindent In order to simplify the analysis of the system, the process of linking/unlinking is supposed to occur at a very fast time scale: this is the meaning of the $\ep$-scaling of the rates in \eqref{eq:nuscalings}.
 As $\chi_{|x_1 - x_2|\leq R} = \chi_{|\tilde{x}_1 - \tilde{x}_2|\leq \tilde{R}}$, we have (for $S,T=$(A,A)):
 \begin{align}
\varepsilon^3 \partial_{\tilde{t}} \tilde{g}^{AA}& =  \varepsilon^3 D^A \Delta_{\tilde{x}_1} \tilde{g}^{AA} + \varepsilon^3 D^A \Delta_{\tilde{x}_2} \tilde{g}^{AA} \nonumber\\
&+ 2 \sqrt{\varepsilon}\nabla_{\tilde{x}_1} \cdot \bigg( \frac{\varepsilon^2\tilde{g}^{AA}}{\varepsilon\tilde{f}^A}\varepsilon^{3/2}\tilde{F}^{AA}[\tilde{g}^{AA}]\bigg)+  \sqrt{\varepsilon}\nabla_{\tilde{x}_1} \cdot \bigg( \frac{\varepsilon^2\tilde{g}^{AA}}{\varepsilon\tilde{f}^A}\varepsilon^{3/2}\tilde{F}^{AB}[\tilde{g}^{AB}]\bigg) \nonumber\\
&+ 2 \sqrt{\varepsilon}\nabla_{\tilde{x}_2} \cdot \bigg( \frac{\varepsilon^2\tilde{g}^{AA}}{\varepsilon\tilde{f}^A}\varepsilon^{3/2}\tilde{F}^{AA}[\tilde{g}^{AA}]\bigg)+  \sqrt{\varepsilon}\nabla_{\tilde{x}_1} \cdot \bigg( \frac{\varepsilon^2\tilde{g}^{AA}}{\varepsilon\tilde{f}^A}\varepsilon^{3/2}\tilde{F}^{AB}[\tilde{g}^{AB}]\bigg) \nonumber \\
&+ \frac{{\nu}^{AA}_c }{2\varepsilon^2} \varepsilon^2 \tilde{h}^{AA}\chi(|\tilde{x}_1-\tilde{x}_2|\leq \tilde{R})
- \frac{{\nu}^{AA}_d}{\varepsilon^2} \varepsilon^2 \tilde{g}^{AA}\nonumber\\
&{+\sqrt{\ep}\nabla_{\tilde x_1} \cdot \bigg( \ep^2\tilde g^{AA}\sqrt{\ep} \nabla_{\tilde x_1} \tilde \Phi^{AA} \bigg)+\sqrt{\ep}\nabla_{\tilde x_2} \cdot \bigg(\ep^2 \tilde g^{AA} \sqrt{\ep}\nabla_{\tilde x_2}\tilde \Phi^{AA}  \bigg)}\nonumber\\
& =  \varepsilon^3 \bigg[D \Delta_{\tilde{x}_1} \tilde{g}^{AA} +  D \Delta_{\tilde{x}_2} \tilde{g}^{AA} \nonumber\\
&+ 2 \nabla_{\tilde{x}_1} \cdot \bigg( \frac{\tilde{g}^{AA}}{\tilde{f}^A}\tilde{F}^{AA}[\tilde{g}^{AA}]\bigg)+ \nabla_{\tilde{x}_1} \cdot \bigg( \frac{\tilde{g}^{AA}}{\tilde{f}^A}\tilde{F}^{AB}[\tilde{g}^{AB}]\bigg) \label{gkinetrescaled}\\
&+ 2 \nabla_{\tilde{x}_2} \cdot \bigg(\frac{\tilde{g}^{AA}}{\tilde{f}^A}\tilde{F}^{AA}[\tilde{g}^{AA}]\bigg)+ \nabla_{\tilde{x}_1} \cdot \bigg( \frac{\tilde{g}^{AA}}{\tilde{f}^A}\tilde{F}^{AB}[\tilde{g}^{AB}]\bigg) \bigg], \nonumber \\
&{+\ep^3\nabla_{\tilde x_1} \cdot \bigg( \tilde g^{AA} \nabla_{\tilde x_1} \tilde \Phi^{AA} \bigg)+\ep^3\nabla_{\tilde x_2} \cdot \bigg( \tilde g^{AA} \nabla_{\tilde x_2}\tilde \Phi^{AA}  \bigg)}\nonumber\\
&+ \frac{{\nu}^{AA}_c }{2}  \tilde{h}^{AA}\chi(|\tilde{x}_1-\tilde{x}_2|\leq \tilde{R})
- {\nu}^{AA}_d \tilde{g}^{AA}. \nonumber
\end{align}
\noindent Now, we aim to pass to the limit $\varepsilon \rightarrow 0$. Denoting $f_{\varepsilon}^S = \tilde{f}^S$, $g_{\varepsilon}^{ST} = \tilde{g}^{ST}$ and $h_{\varepsilon}^{ST} = \tilde{h}^{ST}$, we
want to derive the same system of macroscopic equations as in Proposition \ref{prop1} of section \ref{sec:macro}.

\begin{pf}
From Eq.\eqref{gkinetrescaled} generalized to $S,T$ and using the assumption $h^{ST}_{\varepsilon} (x_1,x_2,t) = f_{\varepsilon}^S(x_1,t)f_{\varepsilon}^T(x_2,t)$, and dropping the tildas, we have:
\begin{align*}
\frac{\nu^{SS}_c}{2} f^{S}(x_1,t)f^{S}(x_2,t)\chi_{|x_1-x_2|\leq R} - \nu^{SS}_d g^{SS}(x_1,x_2,t) = O(\varepsilon^3),\\
\nu^{ST}_c f^{S}(x_1,t)f^{T}(x_2,t)\chi_{|x_1-x_2|\leq R} - \nu^{ST}_d g^{ST}(x_1,x_2,t) = O(\varepsilon^3).
\end{align*}
\noindent Therefore in the limit $\varepsilon \rightarrow 0$, we have that
\begin{align*}
g^{SS}(x_1,x_2,t) = \frac{\nu^{SS}_c}{2\nu^{SS}_d} f^{S}(x_1,t)f^{S}(x_2,t)\chi_{|x_1-x_2|\leq R}\\
g^{ST}(x_1,x_2,t) = \frac{\nu^{ST}_c}{\nu^{ST}_d} f^{S}(x_1,t)f^{T}(x_2,t)\chi_{|x_1-x_2|\leq R}
\end{align*}
\noindent for all $S,T = A,B$. Plugging the expression of $g^{ST}$ into the equation for $f^S$, we have
\eq{\label{fg_app0}
\partial_t f^S(x_1,t) = 2 \nabla_x \cdot F^{SS}[g^{SS}](x,t) + \nabla_x F^{ST}[g^{ST}](x,t) + D \Delta f^S,
}
\noindent where
\eq{\label{fg_app1}
F^{SS}[g^{SS}](x,t) &= \frac{\nu^{SS}_c}{2\nu^{SS}_d}f^S(x,t)\int f^S(y,t) \nabla_{x} {\Phi}^{SS}(x,y) \chi_{|x-y|\leq R}dy \\
F^{ST}[g^{ST}](x,t) &= \frac{\nu^{ST}_c}{\nu^{ST}_d}f^S(x,t)\int f^T(y,t) \nabla_{x} {\Phi}^{ST}(x,y) \chi_{|x-y|\leq R}dy.
}
\noindent Therefore, if the potentials $\Phi^{ST}(x_1,x_2)=U^{ST} (|x_1-x_2|)$, we can write 
\eq{\label{fg_app2}
F^{SS}[g^{SS}](x,t) &= \frac{1}{2}f^S(x,t)\int f^S(y,t) \nabla_x \tilde{\Phi}^{SS}(|x-y|)dy \\
F^{ST}[g^{ST}](x,t) &= f^S(x,t)\int f^T(y,t) \nabla_{x} \tilde{\Phi}^{ST}(|x-y|)dy,
}
\noindent for some potentials $\tilde{\Phi}^{ST}$ such that: 
\eq{\label{fg_app3}
\nabla_i \tilde{\Phi}^{ST}(x) = \frac{\nu^{ST}_c}{\nu^{ST}_d}\big(U^{ST}\big)'(|x|)\chi_{|x|\leq R}\vec{e}_i, \quad i=1,2.
}
Recall that due to the scaling of linking frequencies with $N_A, N_B$ we obtain in the limit that 
 $\frac{\nu^{BA}_c}{\nu^{BA}_d} =r_{AB} \frac{\nu^{AB}_c}{\nu^{AB}_d}$. With this observation the proof of Proposition \ref{prop1} is complete.
\end{pf}

\section{Numerical data visualisation}
\label{AppendixC}
\subsection{Macroscopic model}
For the numerical simulations of the  macroscopic model we consider a periodic square $[-7.5, \; 7.5] \times [ -7.5, \; 7.5]$ discretized with a space step $\Delta x = 0.3$ and the time step $\Delta t = 5\cdot 10^{-3}$. As in the microscopic case, we fix the radii of all types of interactions to $R=1$ and consider the same diffusion coefficient for the species $D=10^{-4}$. We use the numerical scheme introduced in \cite{CaChHu} that was developed in our recent work \cite{Barre2017} to study aggregation-diffusion equation for single species. 

We consider all four cases of  interaction intensities $\kappa^{ST}$ as in the microscopic case, as summarized in Table \ref{Table1}. The initial densities for species $A$ and species $B$ are random perturbations of constant functions, The constants are chosen such that the total mass is equal to 1, i.e.
$$f^S_0(x,y)=\frac{1+0.01\cdot X_S(x,y)}{\int_{-7.5}^{7.5}\int_{-7.5}^{7.5}(1+0.01\cdot X_S(x,y))dx\,dy}, \quad S=A,B,$$
where $X_S(x,y)$ is uniformly distributed random variable between 0 and 1.

To visualise the macroscopic simulations we plot green or red balls in the regions where the densities of green or red cells, respectively, dominate. The balls are of radius $0.3$ with a center at the center of the corresponding pixel perturbed by uniformly distributed random number from interval $[-0.2,0.2]$  in $x$ and $y$ direction. We determine that the concentration of green cells dominates over the red ones if the difference between their densities is larger than parameter $tres_1$. It is equal to $5.0e-4$ times the maximum of the densities of green or red cells at the final time of simulation.  If the difference between the densities of  cells  is less than $tres_1$ but more than $tres_2=tres_1-5.0e-6$ we plot randomly red or green ball. The "empty" black regions on the figures (see for example Figure \ref{fig:step} and last row of Figure \ref{fig:kBB1}) correspond to the case when the difference between the densities of the cells is less than $tres_2$.
\subsection{Image processing}\label{AppendixC2}
In order to compare quantitatively the microscopic and macroscopic simulations, we use image processing tools to define quantifiers of the structures observed in the simulation images. As we will restrict ourselves to regimes where the $B$-family clusterize (represented in green), the developed tools aim to describe type B clusters but can easily be adapted to detect red clusters. Given a RGB image transformed into a binary image via thresholding of the intensity, we use morphological operators to erase the isolated pixels and dilate the image using a binary gradient mask. Interior gaps are then filled and the new image is smoothed via morphological operations. The boundaries of the isolated clusters are then detected using the Moore-Neighbor tracing algorithm modified by Jacob's stopping criteria, implemented in the intrinsic matlab function bwconncomp that we couple with a function that appropriately converts the output of bwconncomp to take into account periodic boundary conditions. Once the clusters are appropriately separated and borders are detected, we finally use the matlab intrinsic function regionprops to measure each 8-connected object (region) of the image. This image processing enables us to compute the number $N_g$ of green clusters in the image. The elongation of each cluster is given by the eccentricity of the ellipse that has the same second-moments as the region (cluster). Finally, we define the overall overlapping amount $Q$ as:
\begin{equation}\label{Q}
Q = \frac{\sum_{i=1}^{N_P} P^i_g P^i_r + (1-max(P^i_g,P^i_r))}{N_P},    
\end{equation}
\noindent where $P^i_g$ (resp. $P^i_r$) is equal to 1 if pixel $i$ has a non-zero green component (resp. red) and $N_P$ is the total number of pixels in the image. Thus defined, $Q\approx 1$ when the two families are perfectly mixed (corresponding to all pixels having both red and green components), and $Q=0$ correspond to completely separated phases (where each pixel is either green or red). Therefore, parameter $Q$ enables to distinguish between homogeneous and segregated states. However, note that $Q$ does not give any information on the form of the clusters.

\bibliographystyle{ieeetr}

\end{document}